\documentclass[onefignum,onetabnum]{siamart171218}

\usepackage{amsmath,amsfonts,amssymb,epsfig,graphicx,latexsym,float,color}

\usepackage{enumerate,relsize,mathrsfs}   
\usepackage{sidecap}
\usepackage{multirow}
\usepackage{placeins}

\newsiamthm{remark}{Remark}
\newsiamthm{myalgorithm}{Algorithm}

\renewcommand{\theequation}{\arabic{section}.\arabic{equation}}

\DeclareMathOperator{\myim}     {image}

\newcommand{\tensorbj}[1]{{\mathbf #1}}
\newcommand{\vecbj}[1]{{\mathbf #1}} 
\newcommand{\tensorcj}[1]{{\mathcal #1}} 
\newcommand{\veccj}[1]{{\mathfrak #1}}

\newenvironment{mymatrix}{\begin{bmatrix}}{\end{bmatrix}}

\newcommand{\Transfer}[1]{\mathscr{#1}}

\DeclareMathOperator{\moma}   {\tensorbj{m}}
\DeclareMathOperator{\momb}   {\boldsymbol{\mu}}
\DeclareMathOperator{\momc}   {\boldsymbol{\eta}}

\newcommand{\circledmy}[1]{\normalfont{\textcircled{\relsize{-0.5}{#1}}}} 
\newcommand{\supertens}[1]{^{\normalfont{\scriptsize\textcircled{{\relsize{-1}#1}}}}} 
\newcommand{\spacevar}{\xi}
\newcommand{\statemap}{\tensorbj{S}}
\newcommand{\statemapAut}{\mathcal{S}}
\newcommand{\siggen}{\tensorbj{T}}

\newcommand{\sigstate}{\vecbj{z}}
\newcommand{\sigstatescal}{{z}}

\def \tikzcolor {blue}
\usepackage{tikz}
\usetikzlibrary{calc,positioning,shapes,arrows,decorations.shapes} 
\usetikzlibrary{decorations.text}

\tikzstyle{arr}=[-latex, black, line width=0.5pt]
\tikzstyle{doublearr}=[latex-latex, black, line width=0.5pt]
\tikzstyle{input}=[font=\small\sffamily\bfseries]
\tikzstyle{rect}=[
   rectangle, draw=black, font=\small\sffamily\bfseries, inner sep=9pt]
   \tikzstyle{rectblue}=[
   rectangle, draw=black, font=\small\sffamily\bfseries, inner sep=9pt]
\tikzstyle{rect2}=[
   rectangle, draw=black, font=\small\sffamily\bfseries, inner sep=4pt]
\tikzstyle{rect3}=[
   rectangle, draw=black, font=\small\sffamily\bfseries, inner sep=3pt]

\headers{Input-tailored system-theoretic MOR}{B.~Liljegren-Sailer and N. Marheineke}

\title{Input-tailored system-theoretic model order reduction for quadratic-bilinear systems
\thanks{Submitted to the editors XXX \funding{The financial support of the German Federal Ministry for Economic Affairs and Energy, Project MathEnergy, is acknowledged.}}}

\author{Bj\"orn Liljegren-Sailer\thanks{Universit\"at Trier, FB IV - Mathematik, Lehrstuhl Modellierung und Numerik, D-54286 Trier, Germany
  (\email{Corresponding author: bjoern.sailer@uni-trier.de}).}
\and Nicole Marheineke\footnotemark[2]}

\usepackage{amsopn}

\begin{document}

\maketitle

\begin{abstract}
In this paper we suggest a moment matching method for qua\-dra\-tic-bilinear dynamical systems. Most system-theoretic reduction methods for nonlinear systems rely on multivariate frequency representations. Our approach instead uses univariate frequency representations tailored towards user-pre-defined families of inputs. Then moment matching corresponds to a one-dimensional interpolation problem, not to multi-dimensional interpolation as for the multivariate approaches, i.e., it also involves fewer interpolation frequencies to be chosen. Compared to former contributions towards nonlinear model reduction with univariate frequency representations, our approach shows profound differences: Our derivation is more rigorous and general and reveals additional tensor-structured approximation conditions, which should be incorporated. Moreover, the proposed implementation exploits the inherent low-rank tensor structure, which enhances its efficiency.
In addition, our approach allows for the incorporation of more general input relations  in the state equations -- not only affine-linear ones as in existing system-theoretic methods -- in an elegant way. As a byproduct of the latter, also a novel modification for the multivariate methods falls off, which is able to handle more general input-relations.
\end{abstract}

\noindent
\textsc{AMS-Classification:} 93Axx, 37N30, 41Axx
 
\section{Introduction}

\subsection*{The approach in a nutshell} \label{subsec:prob-formul-assoc}
In this paper we introduce a new system-theoretic model order reduction method for quadratic-bilinear dynamical systems of the form
\begin{align*} 
	\tensorbj{E} \dot{\vecbj{x}} &= \tensorbj{A} \vecbj{x} + \tensorbj{G} (\vecbj{x} \otimes \vecbj{x}) + \tensorbj{D} (\vecbj{x} \otimes \vecbj{u}) + \tensorbj{B} \vecbj{u}, \quad t \geq 0 \\
	\vecbj{y} &= \tensorbj{C} \vecbj{x}, \qquad  \vecbj{x}(0)= \vecbj{x}_0 \in \mathbb{R}^N, \qquad \vecbj{u}(t) \in \mathbb{R}^p, \quad t \geq 0
\end{align*}
with nonsingular matrix $\tensorbj{E}$ and Kronecker-tensor product $\otimes$, i.e., $\vecbj{x}\otimes \vecbj{x} \in \mathbb{R}^{N^2}$ and  $\vecbj{x} \otimes \vecbj{u} \in \mathbb{R}^{Np}$. The system characterizes a map $\vecbj{u}\mapsto \vecbj{y}$ from typically low-dimensional input $\vecbj{u}$ to low-dimensional output $\vecbj{y}$ via a high-dimensional state $\vecbj{x}$. For a cheaper-to-evaluate reduced model, we seek for an appropriate basis matrix $\tensorbj{V} \in \mathbb{R}^{N,n}$, $n\ll N$, and define the reduced model as
\begin{align*} 
	\tensorbj{E}_r \dot{\vecbj{x}}_r &= \tensorbj{A}_r \vecbj{x}_r + \tensorbj{G}_r ( \vecbj{x}_r  \otimes \vecbj{x}_r) + \tensorbj{D}_r (\vecbj{x}_r \otimes \vecbj{u}) + \tensorbj{B}_r \vecbj{u} 	\\
	\tilde{\vecbj{y}} &= \tensorbj{C}_r \vecbj{x}_r, \qquad \vecbj{x}_r(0)= \tensorbj{V}^T \vecbj{x}_0 \in \mathbb{R}^n 
\end{align*}
with $\tensorbj{E}_r = \tensorbj{V}^T \tensorbj{E} \tensorbj{V}$, $\tensorbj{A}_r = \tensorbj{V}^T \tensorbj{A} \tensorbj{V}$, $\tensorbj{G}_r = \tensorbj{V}^T \tensorbj{G} (\tensorbj{V} \otimes \tensorbj{V})$, $\tensorbj{D}_r = \tensorbj{V}^T \tensorbj{D} (\tensorbj{V} \otimes \tensorbj{I}_p)$, $\tensorbj{B}_r = \tensorbj{V}^T \tensorbj{B}$, $\tensorbj{C}_r = \tensorbj{C} \tensorbj{V}$ and unit matrix $\tensorbj{I}_p$ of dimension $p$. Many system-theoretic methods for linear systems are based on the frequency representation of the input-output map, which is a univariate algebraic mapping, called transfer function. For moment matching, the reduction basis $\tensorbj{V}$ is developed such that the transfer function of the reduced model fulfills certain interpolation conditions. In the nonlinear case, the input-output map generally does not have a univariate frequency representation. Relaxations of the linear notions are needed to generalize it to the nonlinear case, see, e.g., recent multi-moment matching methods for multivariate frequency representations \cite{phd:Gu-mor-nonlinear}, \cite{art:two-sided-hermm}, \cite{art:quadratic-bilinear-regular-krylov}, \cite{inproc:conf/eucc/GoyalAB15}, \cite{inproc:morBenB12}. In our approach we pursue an other idea by using the following three relaxation steps:

\begin{figure}\label{fig:io-sys-siggendriv}
\begin{tabular}{ccc}
\begin{minipage}{0.35\textwidth}
\begin{tikzpicture}[scale=0.7, every node/.style={scale=0.8}]
\node 			(u) at (0.0, 0) {$\vecbj{u}$};
\node 		 	(y) at (5., 0) {$\vecbj{y}$};

\node [rectangle, draw=black,align=center]	(iosys) at (2.4, 0) {Input-output\\ system};
\node 			(Ghostnode) at (0, 1.3) {};
\draw [ ->] (0.45, 0) -- (1.0,0);
\draw [black, ->] (3.8, 0) -- (4.5,0);
\end{tikzpicture} 
\end{minipage}
\hspace{-1.4cm}
&
\begin{minipage}{0.2\textwidth}
\begin{tikzpicture}[scale=0.5, every node/.style={scale=0.8}]
\node 			(VV1) at (-0.7, 0.3) {};
\node 			(VV2) at (2.2, 0.3) {};
\node 			(Ghostnode) at (0, -1) {};
\draw [thick,dashed,bend angle=25, bend left, ->]  (VV1) to (VV2); 
\node [rotate=0,align = center]			(VV12) at (0.4,1.3) {\small{\textit{Drive by signal}}\\ \small{\textit{generator}}};
\end{tikzpicture} 
\end{minipage}
\hspace{-1cm}
&
\begin{minipage}{0.5\textwidth}
\begin{tikzpicture}[scale=0.8, every node/.style={scale=0.8}]
\node [rectangle, draw=black, align=center]	(SG) at (-0.7, 0) {Signal\\ generator};
\node 		 	(y) at (5.6, 0) {$\vecbj{y}$};

\node [rectangle, draw=black, align=center]	(iosys) at (2.7, 0) {Input-output\\ system};
\node [align=center]	 	(AutSys) at (1.7, 1.3) {Autonomous output\\ system};

\draw [ ->] (0.45, 0) -- (1.0,0);
\draw [black, ->] (4.3, 0) -- (5.,0);
\draw (-1.8,-0.8) rectangle (4.1,0.7);
 \\
\end{tikzpicture} 
\end{minipage}
\end{tabular}
\caption{\textit{Left:} Classical input-output modeling. \textit{Right:} Modeling the same situation with an autonomous output system by replacing the external input with a respective signal generator.}
\end{figure}

\begin{enumerate}
\item Instead of considering the input-output map $\vecbj{u}\mapsto \vecbj{y}$ for arbitrary $\vecbj{u}$, we assume the input itself to be described by an autonomous quadratic differential system, the \textit{signal generator}. The input-output system \textit{driven by the signal generator} can then also be characterized by an enlarged autonomous output system, i.e., a system without any input, see Fig.~\ref{fig:io-sys-siggendriv} for an illustration.
\item We construct a \textit{variational expansion} of the autonomous signal generator driven system with respect to (w.r.t.) its initial conditions. This results in an infinite series of linear systems.
\item For the first few terms of the variational expansion we construct univariate frequency representations and perform an \textit{approximate moment matching}. This means the determination of the reduction basis $\tensorbj{V}$ corresponds to approximating certain interpolation conditions for the univariate representations.
\end{enumerate}
Signal generator driven systems have been used in literature to analyze systems for specifically described input signals. In \cite{art:isodori20081}, \cite{book:ast2008} an asymptotic limit behavior, a generalization of a steady state, is studied. The idea of using signal generators for model reduction can be found in  \cite{art:steadystate-mor-2010}, \cite{inproc:mor-steadystate-response}, \cite{inproc:families-nonlin-mom-match}. These works focus on the generalized concept of steady states and aim at identifying surrogate models that match the steady states for certain scenarios. 
Apart from the use of a signal generator our approach is different, it relies on the development of a new input-tailored variational expansion of the state. The work that probably shares most similarities with ours, and which initially inspired us to look deeper into the subject, is \cite{inproc:fast-nonlin-mor-assoc-trafo}, \cite{art:mor-associated-transform-2016}. The common feature is the univariate frequency representation derived for a variational expansion. Nonetheless, our approach exhibits profound differences to the former: 
Using the concept of signal generators we develop a framework that allows us to derive the variational expansion more rigorously and more generally. Our analysis suggests the incorporation of additional tensor-structured approximation conditions. Regarding the cascade- and low-rank tensor-structure in the approximation problems results in a more efficient implementation. The latter point is crucial for practical usage, as the involved univariate frequency representations grow vastly in dimension when considered as unstructured linear ones. It turns out that the \textit{exact} moment matching idea pursued classically in model reduction has to be relaxed to an \textit{approximate} moment matching owed to the tensor structure of the problem. In this respect, our input-tailored moment matching is more involved as the multi-moment approaches \cite{art:quadratic-bilinear-regular-krylov}, \cite{art:benner-goyal-gugercin2018}, \cite{phd:Gu-mor-nonlinear}, \cite{art:two-sided-hermm}, \cite{inproc:morBenB12}. However, our method corresponds to a one-dimensional interpolation problem unlike the multi-moment approaches corresponding to multi-dimensional interpolation problems. The latter consequently involve the choice of more expansion frequencies in multi-dimensional frequency spaces compared to ours, which requires fewer expansion frequencies from a one-dimensional frequency space.
A further difference to other system-theoretic reduction approaches is that ours extends very naturally to systems with more general input relations, such as, e.g., nonlinear functions and time derivatives. In this respect it is similarly flexible as the trajectory-based reduction methods like proper orthogonal decomposition \cite{art:morKunV01}, \cite{art:amsallam-gal-wave}.
As a byproduct of the extension of our method to more general input relations, we also derive a respective extension for system-theoretic methods relying on multivariate frequency representations by incorporating input-weights. Although the use of input-weights in model reduction is not new \cite{art:varga-freq-weighted-bal-rel}, \cite{art:breiten-near-optimal-freq-weighted}, they have -- to the best of the authors' knowledge -- not been applied for this purpose before.

\begin{figure}
\begin{tabular}{ccc}
\begin{minipage}{0.38\textwidth}

{\hspace*{0.08cm} \begin{tikzpicture}[scale=0.7, every node/.style={scale=0.8}]
\node 			(u) at (0.0, 0) {$\vecbj{u}$};
\node [rect2] 	(Sw) at (1.5, 0) {$\statemap$};
\node 		 	(x) at (3., 0) {$\vecbj{x}$};
\node [rect2] 	(Cx) at (4., 0) {$\tensorbj{C}$};
\node 		 	(y) at (5.4, 0) {$\vecbj{y}$};

\node [above right=0.33 and -2.0 of Cx, align = center]  {Input-output\\ system};

\draw [ ->] (0.45, 0) -- (1.0,0);
\draw [\tikzcolor, ->] (2., 0) -- (2.8,0);
\draw [\tikzcolor, ->] (3.2, 0) -- (3.6,0);
\draw [black, ->] (4.5, 0) -- (5.2,0);
\draw (1.1,-0.6) rectangle (4.4,0.7); \\

\node 			(VV1) at (1.4, -0.8) {};
\node 			(VV2) at (1.4, -3.3) {};
\draw [thick,dashed,bend angle=25, bend right, ->]  (VV1) to (VV2); 
\node [rotate=+90,align = center]			(VV12) at (0.5,-2.1) {\small{\textit{Drive by signal}}\\ \small{\textit{generator} $\siggen$}};
\end{tikzpicture} 
}

{\vspace*{-0.5cm}
\begin{tikzpicture}[scale=0.7, every node/.style={scale=0.8}]
\node [rect2]	(Su) at (0, 0) {$\statemapAut_w$};
\node 			(w) at (1., 0) {$\veccj{w}$};
\node [rect2] 	(Cw) at (2.1, 0) {$\tensorcj{P}_x$};
\node 		 	(x2) at (3., 0) {$\vecbj{x}$};
\node [rect2] 	(Cx2) at (4., 0) {$\tensorbj{C}$};
\node 		 	(y2) at (5.4, 0) {$\vecbj{y}$};
\node [above right=0.4 and -1.9 of Cx2, align = center]  {Autonomous output\\ system};
\node [above right= -0.04 and 0.05 of Cw, align = center]  {\Large{$\statemapAut$}};
\draw [\tikzcolor, ->] (0.45, 0) -- (0.8,0);
\draw [\tikzcolor, ->] (1.3, 0) -- (1.7,0); 
\draw [\tikzcolor, ->] (2.45, 0) -- (2.8,0);
\draw [\tikzcolor, ->] (3.2, 0) -- (3.6,0);
\draw [black, ->] (4.5, 0) -- (5.2,0);
\draw (-0.5,-0.6) rectangle (4.5,0.9);
\draw (-0.42,-0.4) rectangle (2.5,0.5);
\end{tikzpicture}
}
\end{minipage}

& 
\begin{minipage}{0.15\textwidth}
\hspace*{-0.9cm}
{\begin{tikzpicture}
\draw [line width=0.4mm, ->>] (0., 2.4) -- (2.4,2.4) node [pos=0.5, above, black,align = center] {${\tensorbj{V}}$} node [pos=0.5, below, black] {\tiny{\textit{Project}}};
\draw [line width=0.4mm, ->>] (0., -0.9) -- (2.4,-0.9) node [pos=0.5, above, black] {${\tensorcj{V}}$ } node [pos=0.5, below, black] {\tiny{\textit{Project}}};

\node 			(VV1) at (1.24, -0.6) {};
\node 			(VV2) at (1.3,2.7) {};
\draw [thick,dashed,bend angle=85, bend right, ->]  (VV1) to (VV2); 
\node [rotate=+90]			(VV12) at (2.0,1.0) {\tiny{\textit{Extract} $\tensorbj{V}$}};
\end{tikzpicture}
}
\end{minipage}

& 
\begin{minipage}{0.35\textwidth}

{
\begin{tikzpicture}[scale=0.7, every node/.style={scale=0.8}]
\node 			(u) at (0.0, 0) {$\vecbj{u}$};
\node [rect2] 	(Sw) at (1.5, 0) {$\statemap_{r}$};
\node 		 	(x) at (3., 0) {$\vecbj{x}_{r}$};
\node [rect2] 	(Cx) at (4., 0) {$\tensorbj{C}_{r}$};
\node 		 	(y) at (5.5, 0) {$\tilde{\vecbj{y}}$};

\node [above right=0.33 and -2.0 of Cx, align = center]  {Input-output\\ system (reduced)};

\draw [ ->] (0.45, 0) -- (1.0,0);
\draw [\tikzcolor, ->] (2., 0) -- (2.65,0);
\draw [\tikzcolor, ->] (3.2, 0) -- (3.6,0);
\draw [black, ->] (4.6, 0) -- (5.3,0);
\draw (1.1,-0.6) rectangle (4.5,0.7);

\node 			(VV1) at (1.8, -0.8) {};
\node 			(VV2) at (1.8, -3.3) {};
\draw [thick,dashed,bend angle=25, bend right, ->]  (VV1) to (VV2); 
\node [rotate=+90,align = center]			(VV12) at (0.9,-2.1) {\small{\textit{Drive by signal}}\\ \small{\textit{generator} $\siggen$}};
\end{tikzpicture} 
}

{\vspace*{-0.5cm}
\begin{tikzpicture}[scale=0.7, every node/.style={scale=0.8}]
\node [rect3]	(Su) at (0, 0) {$\statemapAut_{w,r}$};
\node 			(w) at (1., 0) {$\veccj{w}_r$};
\node [rect3] 	(Cw) at (2., 0) {$\tensorcj{P}_{x_r}$};
\node 		 	(x2) at (3., 0) {$\vecbj{x}_r$};
\node [rect3] 	(Cx2) at (4., 0) {$\tensorbj{C}_r$};
\node 		 	(y2) at (5.5, 0) {$\tilde{\vecbj{y}}$};
\node [above right=0.4 and -1.9 of Cx2, align = center]  {Autonomous output\\ system (reduced)};
\node [above right= -0.1 and 0.0 of Cw, align = center]  {\Large{$\statemapAut_r$}};
\draw [\tikzcolor, ->] (0.4, 0) -- (0.7,0);
\draw [\tikzcolor, ->] (1.25, 0) -- (1.55,0); 
\draw [\tikzcolor, ->] (2.45, 0) -- (2.8,0);
\draw [\tikzcolor, ->] (3.2, 0) -- (3.6,0);
\draw [black, ->] (4.6, 0) -- (5.3,0);
\draw (-0.65,-0.6) rectangle (4.5,0.9);
\draw (-0.5,-0.4) rectangle (2.45,0.5);
\end{tikzpicture}
}
\end{minipage}

\end{tabular}
\caption{Sketch for input-tailored moment matching based on the signal generator driven system and reduction via Galerkin projection.}
\label{fig:input-tailored-sketch}
\end{figure}

\subsection*{Outline}
The outline of this manuscript is as follows: The concept of a signal generator driven system as well as the proposed variational expansion and associated univariate frequency representation of the resulting autonomous system are presented in Section~\ref{sec:tf-qbsys}. We refer to the expansion and frequency representations as input-tailored, as they take into account the input described by the signal generator. The approximation conditions, which our reduction method aims for, resembles an approximate moment matching condition of the input-tailored frequency representations (Section~\ref{sec:towards-input-tailored-syst}). In this context the commuting diagram of Fig.~\ref{fig:input-tailored-sketch} takes a prominent role. Our numerical realization is presented in Section~\ref{sec:approx-mom-match}.  In Section~\ref{sec:ext-var-appr-app} we discuss the ability of handling non-standard input dependencies in our method. Moreover we suggest an extension for other system-theoretic methods to handle non-standard input maps, which falls off as a byproduct of the discussion of our approach. The performance of our input-tailored moment matching method in comparison to the system-theoretic multi-moment matching and the trajectory-based proper orthogonal decomposition as well as the proposed handling of non-standard input maps are numerically studied in Section~\ref{sec:assocmor-num-validation}. The four appendices provide expressions for higher-order univariate frequency representations and details for the derivation of the variational expansion and generalizations.

\subsection*{Notation}
Throughout this paper, matrices/tensors, vectors and scalars are indicated by capital boldfaced, small boldfaced and normal letters, respectively. Moreover, in the typeface we distinguish between the quantities associated to the original input-output system (e.g., $\statemap$, $\tensorbj{A}$, $\vecbj{x}$) and the ones associated to the signal generator driven system (e.g., $\statemapAut$, $\tensorcj{A}$, $\veccj{w}$).  Frequency representations are written in a curved font (e.g., $\Transfer{X}$, $\Transfer{W}$). The subscript $_r$ indicates a reduced quantity gotten by Galerkin projection (cf.\ Fig.~\ref{fig:input-tailored-sketch}).

Moreover, tensor notation is used within the paper, cf.\ \cite{book:nonlinear-system-theory-rugh}, \cite{art:kressner-kryl-tensor}, \cite{book:hackbusch-tensor}. The Kronecker-tensor product is denoted by $\otimes$, it is defined as 
\begin{align*}
		\tensorbj{P} \otimes \tensorbj{Q} = 
	\begin{mymatrix}
		p_{11}\tensorbj{Q} & p_{12}\tensorbj{Q} & \ldots & p_{1N}\tensorbj{Q} \\
		\ldots \\
		p_{M,1}\tensorbj{Q} & p_{M2}\tensorbj{Q}  & \ldots & p_{MN}\tensorbj{Q}
	\end{mymatrix} \quad \text{for } \tensorbj{P}, \tensorbj{Q} \in \mathbb{R}^{M,N}, \, \tensorbj{P}=(p_{ij}).
\end{align*}
We abbreviate $\vecbj{P}\supertens{2} = \vecbj{P} \otimes \vecbj{P}$, $\vecbj{P}{\supertens{3}} = \vecbj{P} \otimes \vecbj{P} \otimes \vecbj{P}$. Additionally, we introduce the notation
\begin{align*}
\circledmy{2}_{\tensorbj{P}}\tensorbj{Q} &= \tensorbj{Q} \otimes \tensorbj{P} + \tensorbj{P} \otimes \tensorbj{Q} &\in \mathbb{R}^{M^2,N^2} \\
	\circledmy{3}_{\tensorbj{P}}\tensorbj{Q} &= \circledmy{2}_{\tensorbj{P}}\tensorbj{Q} \otimes \tensorbj{P} + \tensorbj{P}\supertens{2} \otimes \tensorbj{Q}  &\\
		& = \tensorbj{Q} \otimes \tensorbj{P} \otimes \tensorbj{P} + \tensorbj{P} \otimes \tensorbj{Q} \otimes \tensorbj{P} + \tensorbj{P} \otimes \tensorbj{P} \otimes \tensorbj{Q} &\in \mathbb{R}^{M^3,N^3}.
\end{align*}
The expressions $\vecbj{P}\supertens{i}$ and $\circledmy{{\relsize{-1.5}i}}_{\tensorbj{P}}\tensorbj{Q}$ are defined analogously for $i >3$.
By definition it holds
\begin{align*}
	\begin{mymatrix}
		\tensorbj{A} & \tensorbj{B} \\
		\tensorbj{C} & \tensorbj{D}		
	\end{mymatrix}
	\otimes \tensorbj{P} &=
		\begin{mymatrix}
		\tensorbj{A} \otimes \tensorbj{P} & \tensorbj{B} \otimes \tensorbj{P} \\
		\tensorbj{C} \otimes \tensorbj{P} & \tensorbj{D} \otimes \tensorbj{P}
	\end{mymatrix} \\
	(\tensorbj{A} \tensorbj{B}) \otimes (\tensorbj{C}\tensorbj{D}) &=  (\tensorbj{A} \otimes \tensorbj{C}) \, (\tensorbj{B} \otimes \tensorbj{D})
\end{align*}
for matrices of appropriate dimensions. The unit matrix and the zero matrix are denoted by $\vecbj{I}_N\in \mathbb{R}^{N,N}$ and $\vecbj{0}_{M,N}\in \mathbb{R}^{M,N}$, respectively, whereby the sub-index of the dimensions is omitted if they are clear from the context. For vectors $\vecbj{p}\in \mathbb{R}^M$, $\vecbj{q}\in \mathbb{R}^N$, we often use the notation
\begin{align*} 
[\vecbj{p}; \vecbj{q}] = 
	\begin{mymatrix}
		\vecbj{p}\\
		\vecbj{q}
	\end{mymatrix} \in \mathbb{R}^{M+N}.
\end{align*}

\section{Input-tailored expansion and frequency representation} \label{sec:tf-qbsys}

In this section we develop our input-tailored variational expansion and frequency representation our reduction method is based on. Starting point is the concept of signal generator driven systems (Section~\ref{subsec:sig-gen-driven-sys}). These signal generator driven systems are, by construction, autonomous. Variational expansions of autonomous systems and associated univariate frequency representations are the topic of Section~\ref{subsec-variational-wrt-init}. Our input-tailored expansion and frequency representation are presented in Section~\ref{subsec:input-tailor-var-exp}. In Section~\ref{subsec:mult-tf}, we embed our proposed expansion in the existing literature by relating it to the Volterra series and its frequency representations.

\subsection{Signal generator driven system} \label{subsec:sig-gen-driven-sys}

In focus of this paper are quadratic-bilinear dynamical systems of the form
\begin{subequations} \label{eq:qldae}
\begin{align} 
	\statemap: \quad \tensorbj{E} \dot{\vecbj{x}} &= \tensorbj{A} \vecbj{x} + \tensorbj{G} \vecbj{x}\supertens{2} + \tensorbj{D} (\vecbj{x} \otimes \vecbj{u}) + \tensorbj{B} \vecbj{u}, \quad \vecbj{x}(0)= \vecbj{x}_0 \in \mathbb{R}^N \label{eq:qldae-a} \\
	\vecbj{y} &= \tensorbj{C} \vecbj{x}, \quad  \vecbj{u}(t) \in \mathbb{R}^p, \quad t \geq 0
\end{align}
with $\tensorbj{E}$ nonsingular, all system matrices constant and $\tensorbj{G} \in \mathbb{R}^{N,N^2}$, $\tensorbj{D} \in \mathbb{R}^{N,Np}$. By slight abuse of notation, we identify throughout the paper the realization of the state equation $\statemap$ with its input-to-state map $\statemap:\vecbj{u}\mapsto \vecbj{x}$.

Instead of considering $\statemap$ directly as abstract map, we use the concept of a signal generator driven system. A signal generator is an autonomous differential system describing the input $\vecbj{u}$. We employ here the class of signal generators with quadratic nonlinearities given as
	\begin{align} \label{eq:siggen}
			\siggen:& \qquad \vecbj{u} = \tensorbj{C}_{\sigstatescal} \sigstate, \qquad \dot{\sigstate} = \tensorbj{A}_{\sigstatescal} \sigstate + \tensorbj{G}_{\sigstatescal}\sigstate\supertens{2}, \qquad \sigstate(0) = \sigstate_{0} \in \mathbb{R}^q.
	\end{align}
\vspace{-0.3cm}	
\end{subequations}
\begin{remark}[Signal generators]\label{rem:siggen-ex}
For example, an oscillation $u(t) = a \sin{(\lambda t})$ for $t\geq 0$ and $a, \lambda \in \mathbb{R}$ is readily given by the signal generator
\begin{align*}
				u = [1\,|\,0] \sigstate, \qquad \dot{\sigstate} = 
	\lambda \begin{mymatrix}
		 & 1 \\
		-1 &
\end{mymatrix}	 \sigstate,
	  \qquad \,\sigstate(0) = 	 
\begin{mymatrix}
	0 \\
	a
\end{mymatrix}.
\end{align*}
More generally, any linear combination of exponential pulses and sine- and cosine-oscillations can be described by a linear signal generator (as in \eqref{eq:siggen} with $\tensorbj{G}_{\sigstatescal}~=~\tensorbj{0}$) by superposition of simple signal generators. Taking, e.g., $u(t) = a_1 \exp{(\lambda_1 t)} + a_2 \cos{(\lambda_2 t)}$, the associated signal generator reads
\begin{align*}
				u = [1\,|\,0\,|1] \sigstate, \qquad \dot{\sigstate} = 
	 \begin{mymatrix}
		\lambda_1 & & \\
		& & \lambda_2\\
		& -\lambda_2 &
\end{mymatrix}		
	 \sigstate,
	  \qquad \,\sigstate(0) = \begin{mymatrix}
	a_1\\
	0 \\
	a_2
\end{mymatrix}.
\end{align*}
Arbitrary derivatives in frequency space of the above mentioned functions, such as e.g., $u(t) = t^k \exp(t)$ for $k \in \mathbb{N}$, can also be described with linear signal generators. With nonlinear signal generators an even larger class of inputs can be represented, see e.g., \cite{book:ast2008}, \cite{art:steadystate-mor-2010} for a broader discussion and some applications. We also refer to Section~\ref{sec:burgers-eq}, Case~2 for an example of a quadratic signal generator.
\end{remark}
Similarly to \cite{art:steadystate-mor-2010}, \cite{inproc:mor-steadystate-response}, \cite{inproc:families-nonlin-mom-match}, the notion of a signal generator driven system is defined in the upcoming. It results from inserting a signal generator for the input $\vecbj{u}$ in system $\statemap$.
\begin{definition}[Signal generator driven system] \label{def:sg-ode-casc}
Let a quadratic-bilinear system $\statemap$ with an input $\vecbj{u}$ described by the signal generator $\siggen$ as in \eqref{eq:qldae} be given,
\begin{align*}
	\statemap:& \quad \tensorbj{E} \dot{\vecbj{x}} = \tensorbj{A} \vecbj{x} + \tensorbj{G}\vecbj{x}\supertens{2} + \tensorbj{D}  (\vecbj{x} \otimes \vecbj{u})  +  \tensorbj{B} \vecbj{u} , \quad &&
	\vecbj{x}(0) = \vecbj{x}_0 \in \mathbb{R}^{N}\\
	\siggen:& \quad \,\,\,\,\, \vecbj{u} = \tensorbj{C}_{\sigstatescal} \sigstate, \qquad \dot{\sigstate} = \tensorbj{A}_{\sigstatescal} \sigstate + \tensorbj{G}_{\sigstatescal} \sigstate\supertens{2}, \qquad && \, \sigstate(0) = \sigstate_{0} \in \mathbb{R}^q.
\end{align*}
Let $\tensorbj{Q}$ be the constant matrix such that
\begin{align*}
	\tensorbj{Q} \begin{mymatrix}
		\bar{\vecbj{x}} \\
		\bar{\sigstate}
	\end{mymatrix}\supertens{2} = 
	\begin{mymatrix}
		\bar{\vecbj{x}}\supertens{2} \\
		\bar{\vecbj{x}} \otimes \bar{\sigstate} \\
		\bar{\sigstate}\supertens{2}
	\end{mymatrix}  \nonumber \qquad  \text{for arbitrary  } \bar{\vecbj{x}} \in \mathbb{R}^N,\,  \bar{\sigstate}\in \mathbb{R}^q.
\end{align*}
Then we call the autonomous system
\begin{align*}
\parbox{0.5cm}{{\vspace{0.4cm}$\statemapAut:$}} 
& \quad	\,	\tensorcj{E} \dot{\veccj{w}} = \tensorcj{A} \veccj{w} + \tensorcj{G} \veccj{w}\supertens{2},  && \veccj{w}(0) = \veccj{b} \\
& \qquad  \vecbj{x} = \tensorcj{P}_x \, \veccj{w} &&
\end{align*}
with
\begin{align*}
	\tensorcj{E} &= \begin{mymatrix}
		\tensorbj{E} & \\
					& \tensorbj{I}_q
	\end{mymatrix}
	, \qquad 
	\tensorcj{A} =
		\begin{mymatrix}
		\tensorbj{A} & \tensorbj{B} \tensorbj{C}_{\sigstatescal}  \\
					& \tensorbj{A}_{\sigstatescal}
	\end{mymatrix}, \qquad
	\tensorcj{P}_x = [\tensorbj{I}_N, \, \tensorbj{0}],
	\qquad \veccj{b} = 
	\begin{mymatrix}
		\vecbj{x}_0 \\
		 \sigstate_0
	\end{mymatrix},
\\
	\tensorcj{G} &=
		 \begin{mymatrix}
		\tensorbj{G} & \tensorbj{D}(\tensorbj{I}_N \otimes\tensorbj{C}_{\sigstatescal}) & \\
					& 													  & \tensorbj{G}_{\sigstatescal}
	\end{mymatrix}	\tensorbj{Q}  
\end{align*}
the signal generator driven system $\statemapAut$.
\end{definition}
By definition, the solution $\vecbj{x}$ of system $\statemap$ for input $\vecbj{u}$ described by the signal generator $\siggen$ and the output $\vecbj{x}$ of the signal generator driven system $\statemapAut$ coincide. For an illustration, we refer to Fig.~\ref{fig:input-tailored-sketch}, left column. Note that the state equation of $\statemapAut$ (denoted by $\statemapAut_w$ in Fig.~\ref{fig:input-tailored-sketch}) is autonomous.

\subsection{Variational expansion of autonomous systems and associated univariate frequency representations} \label{subsec-variational-wrt-init}

Our approach employs a variational expansion of the autonomous system $\statemapAut$ from Definition~\ref{def:sg-ode-casc} and associated univariate frequency representations. The theoretical basis is given by the following theorem.

\begin{theorem} \label{theor:Pxu_quadratic}
Let an $\alpha$-dependent initial value problem of the autonomous quad\-ratic differential equation
\begin{align*}
	\tensorcj{E} \dot{\veccj{w}}(t;\alpha) &= \tensorcj{A} \veccj{w}(t;\alpha) + \tensorcj{G} \, (\veccj{w}(t;\alpha))\supertens{2}, \qquad t \in (0,T) \\
		 \veccj{w}(0;\alpha) &= \alpha \veccj{b}
\end{align*}
be given for $T >0$ and constant system matrices $\tensorcj{E},\tensorcj{A} \in \mathbb{R}^{M,M}$, $\tensorcj{G} \in \mathbb{R}^{M,M^2}$ and $ \veccj{b}\in \mathbb{R}^M$ with $\tensorcj{E}$ nonsingular. For parameter $\alpha \in I$, $0\in I \subset \mathbb{R}$ bounded interval, the family of $\alpha$-dependent solutions $\veccj{w}(\cdot,\alpha)$ can then be expanded as
\begin{align}\label{eq:series}
	\veccj{w}(t;\alpha) =  \sum_{i=1}^{N} \alpha^i \veccj{w}_i(t)  + \text{O}(\alpha^{N+1}), \qquad t \in [0,T), \quad \alpha \in I.
\end{align}
The univariate frequency representations $\breve{\Transfer{W}}_i$ of the first three functions $\veccj{w}_i$ for $s \in \mathbb{C}$ are 
 \begin{subequations} \label{eq:assoc-tf}
	\begin{align} 
\breve{\Transfer{W}}_1(s) &= (s \tensorcj{E} - \tensorcj{A} )^{-1} \veccj{b}\\
\breve{\Transfer{W}}_2(s) &= (s \tensorcj{E} - \tensorcj{A} )^{-1} \tensorcj{G} \left(s \tensorcj{E}\supertens{2} - \circledmy{2}_\tensorcj{E} \tensorcj{A}\right)^{-1} \veccj{b}\supertens{2} \label{eq:assoc-tf-2}\\
	\breve{\Transfer{W}}_3(s) &= 2 (s \tensorcj{E} - \tensorcj{A} )^{-1} \tensorcj{G}  \left(s \tensorcj{E}\supertens{2} - \circledmy{2}_\tensorcj{E} \tensorcj{A} \right)^{-1}  (\tensorcj{G} \otimes \tensorcj{E})
 \left(s \tensorcj{E}{\supertens{3}} - \circledmy{3}_\tensorcj{E} \tensorcj{A} \right)^{-1} \veccj{b}{\supertens{3}}. \label{eq:assoc-tf-3}
	\end{align}
\end{subequations}
\end{theorem} 
The proof of Theorem~\ref{theor:Pxu_quadratic} relies on a variational expansion w.r.t.\ the initial conditions and on frequency space formulations using the so-called Associated Transform \cite{book:nonlinear-system-theory-rugh}. It is provided in detail in Appendix \ref{app:proof-theorem-var-input}. Formal similarities to univariate frequency representations of \cite{inproc:fast-nonlin-mor-assoc-trafo}, \cite{art:mor-associated-transform-2016} are addressed and exploited within our proof. Let us emphasize that the cited works rely on variational expansions w.r.t.\ the inputs, which distinguishes their approach from ours.

\begin{remark}
Certainly, the series in \eqref{eq:series} can be formulated regarding terms of arbitrarily high order in $\alpha$. The tensor-structured explicit representations, however, get lengthy for high orders and the calculations more technical. In the main body of the paper, we restrict ourselves from now on to terms up to order two to keep it more comprehensible. The tensor structure patterns that are observed and exploited for order two are preserved for the higher-order expressions as well. For order three this can be seen in Theorem~\ref{theor:Pxu_quadratic} and in respective generalizations of other important results provided in Appendix~\ref{app:lem:mom-assoc-tf3}.
\end{remark}
Another point of view on the associated univariate frequency representation $\breve{\Transfer{W}}_2$ is highlighted in the following lemma that results from straightforward calculus (cf.\ Lemma~\ref{lem:assoc-trafo-repres_A} for $\breve{\Transfer{W}}_3$).

\begin{lemma} \label{lem:assoc-trafo-repres}
Assume that the requirements of Theorem~\ref{theor:Pxu_quadratic} hold true. Then the associated frequency representation $\breve{\Transfer{W}}_2$ can be formulated with the linear representation
\begin{align*}
\breve{\Transfer{W}}_2(s) &= \breve{\tensorcj{C}}_2 \left(s \breve{\tensorcj{E}}_2 - \breve{\tensorcj{A}}_2 \right)^{-1} \breve{\veccj{b}}_2,\\
			 & \text{with} \quad \breve{\tensorcj{E}}_2 =
			\begin{mymatrix}
				\tensorcj{E} & \\
				  		& \tensorcj{E}\supertens{2}
			\end{mymatrix}	
	,\quad	
				  		\breve{\tensorcj{A}}_2 =
			\begin{mymatrix}
				\tensorcj{A} & \tensorcj{G} \\
				  	&  		\circledmy{2}_{\tensorcj{E}}\tensorcj{A}
			\end{mymatrix}, \quad 
			\breve{\veccj{b}}_2 = 
			\begin{mymatrix}
				\tensorbj{0} \\
				\veccj{b}\supertens{2}
			\end{mymatrix}, \quad 
						\breve{\tensorcj{C}}_2 = 
			\begin{mymatrix}
				\tensorbj{I}_M & \tensorbj{0}
			\end{mymatrix}.
\end{align*}
\end{lemma}
\begin{remark}[Cascade- and tensor-structure of associated frequency representations] \label{rem:assoc-cascaded-structure}
The frequency representation $\breve{\Transfer{W}}_1$ associated to the first order term of the variational expansion is a usual linear input-to-state transfer function with dimension $M$  equal to the dimension of the state $\veccj{w}$. According to Lemma~\ref{lem:assoc-trafo-repres} (and Lemma~\ref{lem:assoc-trafo-repres_A}), also the higher-order terms possess linear state representations, which will strongly motivate our subsequently proposed procedure for setting up the approximation conditions in the approximate moment matching. However, since the frequency representations are of growing dimension, $\mathbb{R}^{M+M^2}$ for $\breve{\Transfer{W}}_2$ ($\mathbb{R}^{M+M^2+M^3}$ for $\breve{\Transfer{W}}_3$), operating directly on them -- as done in \cite{inproc:fast-nonlin-mor-assoc-trafo}, \cite{art:mor-associated-transform-2016} -- is unpractical for medium- to large-scale problems. For the development of a numerically tractable method, we instead exploit their special cascade- and tensor-structure that is revealed in Theorem~\ref{theor:Pxu_quadratic}. For example, $\breve{\Transfer{W}}_2$ can be interpreted as the cascade of the transfer functions $\tensorcj{G} \left(s \tensorcj{E}\supertens{2} - \circledmy{2}_\tensorcj{E} \tensorcj{A} \right)^{-1}\veccj{b}\supertens{2}$ and $(s \tensorcj{E} - \tensorcj{A} )^{-1} $, where the former has low-rank tensor structure. 
\end{remark}

\subsection{Input-tailored variational expansion} \label{subsec:input-tailor-var-exp}

Based on the signal generator driven system we can now formulate our input-tailored expansion.

\begin{definition}[Input-tailored variational expansion] \label{def:input-tailor-var-exp}
Let the signal generator driven system $\statemapAut$ with enlarged state $\veccj{w}$ be as in Definition~\ref{def:sg-ode-casc}. Let
\begin{align*}
	\veccj{w}(t;\alpha) = \sum_{i=1}^{N} \alpha^i \veccj{w}_i(t) + \text{O}(\alpha^{N+1}), \qquad t \in [0,T)
\end{align*} 
be the variational expansion of $\veccj{w}$ w.r.t.\ the initial conditions $\veccj{w}(0;\alpha) = \alpha \veccj{b}$. Let $\breve{\Transfer{W}}_i$ be the associated univariate frequency representations of $\veccj{w}_i$ as in Theorem~\ref{theor:Pxu_quadratic}.

Then the input-tailored variational expansion of $\vecbj{x}$ described by $\statemapAut$ (respectively by $\statemap$ and $\siggen$) is defined as
\begin{align*}
	\vecbj{x}(t;\alpha) = \sum_{i=1}^{N} \alpha^i \vecbj{x}_i(t) + \text{O}(\alpha^{N+1}), \qquad \quad \vecbj{x}_i(t) =  \tensorcj{P}_x \, \veccj{w}_i(t)
\end{align*}
with $\tensorcj{P}_x =[\tensorbj{I}_N, \tensorbj{0}_{N,q}]$. The input-tailored frequency representations $\breve{\Transfer{X}}_i$ are given as
\begin{align*}
	\breve{\Transfer{X}}_i(s) = \tensorcj{P}_x \breve{\Transfer{W}}_i(s), \qquad s \in \mathbb{C}.
\end{align*}
\end{definition}
Let us emphasize that our input-tailored variational expansion is not tailored towards a single solution trajectory, but rather towards a family of solutions parametrized in the expansion parameter~$\alpha$. Given, e.g., the signal generator
\begin{align*} 
			u = [1\,| 0] \sigstate, \qquad \dot{\sigstate} = 	\lambda \begin{mymatrix}
		 & 1\\
		 -1 &
\end{mymatrix}	 \sigstate
	  \qquad \,\sigstate(0) = \alpha \begin{mymatrix}
	0 \\
	1
\end{mymatrix}, \qquad \alpha \in \mathbb{R},
\end{align*}
it relates to the inputs $u(t) = \alpha \sin(\lambda t)$, i.e., oscillations of varying amplitude. This clearly distinguishes our approach from snapshot based methods like proper orthogonal decomposition \cite{art:morKunV01}.


\begin{remark}[Possible generalizations]
We point out that the definition of signal generator driven systems, Definition~\ref{def:sg-ode-casc}, and with that our whole approach can be generalized straightforwardly to systems with more sophisticated input maps, e.g., quadratic inputs, time derivatives, see Section~\ref{sec:ext-var-appr-app}.

Moreover, the variational expansion from Theorem~\ref{theor:Pxu_quadratic} itself can be generalized. Instead of considering families of solutions parametrized in initial conditions that dependent only on a single parameter $\alpha$, also families of solutions parametrized in a multidimensional parameter can be treated, see Appendix~\ref{app:var-multdim-init}. This includes solutions parametrized in inputs $u(t) = \sum_j \alpha_j u_j(t)$ for varying $\alpha_j$, where all $u_j$ have a linear signal generator.
\end{remark}

\subsection{Relation to Volterra series expansion} \label{subsec:mult-tf}

In the following, we discuss the relation of our input-tailored variational expansion with the Volterra series, which is a variational expansion of the solution w.r.t.\ the input. The Volterra series has recently been extensively used as a basis for model reduction. For example, multi-moment matching has been discussed in \cite{phd:Gu-mor-nonlinear}, \cite{inproc:morBenB12}, \cite{techrep:two-sided-mm}, hermite multi-moment matching in \cite{art:two-sided-hermm}, \cite{art:quadratic-bilinear-regular-krylov}, \cite{art:benner-goyal-gugercin2018}, and balanced truncation in \cite{art:quadratic-bililinear-bt2017}. We recapitulate the variational ansatz from \cite{book:nonlinear-system-theory-rugh}, \cite{art:volterra-exist-and-unique}, \cite{art:func-exp-gilbert}. As the references are restricted to the scalar input case $u:\mathbb{R} \rightarrow \mathbb{R}$, we also use this restriction for convenience. Consider the state equation $\statemap$ with a scalar-valued input and trivial initial conditions, i.e.,
\begin{align*}
	\tensorbj{E} \dot{\vecbj{x}} &= \tensorbj{A} \vecbj{x} + \tensorbj{G} \vecbj{x}\supertens{2} +  u \,\tensorbj{D} \vecbj{x}  + \tensorbj{b} u ,  \quad \vecbj{x}(0) = \vecbj{0}, \quad u: \mathbb{R} \rightarrow \mathbb{R}
\end{align*}
with $\tensorbj{b} \in \mathbb{R}^{N}$. For appropriate input $u(t) = \alpha v(t)$ with $\alpha \in \mathbb{R}$ sufficiently small and under the assumption that the system is uniquely solvable in an $\alpha$-neighborhood containing zero, a variational expansion in the input exists. That is, the solution can be expanded in $\alpha$ for $N>0$ as
\begin{align} \label{eq:var-ansatz-alpha}
	\vecbj{x}(t;\alpha) = \sum_{i=1}^N \alpha^i \vecbj{x}_i(t) + \text{O}(\alpha^{N+1}) \qquad t \in [0,T)
\end{align}
for some $T>0$. It can be shown, using the multivariate Laplace transform as in \cite{art:benner-goyal-gugercin2018}, that the terms $\vecbj{x}_i$ have multivariate frequency representations $\Transfer{X}_i$ with
\begin{align*}
	\Transfer{X}_1(s_1) &= \Transfer{G}_1(s_1) \Transfer{U}(s_1), \\
	 \Transfer{X}_2(s_1,s_2) &= \Transfer{G}_2(s_1,s_2) \Transfer{U}(s_1) \Transfer{U}(s_2) \\
	\Transfer{X}_i(s_1,s_2,\ldots,s_i) &= \Transfer{G}_i(s_1,s_2,\ldots,s_i) \Transfer{U}(s_1) \Transfer{U}(s_2) \ldots \Transfer{U}(s_i),&& s_i \in \mathbb{C}, \quad i \leq N
\end{align*}
where $\Transfer{U}$ is the Laplace transform of the input $u$ and $\Transfer{G}_i$ are the so-called symmetric transfer functions, see \cite{book:nonlinear-system-theory-rugh}, \cite{art:compact-red-mod-Li}, \cite{art:mor-associated-transform-2016} for details on them. The model reduction methods relying on the Volterra series \eqref{eq:var-ansatz-alpha} typically formulate approximation conditions for the transfer functions ${\Transfer{G}}_i$.

At first glance there seems not to be a connection to our input-tailored variational expansion. The upcoming lemma, however, shows that for inputs described by linear signal generators, both expansions lead to the same result.

\begin{lemma}
Let a quadratic-bilinear differential system $\statemap$ with a linear signal generator $\siggen$ be given as
\begin{align*}	
	 \statemap: \,\, \quad \tensorbj{E} \dot{\vecbj{x}} &= \tensorbj{A} \vecbj{x} + \tensorbj{G} \vecbj{x}\supertens{2} +  u \,\tensorbj{D} \vecbj{x} + \tensorbj{b} {u}, && \vecbj{x}(0) = \vecbj{0}, \\
	\siggen: \quad \quad  {u} &= \mathbf{C}_{\sigstate} \sigstate, \qquad \dot{\sigstate} = \mathbf{A}_{\sigstate} \sigstate, && \sigstate(0) = \sigstate_0.
\end{align*}
 Then the same expansion of the solution
\begin{align*}
 \vecbj{x}(t;\alpha)  &=  \sum_{i=1}^k \alpha^i \vecbj{x}_i(t) + O(\alpha^{k+1})
 \end{align*}
can be obtained by two different approaches:
\begin{enumerate}[a)]
	\item \label{it:var-approach-relation-odesg} By the input-tailored variational expansion of $\vecbj{x}$ as in Definition \ref{def:input-tailor-var-exp}.
	\item \label{it:var-approach-relation-volt} By the Volterra series: Expand the state $\vecbj{x}$ for input $u(t) = \alpha {u}_1(t)$ in $\alpha$, and then set the input $u$ to be as in the linear signal generator.
\end{enumerate}
\end{lemma}
\begin{proof}
Proceeding from Approach~\ref{it:var-approach-relation-odesg}) we show the equality to Approach~\ref{it:var-approach-relation-volt}). In Approach~\ref{it:var-approach-relation-odesg}) we assume for initial value $\veccj{w}_0 = \alpha [\vecbj{0}; \bar{\sigstate}_0]$ that the extended state can be expanded as
\begin{align*}
 \veccj{w}(t;\alpha)  &=
 \begin{mymatrix}
 \vecbj{x}(t;\alpha) \\
 \sigstate(t;\alpha)
\end{mymatrix} 
= \sum_{i=1}^k \alpha^i 
\begin{mymatrix}
 \vecbj{x}_i(t) \\
 \sigstate_i(t)
\end{mymatrix} 
+ O(\alpha^{k+1}).
 \end{align*}
From the signal generator relation $\siggen$ it then follows
\begin{align*}
 {u}(t;\alpha)  &=  \sum_{i=1}^k \alpha^i u_i(t) + O(\alpha^{k+1}), \qquad u_i(t)= \mathbf{C}_{\sigstate} \sigstate_i(t).
 \end{align*}
As the signal generator is linear, it is easily seen that $\sigstate \equiv \alpha \sigstate_1$. It therefore also holds ${u} \equiv \alpha u_1$, and the expansion terms $\alpha^i\vecbj{x}_i$ scale with $u^i \equiv \alpha^i u_1^i$ as in Approach $\ref{it:var-approach-relation-volt})$. Hence, the expansion terms $\vecbj{x}_i$ of both approaches coincide.
\end{proof} 
Inputs described by linear signal generators are an important case. Alternatively to the derivation in  \cite{art:benner-goyal-gugercin2018}, the multivariate symmetric transfer functions $\Transfer{G}_i$ can already be derived by only considering the response to sums of exponential functions
\begin{align*}
	u(t) = \sum_{k=1}^{i} a_k \exp(\lambda_k t), \qquad \text{for arbitrary } a_k, \lambda_k \in \mathbb{R},
\end{align*}
which is, e.g., used in the growing exponential approach, \cite{book:nonlinear-system-theory-rugh}, \cite{morBre13}. Clearly, sums of exponential functions can be described by linear signal generators, cf.\ Remark~\ref{rem:siggen-ex}. Therefore, loosely spoken, the associated univariate input-tailored frequency representation tailored towards the upper growing exponentials for different choices $a_k$, $\lambda_k$ resemble the multivariate transfer functions $\Transfer{G}_i$. The works \cite{art:autonomous-volterra-steady}, \cite{art:mor-associated-transform-2016} indirectly heavily rely on the upper resemblance, but do not explicitly elaborate on it.

Finally, let us comment on the more formal approach by \cite{art:mor-associated-transform-2016}, \cite{inproc:fast-nonlin-mor-assoc-trafo} that leads to similar univariate frequency representations as ours.

\begin{remark} \label{rem:formal-volterra-dirac}
In \cite{art:mor-associated-transform-2016}, \cite{inproc:fast-nonlin-mor-assoc-trafo} the quadratic-bilinear equation of Theorem~\ref{theor:Pxu_quadratic} with zero (pre-)initial conditions but an initial jump is considered, i.e., 
\begin{align*}
	\tensorcj{E} \dot{\veccj{w}} &= \tensorcj{A} \veccj{w} + \tensorcj{G} \veccj{w}\supertens{2} + \veccj{b} u(t), \quad u(t) = \alpha \delta(t),\quad
		  \lim_{\bar{t}\uparrow 0} \veccj{w}(\bar{t}) = \vecbj{0},
\end{align*}
where $\delta(t)$ is the Dirac-impulse. There the solution $\veccj{w}$ is expanded formally as Volterra series with that distributional input $u(t) = \alpha \delta(t)$, yielding the same expansion terms as ours. However, the validity of the Volterra series when the input is a Dirac-impulse is not covered by the classical result on Volterra series expansions -- as far as the authors know (cf., e.g., \cite{book:nonlinear-system-theory-rugh}, \cite{art:volterra-exist-and-unique}, \cite{art:func-exp-gilbert} or \cite{art:volterra-impulse-consid-borys}). This issue is also not further addressed or discussed in the respective works.
\end{remark}
	
\section{Input-tailored system-theoretic model reduction framework} \label{sec:towards-input-tailored-syst}

Aim of our method is to construct a reduced model such that for the input-tailored frequency representations $\breve{\Transfer{X}}_i$ the so-called moments
\begin{align*}
\frac{d^k}{ds^k} \breve{\Transfer{X}}_{i}(s)_{|s=s_0} \qquad \text{for } k, \, i, \, s_0  \text{ given}
\end{align*}
of the full order model are \textit{approximately} matched by their reduced counterparts. This is a relaxation of the linear moment matching idea, which we recapitulate in Section~\ref{subsec:exact-moments}. Our input-tailored moment matching problem is formulated in Section~\ref{subsec-siggendriv_r}. The notion of a signal generator driven system $\statemapAut$ and its reduced counterpart is herefore essential. The structure of the approximation problem is analyzed in Section~\ref{subsec:input-tailored-mom-lin}. From a theoretical point of view, it can be characterized with linear theory. To do so, a change to high-dimensional state representations (cf.\ Lemma~\ref{lem:assoc-trafo-repres}) is needed. Our projection ansatz, however, operates on the lower-dimensional original representation with tensor structure, which is why the relaxation from \textit{exact} to \textit{approximate} moment matching is needed. The proposed conditions aiming for approximate moment matching are presented in Section~\ref{subsec:prop-app-cond-input-tailored}.

\subsection{Moments and linear theory} \label{subsec:exact-moments}
Some basic theory of linear moment matching is recalled here for convenience. For further reading and proofs we refer to, e.g., \cite{book:antoulas2005}, \cite{phd:morGri97}, \cite{art:steadystate-mor-2010}, and references therein.

\begin{definition}[Moments] \label{def:gen-moment}
Given a univariate frequency representation $\Transfer{H}$ being $k$-times differentiable at $s_0 \in \mathbb{C}$, its $k$-th moment at $s_0$ is defined as
\begin{align*}
	\moma_{k} = \frac{(-1)^k}{k!} \frac{d^k}{ds^k} \Transfer{H}(s)_{|s=s_0}.
\end{align*}	
\end{definition}
The moments $\moma_k$ depend on the selected expansion frequency $s_0$. Note that we suppress this dependency in our notation in order to keep it shorter.

\begin{lemma}\label{lem:lin-moment}
Let $\Transfer{H}$ be a frequency representation, $\Transfer{H}(s) = \tensorbj{C} (s \tensorbj{E} - \tensorbj{A} )^{-1} \tensorbj{B}$. Let, for given $s_0$, $\tensorbj{A}_{s_0} = -s_0\tensorbj{E} + \tensorbj{A}  $ be nonsingular. Then, the $k$-th moment of $\Transfer{H}$ at $s_0$ reads 
\begin{align*}
	\moma_{k}  =  -  \tensorbj{C} \left[ \tensorbj{A}_{s_0}^{-1} \tensorbj{E} \right]^{k} \tensorbj{A}_{s_0}^{-1} \tensorbj{B}, \quad \text{for } k\geq 0.
\end{align*}
The moments can be determined as follows: 
Calculate $\vecbj{k}_i$, which are the moments of $s \mapsto(s \tensorbj{E} - \tensorbj{A} )^{-1} \tensorbj{B}$ at $s_0$, by the recursion
\begin{align*}
	i = 0: \quad \tensorbj{A}_{s_0} \vecbj{k}_{0} &= - \tensorbj{B} \\ 
	i > 0: \quad \tensorbj{A}_{s_0} \vecbj{k}_{i} &= \tensorbj{E} \vecbj{k}_{i-1}.
\end{align*}
Then set $\moma_{k} = \tensorbj{C} \vecbj{k}_{k}$.
\end{lemma}
For linear systems, reduced models fulfilling moment matching can be constructed by means of the following lemma.
\begin{lemma}\label{lem:lin-moment-matching}
Assume $\Transfer{H}(s) = \tensorbj{C} (s \tensorbj{E} - \tensorbj{A} )^{-1} \tensorbj{B}$ with $\tensorbj{E},\tensorbj{A} \in \mathbb{R}^{N,N}$. Let for given reduction basis $\tensorbj{V} \in \mathbb{R}^{N,n}$ the reduced system be defined as $\Transfer{H}_r(s) = \tensorbj{C}_r (s \tensorbj{E}_r - \tensorbj{A}_r )^{-1} \tensorbj{B}_r $
with $\tensorbj{E}_r = \tensorbj{V}^T \tensorbj{E} \tensorbj{V}$, $\tensorbj{A}_r = \tensorbj{V}^T \tensorbj{A} \tensorbj{V}$, $\tensorbj{B}_r = \tensorbj{V}^T \tensorbj{B}$ and $\tensorbj{C}_r = \tensorbj{C} \tensorbj{V}$. If, for prescribed $s_0$, the relation $\mathrm{span}\{ \vecbj{k}_0, \vecbj{k}_1,\ldots \vecbj{k}_k \} \subseteq \myim(\tensorbj{V})$
holds for $\vecbj{k}_i$, $0\leq i \leq k$ (as defined in Lemma \ref{lem:lin-moment}), then the (exact) moment matching condition
\begin{align*}
	\frac{d^i}{ds^i} \Transfer{H}(s)_{|s=s_0} = \frac{d^i}{ds^i} \Transfer{H}_r(s)_{|s=s_0}, \qquad i \leq k 
\end{align*}
is satisfied. We say that the moments of the full and the reduced models match (up to $k$-th order at $s_0$). Moreover, it holds 
\begin{align*}
	\vecbj{k}_i = \tensorbj{V} \vecbj{k}_{r,i}, \qquad i \leq k,
\end{align*}
where $\vecbj{k}_{r,i}$ is recursively defined with $\tensorbj{A}_{r,s_0} = -s_0\tensorbj{E}_r + \tensorbj{A}_r $ as 
\begin{align*}
	i = 0: \quad \tensorbj{A}_{r,s_0} \vecbj{k}_{r,0} &= - \tensorbj{B}_r\\
	i > 0: \quad \,\tensorbj{A}_{r,s_0} \vecbj{k}_{r,i} &= \tensorbj{E}_r \vecbj{k}_{r,i-1}.
\end{align*}
\end{lemma}

\subsection{Reduced signal generator driven system} \label{subsec-siggendriv_r}

In this subsection we clarify our notion of a reduced signal generator driven system and its usage. We start by stating the basic result behind the commuting diagram sketched in Fig.~\ref{fig:input-tailored-sketch}.

\begin{lemma} \label{lem:commuting-sg-proj}
Let a quadratic-bilinear system $\statemap$, a signal generator $\siggen$,
\begin{align*}
	\statemap:& \quad \tensorbj{E} \dot{\vecbj{x}} = \tensorbj{A} \vecbj{x} + \tensorbj{G}\vecbj{x}\supertens{2} + \tensorbj{D} (\vecbj{x} \otimes \vecbj{u})  +  \tensorbj{B} \vecbj{u} , \quad &&
	\vecbj{x}(0) = \vecbj{x}_0,\\
	\siggen:& \quad \,\,\,\, \vecbj{u} = \tensorbj{C}_{\sigstatescal} \sigstate, \qquad \dot{\sigstate} = \tensorbj{A}_{\sigstatescal} \sigstate + \tensorbj{G}_{\sigstatescal} \sigstate\supertens{2}, \qquad && \,\sigstate(0) = \sigstate_{0}, \,\,
\end{align*}
and the associated signal generator driven system $\statemapAut$, as in Definition~\ref{def:sg-ode-casc}, be given. Let furthermore, for given reduction basis $\tensorbj{V}\in \mathbb{R}^{N,n}$, $n\ll N$, the reduced state matrices be defined as $\tensorbj{E}_r = \tensorbj{V}^T \tensorbj{E} \tensorbj{V}$, $\tensorbj{A}_r = \tensorbj{V}^T \tensorbj{A} \tensorbj{V}$, $\tensorbj{G}_r = \tensorbj{V}^T \tensorbj{G} (\tensorbj{V} \otimes \tensorbj{V})$, $\tensorbj{D}_r = \tensorbj{V}^T \tensorbj{D} (\tensorbj{V} \otimes \tensorbj{I}_p)$ and  $\tensorbj{B}_r = \tensorbj{V}^T \tensorbj{B}$.
Let $\tensorbj{Q}_r$ be the constant matrix such that
\begin{align*}
	\tensorbj{Q}_r \begin{mymatrix}
		\bar{\vecbj{x}} \\
		\bar{\sigstate}
	\end{mymatrix}\supertens{2} = 
	\begin{mymatrix}
		\bar{\vecbj{x}}\supertens{2} \\
		\bar{\vecbj{x}} \otimes \bar{\sigstate} \\
		\bar{\sigstate}\supertens{2}
	\end{mymatrix}  \nonumber \qquad  \text{for arbitrary  } \bar{\vecbj{x}} \in \mathbb{R}^n,\,  \bar{\sigstate}\in \mathbb{R}^q.
\end{align*}	 
Introducing the reduced system as
\begin{align*}
 \statemap_r:& \qquad	\tensorbj{E}_r \dot{\vecbj{x}}_r = \tensorbj{A}_r \vecbj{x}_r + \tensorbj{G}_r \vecbj{x}_r\supertens{2} + \tensorbj{D}_r (\vecbj{x}_r \otimes \vecbj{u}) + \tensorbj{B}_r \vecbj{u},  && \vecbj{x}_r(0)= \tensorbj{V}^T \vecbj{x}_0,
\end{align*}
and setting up the signal generator driven system for $\statemap_r$ and $\siggen$ gives
\begin{align*}
\parbox{0.6cm}{{\vspace{0.4cm}$\statemapAut_r:$}} 
\,\,& \qquad	\,\, \tensorcj{E}_{r} \dot{\veccj{w}}_r = \tensorcj{A}_{r} \veccj{w}_r + \tensorcj{G}_{r} \veccj{w}_r\supertens{2},  && \veccj{w}_r(0) = \veccj{b}_r,  \\
& 	\quad \qquad \quad   \vecbj{x}_r = \tensorcj{P}_{x_r} \, \veccj{w}_r. & &
\end{align*}
with
\begin{align*}
	\tensorcj{E}_r &= \begin{mymatrix}
		\tensorbj{E}_r & \\
					& \tensorbj{I}_q
	\end{mymatrix}
	, \quad 
	\tensorcj{A}_r =
		\begin{mymatrix}
		\tensorbj{A}_r & \tensorbj{B}_r \tensorbj{C}_{\sigstatescal}  \\
					& \tensorbj{A}_{\sigstatescal}
	\end{mymatrix},
	\quad \tensorcj{G}_r =
		 \begin{mymatrix}
		\tensorbj{G}_r & \tensorbj{D}_r(\tensorbj{I}_n \otimes\tensorbj{C}_{\sigstatescal}) & \\
					& 													  & \tensorbj{G}_{\sigstatescal}
	\end{mymatrix}	\tensorbj{Q}_r\\
	\veccj{b}_r &= 
	\begin{mymatrix}
		\vecbj{V}^T \vecbj{x}_{0} \\
		 \sigstate_0
	\end{mymatrix},
	\quad \, \tensorcj{P}_{x_r} = [\tensorbj{I}_n, \, \tensorbj{0}].
\end{align*}
Projecting the realization of $\statemapAut$ as
\begin{align*}
& \tensorcj{E}_r = \tensorcj{V}^T \tensorcj{E} \tensorcj{V}, \qquad \tensorcj{A}_r = \tensorcj{V}^T \tensorcj{A} \tensorcj{V}, \qquad \tensorcj{G}_r = \tensorcj{V}^T \tensorcj{G} (\tensorcj{V} \otimes \tensorcj{V}) 
	 \qquad \veccj{b}_r = \tensorcj{V}^T \veccj{b}, \\
& \text{with reduction basis } \,\tensorcj{V} =
	\begin{mymatrix}
		\tensorbj{V} & \\
					 & \tensorbj{I}_q		
	\end{mymatrix}
\end{align*}
leads to the same reduced signal generator driven system $\statemapAut_r$.
\end{lemma}
The lemma is quite obvious but nonetheless of high importance for our approach. The input-tailored frequency representations $\breve{\Transfer{X}}_{r,i}$ of $\statemapAut_r$ are accordingly obtained by $\breve{\Transfer{X}}_{r,i}(s) = \tensorcj{P}_{x_r} \, \breve{\Transfer{W}}_{r,i}(s)$ for $s \in \mathbb{C}$
with $\breve{\Transfer{W}}_{r,i}$ being the frequency representation of the variational expansion term $\veccj{w}_{r,i}$ (cf.\ Definition \ref{def:input-tailor-var-exp}). For prescribed order $L_i$, index $\bar{i} \in \mathbb{N}$ and frequency $s_0 \in \mathbb{C}$, the \textit{approximate} moment matching conditions we require on the reduced model $\statemap_r$ are
\begin{align} \label{eq:approx-abstract-cond-assoc-tfx}
	\tensorbj{V} \frac{d^k}{ds^k} \breve{\Transfer{X}}_{r,i}(s)_{|s=s_0} \stackrel{!}{\approx} \frac{d^k}{ds^k} \breve{\Transfer{X}}_{i}(s)_{|s=s_0} \qquad \text{for } k \leq L_i, \quad i \leq \bar{i}.
\end{align}

\begin{remark}[Extracting reduction basis from extended problem] \label{rem:connection-xr-wr}
Note that the signal generator itself is not reduced in the construction of Lemma~\ref{lem:commuting-sg-proj}. This is reflected in the block structure of $\tensorcj{V}$ with a unit matrix block $\tensorbj{I}_q$. The lemma shows that projection and driving by a signal generator commute. Thus, the input-tailored moment matching \eqref{eq:approx-abstract-cond-assoc-tfx} can be approached in a two-step procedure:
\begin{itemize}
\item Find basis $\tensorcj{V}$ such that
\begin{align} \label{eq:approx-abstract-cond-assoc-tf}
	\tensorcj{V} \frac{d^k}{ds^k} \breve{\Transfer{W}}_{r,i}(s)_{|s=s_0} \stackrel{!}{\approx} \frac{d^k}{ds^k} \breve{\Transfer{W}}_{i}(s)_{|s=s_0} \qquad \text{for } k \leq L_i, \quad i \leq \bar{i}
\end{align}
for prescribed $L_i$, $\bar{i}, s_0$, where $\breve{\Transfer{W}}_{i}$, $\breve{\Transfer{W}}_{r,i}$ are as in Definition~\ref{def:input-tailor-var-exp} given $\statemapAut$, $\statemapAut_r$.
\item  Extract the basis $\tensorbj{V}$ from $\tensorcj{V}$.
\end{itemize}
\end{remark}

\subsection{Input-tailored moments and projection} \label{subsec:input-tailored-mom-lin}
Up to now, the input-tailored moment matching problem has been tracked back to the extended problem  (Remark~\ref{rem:connection-xr-wr}), and it has been shown that the reduced signal generator driven system $\statemapAut_r$ can be seen as the projection of $\statemapAut$, Lemma~\ref{lem:commuting-sg-proj}. In this subsection the actual structure of the extended problem \eqref{eq:approx-abstract-cond-assoc-tf} is investigated and structural properties are highlighted. 

\begin{lemma}[Reduced associated frequency representation]\label{lem:red-w2r-proj}
Given the full order signal generator $\statemapAut$ and its reduced counterpart $\statemapAut_r$ as in Lemma~\ref{lem:commuting-sg-proj}, the reduced associated frequency representation $\breve{\Transfer{W}}_{r,2}$ is the Galerkin-projection of $\breve{\Transfer{W}}_{2}$ written in its high-dimensional linear representation of Lemma~\ref{lem:assoc-trafo-repres}, i.e.,
\begin{align*}
\breve{\Transfer{W}}_{r,2}(s) &= \breve{\tensorcj{C}}_{r,2} \left(s \breve{\tensorcj{E}}_{r,2} - \breve{\tensorcj{A}}_{r,2} \right)^{-1} \breve{\veccj{b}}_{r,2}\\
& \text{with} \quad
	\breve{\tensorcj{E}}_{r,2} = \breve{\tensorcj{V}}_2^T \breve{\tensorcj{E}}_2 \breve{\tensorcj{V}}_2, \quad \breve{\tensorcj{A}}_{r,2} = \breve{\tensorcj{V}}_2^T \breve{\tensorcj{A}}_2 \breve{\tensorcj{V}}_2,
	 \quad \breve{\veccj{b}}_{r,2} = \breve{\tensorcj{V}}_2^T \breve{\veccj{b}}_2, \quad \breve{\tensorcj{C}}_{r,2} = \breve{\tensorcj{C}}_2 \breve{\tensorcj{V}}_2, \\
	 & \text{and} \quad
		\breve{\tensorcj{V}}_2 = 
	\begin{mymatrix}
		\tensorcj{V} & \\
		& \tensorcj{V}\supertens{2}
	\end{mymatrix}.
\end{align*}
\end{lemma}
The proof is straightforward using the properties of the Kronecker product. From a theoretical point of view, the moment matching problem \eqref{eq:approx-abstract-cond-assoc-tf} can be embedded into the linear theory when changing into a high-dimensional space. The inherent tensor-structure of the problem is handed over via the reduction basis $\breve{\tensorcj{V}}_2$. Our method makes use of this special cascade- and tensor-structure that is also present in the moments, as we show in the upcoming.

\begin{lemma} \label{lem:As0-tensor}
	Given $\tensorcj{A}_{s_0} = - s_0 \tensorcj{E} + \tensorcj{A}$ for $s_0 \in \mathbb{C}$, $i >0$ and quadratic matrices $\tensorcj{E},\tensorcj{A}$, the following relation holds $\circledmy{{\relsize{-1.5}i}}_{\tensorcj{E}} \tensorcj{A}_{s_0/i} = - s_0 \tensorcj{E}\supertens{i} + \circledmy{{\relsize{-1.5}i}}_{\tensorcj{E}} \tensorcj{A}.$
\end{lemma}
\begin{proof}
For $0\leq k,m \leq i-1$ with $k+m+1=i$ it holds
\begin{align*}
	\tensorcj{E}\supertens{k}  \otimes \left(- \frac{s_0}{i}\tensorcj{E} + \tensorcj{A} \right) \otimes  \tensorcj{E}\supertens{m} = -  \frac{s_0}{i} \tensorcj{E} \supertens{i}+  \tensorcj{E}\supertens{k}  \otimes \tensorcj{A}  \otimes \tensorcj{E}\supertens{m} .
\end{align*}
Since $\circledmy{{\relsize{-1.5}i}}_{\tensorcj{E}} \tensorcj{A}_{s_0/i}$ can be written as sum of $i$ such expressions with $k= 0, \ldots, i-1$, and $m=i-k-1$, the lemma follows.
\end{proof}
A recursion formula for the moments $\moma^{(2)}_i$ of $\breve{\Transfer{W}}_2$ can now be stated (cf.\ Theorem~\ref{lem:mom-assoc-tf3} for $\breve{\Transfer{W}}_3$). Note that the super-index $^{(j)}$ is used throughout to indicate the correspondence to the $j$-th frequency representation $\breve{\Transfer{W}}_j$, $j=2,3$. 

\begin{theorem}[Extended input-tailored moments] \label{lem:mom-assoc-tf2}
	Assume that the requirements of Theorem~\ref{theor:Pxu_quadratic} and Lemma~\ref{lem:assoc-trafo-repres} hold and let for given $s_0 \in \mathbb{C}$ the matrix $\tensorcj{A}_{s_0} = - s_0 \tensorcj{E} + \tensorcj{A}$ be nonsingular. Then the moments $\moma^{(2)}_i$ of $\breve{\Transfer{W}}_2$ at $s_0$ are characterized by the recursion:
		\begin{eqnarray*}
	i = 0:\,\, &	\circledmy{2}_{\tensorcj{E}} \tensorcj{A}_{s_0/2} \momb^{(2)}_{0} &=\, -  \veccj{b} \supertens{2}\\
			&  \qquad \tensorcj{A}_{s_0} \moma^{(2)}_{0} &=\, -\tensorcj{G}  \momb^{(2)}_{0} \\
	i >0:\,\,
		&\circledmy{2}_{\tensorcj{E}} \tensorcj{A}_{s_0/2} \, \momb^{(2)}_i &=\,  \tensorcj{E}\supertens{2} \, \momb^{(2)}_{i-1}\\
		& \qquad \tensorcj{A}_{s_0} \, \moma^{(2)}_i &=\,  \tensorcj{E}  \moma^{(2)}_{i-1} -\tensorcj{G}  \momb^{(2)}_{i}.
	\end{eqnarray*}
Moreover, $\tensorbj{k}_i^{(2)} = [\moma^{(2)}_i;\momb^{(2)}_i]$ are the moments of $s \mapsto\left(s \breve{\tensorcj{E}}_2 - \breve{\tensorcj{A}}_2 \right)^{-1} \breve{\veccj{b}}_2$ at $s_0$.
\end{theorem}	
\begin{proof}
The representation of Lemma \ref{lem:assoc-trafo-repres} for $\breve{\Transfer{W}}_2$ is a linear state representation. Therefore, following Lemma \ref{lem:lin-moment}, the factors $ \tensorbj{k}_i^{(2)}$, which we recursively define by
\begin{align*}
	i = 0: \quad (-s_0 \breve{\tensorcj{E}}_2 + \breve{\tensorcj{A}}_2) \vecbj{k}_{0}^{(2)} &= -  \breve{\veccj{b}}_2\\
	i > 0: \quad (-s_0 \breve{\tensorcj{E}}_2 + \breve{\tensorcj{A}}_2)  \vecbj{k}_{i}^{(2)} &= \breve{\tensorcj{E}}_2 \, \vecbj{k}_{i-1}^{(2)},
\end{align*}
are the moments of $s \mapsto\left(s \breve{\tensorcj{E}}_2 - \breve{\tensorcj{A}}_2 \right)^{-1} \breve{\veccj{b}}_2$ at $s_0$. Let us introduce the following notation for the upper and lower blocks
\begin{align*}
	\tensorbj{k}_i^{(2)} = \begin{mymatrix}
			\moma^{(2)}_i \\
			\momb^{(2)}_i
		\end{mymatrix},  \qquad  \text{where $\moma^{(2)}_i \in \mathbb{R}^{M}, \quad \momb^{(2)}_i \in \mathbb{R}^{M^2}$}.
\end{align*}
Then these blocks fulfill for $i>0$
\begin{align*}
	(- s_0 \tensorcj{E} + \tensorcj{A}) \moma^{(2)}_i + \tensorcj{G} \momb^{(2)}_i &= \tensorcj{E} \moma^{(2)}_{i-1} \\
	\left( - s_0 \tensorcj{E}\supertens{2} + \circledmy{2}_{\tensorcj{E}} \, \tensorcj{A}  \right) \momb^{(2)}_i &=  \tensorcj{E}\supertens{2} \momb^{(2)}_{i-1}.
\end{align*}
Using Lemma \ref{lem:As0-tensor}, we get the recursive expression for $\moma^{(2)}_i$ for $i>0$. The initial step $i=0$ follows similarly. In fact, $\moma^{(2)}_i$ is the $i$-th moment of $\breve{\Transfer{W}}_2$ at $s_0$, as it equals $\breve{\tensorcj{C}}_2 \tensorbj{k}_i^{(2)}$, which corresponds to the expression for the moments of Lemma \ref{lem:lin-moment}.
\end{proof}
According to the linear theory, exact moment matching requires 
\begin{align}\label{eq:exact}
\vecbj{k}^{(2)}_i\in \text{ image}(\breve{\tensorcj{V}}_2).
\end{align} 
This corresponds to a condition in a $(N+q)^2+(N+q)$-dimensional space. However, this condition cannot be fulfilled exactly because of the specific form our reduction basis has.

\subsection{Proposed approximation conditions}\label{subsec:prop-app-cond-input-tailored}
We propose an approximate moment matching that accounts for the special tensor structure of the problem. Considering the reduction basis for $\breve{\Transfer{W}}_2$
\begin{align*}
		\breve{\tensorcj{V}}_2 = 
	\begin{mymatrix}
		\tensorcj{V}_2 & \\
		& \tensorcj{V}\supertens{2}_2
	\end{mymatrix},
	\end{align*}
we solve the following splitted problem: Find $\tensorcj{V}_2$ such that it holds
\begin{subequations} \label{eq:tf2-approxspace}
\begin{align} 
	||(\tensorbj{I}_{N+q} - \tensorcj{V}_2 \tensorcj{V}_2^T )\moma^{(2)}_{i}|| / || \moma^{(2)}_{i}|| \text{ \textit{ small} for } i=0,1,\ldots L \label{eq:tf2-approxspace-a} \\
	 ||(\tensorbj{I}_{(N+q)^2} - \tensorcj{V}\supertens{2}_2 (\tensorcj{V}_2\supertens{2})^T ) \momb^{(2)}_{i}|| / || \momb^{(2)}_{i}|| \text{ \textit{ small} for } i=0,1,\ldots L \label{eq:tf2-approxspace-b} 
\end{align}
\end{subequations}
for $\moma^{(2)}_i$, $\momb^{(2)}_i$ from Lemma~\ref{lem:mom-assoc-tf2}. 
This aims for small projection errors
\begin{align*}
	\left(\tensorbj{I} - \breve{\tensorcj{V}}_2 \breve{\tensorcj{V}}_2^T \right)  \vecbj{k}^{(2)}_i  \qquad \text{ with } \vecbj{k}^{(2)}_i =[\moma^{(2)}_i;\momb^{(2)}_i],
\end{align*}
which is a relaxation of the exact moment matching in \eqref{eq:exact}.

In the assembly of the global reduction basis $\tensorcj{V}$ that corresponds to all considered frequency representations $\breve{\Transfer{W}}_i$, $i\leq \bar i$, cf.\ \eqref{eq:approx-abstract-cond-assoc-tf}, we provide a block structure of the form
\begin{align*}
 \tensorcj{V} =
	 \begin{mymatrix}
	 	\tensorbj{V} & \\
	 				 & \tensorbj{I}_q
	 \end{mymatrix}.
 \end{align*}
This reflects that the signal generator itself is not reduced and gives the desired reduction basis $\tensorbj{V}$ of the original system.

\begin{remark} \label{rem-compare-to-assoctraf}
Let us stress the difference to former work on model reduction using univariate frequency representations for nonlinear systems. Comparing our approach with the one from \cite{inproc:fast-nonlin-mor-assoc-trafo}, \cite{art:mor-associated-transform-2016} there are, besides the more rigorous treatment of the variational expansion (cf.\ Remark~\ref{rem:formal-volterra-dirac}), three major differences: The first and most important one is that our analysis reveals an additional tensor-structured approximation condition \eqref{eq:tf2-approxspace-b} to naturally appear when aiming for approximate moment matching. Such a condition is not present in the former approach. Second, our framework using the concept of signal generator driven systems enables us to consider a larger class of input scenarios within the process. And finally, the inherent cascade- and sparse-tensor-structure has not been exploited in the former algorithmic implementation. It will be seen in Section~\ref{sec:approx-mom-match} that the appearing tensor-structured problems can be formulated as Lyapunov-type equations with 'sparse right hand sides'. We deal with them using recently proposed low-rank solvers from literature, which is known to save memory- and time-effort by orders of magnitude, cf.~\cite{SaaKB16-mmess-1.0.1}, \cite{art:Simoncini_anew}, \cite{art:kressner-kryl-tensor}.
\end{remark}

\section{Numerical realization of approximate input-tailored moment matching} \label{sec:approx-mom-match}

In this section we present and discuss the algorithms for the numerical realization of our input-tailored moment matching method.

\subsection{Low-rank calculations of input tailored moments} \label{subsec:hope-eff-method}

The main part of the numerics consists in constructing the subspace for basis $\tensorcj{V}$ such that \eqref{eq:tf2-approxspace} holds. Clearly, it is easy to construct a basis matrix $\tensorcj{V}$ fulfilling \eqref{eq:tf2-approxspace-a} exactly, namely just use the matrix composed of the moments $\moma^{(2)}_{i}$ itself. The question remains, why a low-rank basis fulfilling \eqref{eq:tf2-approxspace-b} should exist. Let us herefore look at the zeroth auxiliary moment $\momb^{(2)}_0$ around $s_0$. It reads
\begin{align*}
	\left[ \tensorcj{E} \otimes \tensorcj{A}_{s_0/2} + \tensorcj{A}_{s_0/2} \otimes \tensorcj{E} \right]  \momb^{(2)}_0 + \veccj{b} \supertens{2} = \tensorbj{0},
\end{align*}
which is the well-known Lyapunov equation, written in tensor notation, with a sparse 'right hand side' $\veccj{b} \supertens{2}$. Low-rank solutions for this kind of equations exist under reasonable conditions \cite{art:kressner-kryl-tensor}, \cite{art:Simoncini_anew}, \cite{art:H2mor-bilinear-breiten2012}, and take the form
\begin{align} \label{eq:low-rank-approx}
 \sum_{k=1}^{n_i} \vecbj{z}^k_i \otimes \vecbj{z}^k_i \approx {\momb}^{(2)}_i  \quad \text{for small }n_i.
\end{align}
For the higher order terms, e.g., $\momb^{(2)}_1$, we suggest to follow up the iteration with the new  sparse 'right hand side' $\tensorcj{E}\supertens{2}\momb^{(2)}_0$, i.e., the low-rank approximation from the former step, and so on. By that, we do not only have a strategy to efficiently approximate $\momb^{(2)}_i$ and $\moma^{(2)}_i$ up to a certain extend but also a candidate for a low-rank basis, namely the span over all $\vecbj{z}^k_i$. The upcoming Algorithm~\ref{alg:mmh2-simple} summarizes our approach aiming towards \eqref{eq:tf2-approxspace}.

Note that the moments involved are the ones for the signal generator driven system $\statemapAut$. Albeit the reduction basis $\tensorbj{V}$  is constructed for the original system $\statemap$. Thus, the selection matrix $\tensorcj{P}_x: \veccj{w} \mapsto \vecbj{x}$ appears here.

\begin{myalgorithm}[Moment-matching-bases for $\breve{\Transfer{X}}_2$]\text{$ $} \label{alg:mmh2-simple}\\
INPUT:
\begin{itemize}
\item Realization matrices of signal generator driven system $\statemapAut$ (cf.\ Definition~\ref{def:sg-ode-casc}): $\tensorcj{E}$, $\tensorcj{A}$, $\tensorcj{G}$, $\veccj{b}$
\item Dimension of state variable $N$; Dimension of signal generator: $q$
\item Expansion frequencies: $(s_1,s_2,\ldots,s_\mu)$; Number of moments: $(L_1,L_2,\ldots,L_\mu)$
\item Tolerance for low-rank approximations: \text{tol}
\item Basis for space not considered in low-rank approximation: $\tensorbj{V}_\perp$
\end{itemize}
OUTPUT: Reduction bases: $\tensorbj{V}_a$, $\tensorbj{V}_b$. 
\begin{enumerate}
	\item Set $\tensorcj{P}_x = [\tensorbj{I}_N , \ \tensorbj{0}_{N,q}]$. 
	\item for $j = 1,\ldots \mu$
		\begin{enumerate}[a)]
				\item Set $s_0 := s_j$ and $L := L_j$.
				\item \label{it:alg:mmh2-simple-zik} Calculate low-rank factors $\vecbj{z}_i^k$ for $k = 1, \ldots n_i$, $i = 0, \ldots L-1$, see \eqref{eq:low-rank-approx}, i.e.,
				\vspace{-0.3cm}
					\begin{align*}
						\vecbj{z}_i^k \text{ with: } \sum_{k=1}^{n_i} \vecbj{z}_i^k \otimes \vecbj{z}_i^k \approx \left( \left(\circledmy{2}_{\tensorcj{E}} \tensorcj{A}_{s_0/2}\right)^{-1} \, \tensorcj{E}\supertens{2} \right)^{i} \left( \circledmy{2}_{\tensorcj{E}} \tensorcj{A}_{s_0/2}\right)^{-1} \veccj{b} \supertens{2}.
					\end{align*}
				\item Gather all $(\tensorcj{P}_x \vecbj{z}_i^k)$ in \,$\tensorbj{Z}_{s_j}:= \tensorcj{P}_x [\vecbj{z}_0^1,\vecbj{z}_0^2,\ldots \vecbj{z}_0^{n_0}, \vecbj{z}_1^{1},\ldots \vecbj{z}_1^{n_1},\ldots \vecbj{z}_{L-1}^{n_{L-1}}]$.
		\end{enumerate}
	\item[] endfor
		\item Gather all $\tensorbj{Z}_{s_j}$  in \, $\tensorbj{Z}:= [\tensorbj{Z}_{s_1},\tensorbj{Z}_{s_2},\ldots \tensorbj{Z}_{s_\mu}]$.
		\item for $j = 1,\ldots \mu$	
		\begin{enumerate}[a)]
			\item Set $s_0 := s_j$ and $L := L_j$.
			\item Calculate $\moma^{(2)}_i$ for $s_0$ from Lemma \ref{lem:mom-assoc-tf2} (using the low-rank approximations on $\momb^{(2)}_i$ from Step \eqref{it:alg:mmh2-simple-zik})
			\item Gather all $(\tensorcj{P}_x \moma^{(2)}_i)$ in \, $\tensorbj{M}_{s_j}:= \tensorcj{P}_x [\vecbj{m}^{(2)}_0,\vecbj{m}^{(2)}_1,\ldots \vecbj{m}^{(2)}_{L-1}]$.
		\end{enumerate}
	endfor
	\item Construct $\tensorbj{V}_a$ as orthogonal basis of $[\tensorbj{M}_{s_1} , \ldots ,\tensorbj{M}_{s_\mu}]$.
	\item \label{it:alg:mmh2-orthproj} Define $\tensorbj{P}_\perp$ as orthogonal projection onto the orthogonal complement of span of $[\tensorbj{V}_a,\tensorbj{V}_\perp]$. Then calculate $\tensorbj{V}_b$ with column span defined by all left-singular vectors of  $(\tensorbj{P}_\perp \tensorbj{Z})$ with singular value bigger than \text{tol}.
\end{enumerate}
\end{myalgorithm}
In terms of numerical calculation, the most delicate step is the construction of the low-rank factors $\vecbj{z}_i^k$. Note that for each $i$ in Step \eqref{it:alg:mmh2-simple-zik} we actually need to construct a low-rank solution of a Lyapunov equation. The projection step with $\tensorbj{P}_\perp$ removes components of the dominant space already present in the former constructed bases and therefore allows for lower-order truncation in Step \eqref{it:alg:mmh2-orthproj}.

\subsection{Constructing the full reduction basis}

In this subsection we conclude our approach for the construction of a reduced model, which aims at approximate moment matching of the input-tailored frequency representations $\breve{\Transfer{X}}_1$, $\breve{\Transfer{X}}_2$ from Definition \ref{def:input-tailor-var-exp}.

For the basis construction concerning $\breve{\Transfer{X}}_1$, the signal generator does not need to be considered. This is because $\breve{\Transfer{W}}_1$ can be factorized as
\begin{align*}
	\breve{\Transfer{W}}_1(s) = \left[(s \tensorbj{E} - \tensorbj{A} )^{-1} \tensorbj{B} \right] \tensorbj{C}_{\sigstatescal}  (s \tensorbj{I}_q - \tensorbj{A}_{\sigstatescal} )^{-1} \sigstate_0,
\end{align*}
i.e., into the standard linear transfer function  $(s \tensorbj{E} - \tensorbj{A} )^{-1} \tensorbj{B}$ and the signal generator.
As discussed in Section \ref{sec:towards-input-tailored-syst}, the signal generator is not reduced, and thus moment matching of the linear transfer function automatically imposes moment matching on $\breve{\Transfer{W}}_1$. Concluding, the following algorithm for the construction of a reduced model is proposed.

\begin{myalgorithm}[Input-tailored approximate moment matching]\label{alg:mm-complete}
\\ INPUT:
\begin{itemize}
\item Realization matrices of the quadratic-bilinear dynamical system $\statemap$ to reduce: $\tensorbj{E}$, $\tensorbj{A}$, $\tensorbj{G}$, $\tensorbj{D}$, $\tensorbj{B}$, $\tensorbj{C}$
\item Realization matrices of the signal generator $\siggen$: $\tensorbj{A}_{\sigstatescal}$, $\tensorbj{G}_{\sigstatescal}$, $\tensorbj{C}_{\sigstatescal}$
\item Initial value vectors: $\vecbj{x_0}, \sigstate_0$
\item Concerning $\breve{\Transfer{X}}_2$: Expansion frequencies: $(s_1,s_2,\ldots,s_\mu)$; Number of moments to match: $(L_1,L_2,\ldots,L_\mu)$; Tolerance for low-rank approximations in Algorithm~\ref{alg:mmh2-simple}: \text{tol}
\item Concerning $\breve{\Transfer{X}}_1$: Expansion frequencies: $(\tilde{s}_1,\tilde{s}_2,\ldots,\tilde{s}_\nu)$; Number of moments to match: $(\tilde{L}_1,\tilde{L}_2,\ldots,\tilde{L}_\nu)$
\end{itemize}
OUTPUT: Reduced realization: $\tensorbj{E}_r$, $\tensorbj{A}_r$, $\tensorbj{G}_r$, $\tensorbj{D}_r$, $\tensorbj{B}_r$, $\tensorbj{C}_r$. 

\begin{enumerate}
	\item \label{it:h1-mm-step} Construct reduction basis $\tensorbj{V}_{1}$ for $\breve{\Transfer{X}}_1$ as orthonormal basis\\for the union of the Krylov spaces
	$\mathcal{K}_{\tilde{L}_j}(\tensorbj{A}_{\tilde{s}_j}^{-1}\tensorbj{E}, \tensorbj{A}_{\tilde{s}_j}^{-1} \vecbj{b})$ for $j = 1,\ldots \nu$.
	\item Construct realization for signal generator driven system $\statemapAut$ (Definition \ref{def:sg-ode-casc}): $\tensorcj{E}$, $\tensorcj{A}$, $\tensorcj{G}$, $\veccj{b} $
	\item Construct reduction bases $\tensorbj{V}_a, \, \tensorbj{V}_b$ for $\breve{\Transfer{X}}_2$ by Algorithm \ref{alg:mmh2-simple}	for \\
	frequencies $(s_1,s_2,\ldots,s_\mu)$, number of moments $(L_1,L_2,\ldots,L_\mu)$, tolerance \text{tol} and $\tensorbj{V}_\perp = \tensorbj{V}_{1}$.
	\item Construct $\tensorbj{V}$ as orthogonal basis of span of \, $[\tensorbj{V}_{a} , \tensorbj{V}_b,\tensorbj{V}_1]$.
	\item Calculate reduced state representation as $\tensorbj{E}_r = \tensorbj{V}^T \tensorbj{E} \tensorbj{V}$, $\tensorbj{A}_r = \tensorbj{V}^T \tensorbj{A} \tensorbj{V}$, $\tensorbj{G}_r = \tensorbj{V}^T \tensorbj{G} (\tensorbj{V} \otimes \tensorbj{V})$, $\tensorbj{D}_r = \tensorbj{V}^T \tensorbj{D} (\tensorbj{V} \otimes \tensorbj{I}_p)$, $\tensorbj{B}_r = \tensorbj{V}^T \tensorbj{B}$, $\tensorbj{C}_r = \tensorbj{C} \tensorbj{V}$.
\end{enumerate}
\end{myalgorithm}
For Step~\eqref{it:h1-mm-step} in Algorithm~\ref{alg:mm-complete} we just use the standard Krylov method as in \cite{phd:morGri97}, \cite{book:antoulas2005}. Note furthermore that in the calculation of $\tensorbj{G}_r$ it is advisable to avoid the memory-demanding explicit calculation of $\tensorbj{V} \otimes \tensorbj{V}$, see \cite{morBre13}, which we also do.

\begin{remark}
Algorithm~\ref{alg:mmh2-simple} is only assumed to be stable, if the orders $L_j$ of matched moments are chosen moderate. This is, because we actually seek for a special so-called Krylov space. For matrices $\tensorbj{M}$, $\tensorbj{L}$ of appropriate dimension and $L \in \mathbb{N}$ the Krylov space is defined as
\begin{align*}
	\mathcal{K}_{L}(\tensorbj{M}, \tensorbj{L}) := \text{span}\left\{ \left[ \tensorbj{L},\, \tensorbj{M} \tensorbj{L},\ldots,\, \tensorbj{M}^{L-1} \tensorbj{L}\right] \right\}.
\end{align*} 
Step \eqref{it:alg:mmh2-simple-zik}, thought of in $\mathbb{R}^{N^2}$, consists of constructing the Krylov space
\begin{align*}
\mathcal{K}_{L} \left( \left(\circledmy{2}_{\tensorcj{E}} \tensorcj{A}_{s_0/2}\right)^{-1}   \, \tensorcj{E}\supertens{2} ,\, \left(\circledmy{2}_{\tensorcj{E}} \tensorcj{A}_{s_0/2}\right)^{-1}   \veccj{b} \supertens{2}  \right)
\end{align*}
without any orthogonalization between the iteration. This is known to be unstable for high orders, see, e.g., \cite{phd:morGri97}, \cite{book:antoulas2005}. However, orthogonalization in $\mathbb{R}^{N^2}$ destroys our tensor structure. It is possible to recover a low-rank tensor structure by additional truncation, but this goes with further approximation errors \cite{book:tensor-num-methods-quantum}. Therefore, we recommend to match the moments at several frequencies $s_i$ rather than for high-order moments as it is also usual practice for linear moment matching.
\end{remark}

\section{Handling non-standard input dependencies} \label{sec:ext-var-appr-app}

In practical applications the state equation $\statemap$ to reduce may take a more general form as in \eqref{eq:qldae-a}, e.g.,
\begin{align*} 
	\tensorbj{E} \dot{\vecbj{x}} &= \tensorbj{A} \vecbj{x} + \tensorbj{G} \vecbj{x}\supertens{2} + \tensorbj{D} (\vecbj{x} \otimes \vecbj{u}) + \tensorbj{B} \vecbj{u} + \vecbj{K}(\vecbj{u})
\end{align*}
with $\vecbj{K}(\vecbj{u})$ describing input dependencies not affine-linear in $\vecbj{u}$. For example, quad\-ratic terms in the inputs can come from boundary control terms, when systems with quadratic nonlinearities are discretized, as shown for the Burgers' equation in Section~\ref{sec:burgers-eq}. Also time derivatives in the input can appear, when the state equation $\statemap$ originates from an index-reduced differential-algebraic equation \cite{book:dae-mehr}, \cite{book:dae-projector-based}. The usual work-around in system-theoretic model reduction is to introduce artificial augmented inputs for all non-standard terms. Obviously, this enlarges the input and ignores known input-structure, which leads to worse results in model reduction. 

\subsection{Extension of input-tailored method} 
Our input-tailored approach can incorporate a large class of input-relations directly, as we discuss for some cases in the following.
\paragraph{Input map with quadratic term and/or time derivative} For $\vecbj{K}(\vecbj{u}) = \tensorbj{G}_u \vecbj{u}\supertens{2} +  \tensorbj{B}_p \dot{\vecbj{u}}$ our signal generator driven system, and with that the core of our approach, generalizes as follows.

\begin{definition}[Generalization of Definition \ref{def:sg-ode-casc}, Signal generator driven system] \label{def:sg-ode-casc2}
Let a system $\statemap$ with an input $\vecbj{u}$ described by a signal generator $\siggen$ (as in \eqref{eq:siggen}) be given as
\begin{align*}
\statemap:& \quad 	\tensorbj{E} \dot{\vecbj{x}} = \tensorbj{A} \vecbj{x} + \tensorbj{G}\vecbj{x}\supertens{2} + \tensorbj{D}  (\vecbj{x} \otimes \vecbj{u})  +  \tensorbj{B} \vecbj{u} + \tensorbj{G}_u \vecbj{u}\supertens{2} +  \tensorbj{B}_p \dot{\vecbj{u}},  &&
	\vecbj{x}(0) = \vecbj{x}_0 \in \mathbb{R}^{N}\\
\siggen:& \quad \,\,\,\, \,\vecbj{u} = \tensorbj{C}_{\sigstatescal} \sigstate, \qquad \dot{\sigstate} = \tensorbj{A}_{\sigstatescal} \sigstate + \tensorbj{G}_{\sigstatescal} \sigstate\supertens{2}, && \, \sigstate(0) = \sigstate_{0} \in \mathbb{R}^q.
\end{align*}
Let $\tensorbj{Q}$ be the constant matrix such that
\begin{align*}
	\tensorbj{Q} \begin{mymatrix}
		\bar{\vecbj{x}} \\
		\bar{\sigstate}
	\end{mymatrix}\supertens{2} = 
	\begin{mymatrix}
		\bar{\vecbj{x}}\supertens{2} \\
		\bar{\vecbj{x}} \otimes \bar{\sigstate} \\
		\bar{\sigstate}\supertens{2}
	\end{mymatrix}  \nonumber \qquad  \text{for arbitrary  } \bar{\vecbj{x}} \in \mathbb{R}^N,\,  \bar{\sigstate}\in \mathbb{R}^q.
\end{align*}
Then we call the autonomous system
\begin{align*}
\parbox{0.5cm}{{\vspace{0.4cm}$\statemapAut:$}} 
& \quad	\,	\tensorcj{E} \dot{\veccj{w}} = \tensorcj{A} \veccj{w} + \tensorcj{G} \veccj{w}\supertens{2}, && \veccj{w}(0) = \veccj{b} \phantom{\in \mathbb{R^{kkk}}} \\
& \qquad  \vecbj{x} = \tensorcj{P}_x \, \veccj{w} &&
\end{align*}
with
\begin{align*}
	\tensorcj{E} &= \begin{mymatrix}
		\tensorbj{E} & \\
					& \tensorbj{I}_q
	\end{mymatrix}
	, \quad 
	\tensorcj{A} =
		\begin{mymatrix}
		\tensorbj{A} & \tensorbj{B} \tensorbj{C}_{\sigstatescal} + \tensorbj{B}_p\tensorbj{C}_{\sigstatescal} \tensorbj{A}_{\sigstatescal} \\
					& \tensorbj{A}_{\sigstatescal}
	\end{mymatrix},
	\quad \tensorcj{P}_x =[\tensorbj{I}_N, \, \tensorbj{0}],
	\quad \veccj{b}= 
	\begin{mymatrix}
		\vecbj{x}_0 \\
		 \sigstate_0
	\end{mymatrix},
\\
	\tensorcj{G} &=
		 \begin{mymatrix}
		\tensorbj{G} & \tensorbj{D}(\tensorbj{I}_N \otimes\tensorbj{C}_{\sigstatescal}) & \tensorbj{G}_u (\tensorbj{C}_{\sigstatescal} \otimes \tensorbj{C}_{\sigstatescal}) + \tensorbj{B}_p \tensorbj{C}_{\sigstatescal} \tensorbj{G}_{\sigstatescal} \\
					& 													  & \tensorbj{G}_{\sigstatescal}
	\end{mymatrix}	\tensorbj{Q},
\end{align*}
the signal generator driven system $\statemapAut$.
\end{definition}
Note that the solution $\vecbj{x}$ of system $\statemap$ for input $\vecbj{u}$ described by the signal generator $\siggen$ and the output $\vecbj{x}$ of the signal generator-driven system $\statemapAut$ from the definition coincide.

\paragraph{Input map with higher-order time derivatives}
When higher-order time derivatives occur in the input map, the further procedure depends on the signal generator. If the signal generator is linear, we can use that for
\begin{align*}
	\vecbj{u} &= \tensorbj{C}_{\sigstatescal} \sigstate, \qquad \dot{\sigstate} = \tensorbj{A}_{\sigstatescal} \sigstate \qquad \sigstate(0) = \sigstate_{0} \quad \text{it holds }
\frac{d^i}{dt^i} \vecbj{u} = \tensorbj{C}_{\sigstatescal} \tensorbj{A}_{\sigstatescal}^i \sigstate.
\end{align*}
Thus, a signal generator driven system, which is quadratic in the extended state $[\vecbj{x};\sigstate]$, can be directly constructed. Only the system matrices $\tensorcj{A}$, $\tensorcj{G}$ have to be slightly adjusted.

If the signal generator is nonlinear, we suggest to further extend the signal generator driven system. We exemplarily discuss this for the case of second order derivatives $\ddot{\vecbj{u}}$: Introduce $\sigstate_1= \dot{\sigstate}$ as a dependent variable and extend the signal generator driven state to $\veccj{w} = [\vecbj{x};\sigstate;\sigstate_1]$. Add  the additional equation
\begin{align*}
	\dot{\sigstate}_1 = \tensorbj{A}_{\sigstatescal} \sigstate_1 + \tensorbj{G}_{\sigstatescal} (\sigstate_1 \otimes \sigstate + \sigstate \otimes \sigstate_1), \qquad \sigstate_1(0) = \sigstate_{10}
\end{align*}
with $\sigstate_{10}$ chosen consistently to the signal generator driven system. Then proceed as in Definition~\ref{def:sg-ode-casc2} to construct the quadratic signal generator driven system with extended state $\veccj{w} = [\vecbj{x};\sigstate;\sigstate_1]$.

\subsection{Input-weighted concept for input-output type methods} \label{subsec:input-weighted}
At least formally, our input-tailoring shows some similarities to the concept of input-weighting. The latter has been used in system-theoretic model reduction of linear systems to get reduced models with enhanced fidelity in certain frequency ranges. We refer to \cite{art:varga-freq-weighted-bal-rel}, \cite{art:breiten-near-optimal-freq-weighted} and references therein
for details.

Motivated by our approach, we propose the usage of input-weights to incorporate non-standard input maps in the system-theoretic methods like multi-moment matching or balanced truncation \cite{art:quadratic-bililinear-bt2017} based on multivariate frequency representations. 
To the best of the authors' knowledge, this has not been discussed before. 
To stress the formal similarities to our input-tailored approach, we use a similar notation.
\begin{definition}[Input-weighted system] \label{def:iw-sys}
Let a system $\statemap$ and an input-weight $\mathbf{F}$ be given as
\begin{align*}
\statemap:& \quad 	\tensorbj{E} \dot{\vecbj{x}} = \tensorbj{A} \vecbj{x} + \tensorbj{G}\vecbj{x}\supertens{2} + \tensorbj{D}  (\vecbj{x} \otimes \vecbj{u})  +  \tensorbj{B} \vecbj{u} + \tensorbj{G}_u \vecbj{u}\supertens{2} +  \tensorbj{B}_p \dot{\vecbj{u}},  &&
	\vecbj{x}(0) = \vecbj{0} \in \mathbb{R}^{N}\\
\mathbf{F}:& \quad \,\,\, \,\vecbj{u} = \tensorbj{C}_{\sigstatescal} \sigstate, \qquad \dot{\sigstate} = \tensorbj{A}_{\sigstatescal} \sigstate + \tensorbj{G}_{\sigstatescal} \sigstate\supertens{2} + \tensorbj{B}_{\sigstatescal} \vecbj{u}_F, && \, \sigstate(0) = \vecbj{0} \in \mathbb{R}^q.
\end{align*}
Then we call $\statemap_F: \vecbj{u}_F \mapsto \vecbj{x}$
\begin{align*}
\parbox{0.7cm}{{\vspace{0.4cm}$\statemap_F:$}} 
& \quad	\,	\tensorcj{E} \dot{\veccj{w}} = \tensorcj{A} \veccj{w} + \tensorcj{G} \veccj{w}\supertens{2} + \tensorcj{B} \vecbj{u}_F, && \veccj{w}(0) = \vecbj{0} \phantom{\in \mathbb{R^{kkk}}} \\
& \qquad  \vecbj{x} = \tensorcj{P}_x \, \veccj{w} &&
\end{align*}
with
\begin{align*}
	\tensorcj{B} =
	\begin{mymatrix}
		\tensorbj{B}_p \tensorbj{C}_{\sigstatescal} \tensorbj{B}_{\sigstatescal} \\
		 \tensorbj{B}_{\sigstatescal}
	\end{mymatrix} \qquad \text{ and  $\tensorcj{E},\tensorcj{A},\tensorcj{G}$ as in Definition~\ref{def:sg-ode-casc2}},
\end{align*}
the input-weighted system.
\end{definition}
The upper input-weighted system $\statemap_F$ results from the assumption that the inputs of interest $\vecbj{u}$ can be constructed by incorporating the input-weight $\mathbf{F}$ and some auxiliary input $\vecbj{u}_F$ into the input-output description. By construction, $\statemap_F$ has a linear input map. Therefore, any standard system-theoretic model reduction method based on the input-independent multivariate frequency representations can be used on it to construct an extended reduction basis $\mathcal{V}$. The reduction basis $\tensorbj{V}$ for the original system $\statemap$ can be extracted from the extended basis $\tensorcj{V}$ in the same fashion as we do it in our input-tailored approach, cf. Remark~\ref{rem:connection-xr-wr}. Of course, the choice of input-weight $\mathbf{F}$ and its influence on the reduction method is an important issue in this approach, but beyond the scope of this work.

\section{Numerical results} \label{sec:assocmor-num-validation}

In this section we numerically investigate the performance of our new input-tailored approximate moment matching method in comparison to existing Galerkin-type reduction methods, such as the system-theoretic multi-moment matching, the trajectory-based proper orthogonal decomposition, and the method on univariate frequency representations proposed in \cite{inproc:fast-nonlin-mor-assoc-trafo}, \cite{art:mor-associated-transform-2016}. Three benchmark examples are considered, which have been used in literature on nonlinear system-theoretic model reduction methods, e.g., \cite{art:quadratic-bilinear-regular-krylov}, \cite{art:quadratic-bililinear-bt2017}, \cite{phd:Gu-mor-nonlinear}, \cite{techrep:two-sided-mm}, \cite{art:two-sided-hermm}, \cite{art:qlmor-gu-2011}. Apart from a general performance study, certain aspects are further highlighted in the different benchmark tests: The handling of non-standard input maps is demonstrated for the viscous Burgers' equation (Section~\ref{sec:burgers-eq}), the treatment of non-trivial initial conditions is showcased for the Chafee-Infante equation (Section~\ref{sec:chaf-eq}). We illustrate the difference between input-tailoring in our method and the use of training trajectories in proper orthogonal decomposition. On the one hand, different input-scenarios may lead to the same input-tailored expansion, although the solution trajectories differ nonlinearly, as discussed for the Chafee-Infante equation. On the other hand, our method is overall less dependent on the input-scenario, as showcased for the nonlinear RC-ladder (Section~\ref{subsec:RC-eq}). A discussion on the difference and computational advantage of our approach to \cite{inproc:fast-nonlin-mor-assoc-trafo}, \cite{art:mor-associated-transform-2016} concludes this section.



\subsection{Setup for numerical results}
The numerical results have been generated with \texttt{MATLAB} Version 9.1.0.441655 (R2016b) on an Intel Core~i5-7500 CPU with 16.0GB RAM. For an efficient realization of our input-tailored moment matching method (referred to as \textit{AssM}, Algorithm~\ref{alg:mm-complete}) we use the Lyapunov equation solver '\texttt{mess{\_}lyap}' from \texttt{M.E.S.S.} Toolbox \cite{SaaKB16-mmess-1.0.1} with its default settings in Step \textit{\eqref{it:alg:mmh2-simple-zik}} of Algorithm~\ref{alg:mmh2-simple}. For comparison, we consider two one-sided multi-moment matching methods from literature: the first one (referred to as \textit{MultM}) is taken from \cite[Alg.~2]{inproc:morBenB12}, \cite{techrep:two-sided-mm} and the second one (referred to as \textit{MpMo}) from \cite{art:two-sided-hermm}.  \textit{MultM} matches the moments of the transfer function $\Transfer{G}_1(\sigma)$ at given frequencies $(s_1, s_2,\ldots,s_\mu)$ up to order $q_1-1$ and the multi-moments of $\Transfer{G}_2(\sigma_1,\sigma_2)$ at the diagonal frequency pairs $((s_1,s_1),(s_2,s_2), \ldots ,$ $(s_\mu,s_\mu))$ up to order $q_2-1$, where $q_1 \geq q_2$ has to be chosen. \textit{MpMo} only takes the expansion frequencies $(s_1, s_2,\ldots,s_\mu)$ as parameters, while the orders up to which the moments are matched are fixed. As discussed in the reference, only at most two basis vectors per frequency are needed in \textit{MpMo}, whereas \textit{MultM} requires three for $q_1=q_2 =1$. In \textit{AssM} the expansion frequencies for the first and the second order transfer functions could be chosen independently of each other. For convenience, we take here the same frequencies and the same number of moments for all expansion frequencies, i.e., $\tilde{L}_1= \cdots =\tilde{L}_\mu$ and $L_1= \cdots =L_\mu$ abbreviated as $\tilde L$ and $L$, respectively. Apart from that, \textit{AssM} requires an additional parameter, the tolerance \textit{tol}, that determines the truncation error in the approximate moment matching. For given order of moments, its influence is moderate. In the following benchmark tests we choose \textit{tol} such that the resulting reduced models are of same dimension as in the other methods. As a heuristic for the selection of expansion frequencies for all methods, we apply the IRKA algorithm to the linear transfer function $\Transfer{G}_1$, as in \cite{inproc:morBenB12}, \cite{art:two-sided-hermm}, \cite{art:quadratic-bilinear-regular-krylov}. The IRKA-points are in general complex-valued. However, it turns out that for few frequencies real values result in all benchmark settings except for the boundary-controlled Chafee-Infante equation. Hence, in \textit{AssM} and \textit{MultM} we deal with real expansion frequencies. In the case of the Chafee-Infante equation in Section~\ref{subsec:chaf-bc} we particularly take the first real IRKA-points and ignore the imaginary ones for simplicity. Treating complex ones is possible but requires some technical adjustments in our approach.
For the method of proper orthogonal decomposition (\textit{POD}), cf., \cite{art:morKunV01}, \cite{art:amsallam-gal-wave}, we use time snapshots of the solution trajectory as training set, unless otherwise stated. In all benchmark test cases, 300 uniformly distributed snapshots are taken, as we do not experience any improvement in the results when increasing the number. The simulations of full order models (FOM) and reduced order models (ROM) are done using \texttt{MATLAB}'s solver '\texttt{ode15s}', where the tolerances are modified to 'AbsTol = $10^{-8}$' and 'RelTol = $10^{-6}$' and the exact Jacobian matrices are forwarded to the solver.

\subsection{Burgers' equation} \label{sec:burgers-eq}
On the spatial domain $\Omega=(0,1)$ we consider the nonlinear viscous Burgers' equation given by
\begin{align*}
	\partial_t v(\spacevar,t) &= - v(\spacevar,t) \, \partial_\spacevar v(\spacevar,t) + \nu \,\partial_{\spacevar \spacevar} v(\spacevar,t)  && \text{in } (0,1) \times (0,T) \\
	v(0,t) &=  u(t), \qquad \quad
	\partial_\spacevar v(1,t) = 0 && \text{in } (0,T)\\
	v(\spacevar,0) &= 0  && \text{on } [0,1] 
\end{align*}
with viscosity constant $\nu = 0.01$. The input $u$ particularly prescribes a Dirichlet boundary condition on the left boundary ($\spacevar=0$). We choose the output to be the boundary value on the right, $y(t)=v(1,t)$.
The two input-scenario cases we present relate to one linear and one nonlinear signal generator:
\begin{itemize}
\item[\emph{Case 1}] \emph{Linear signal generator}.
\begin{align*} 
	u(t) = 0.5\left( \cos{(1.3 \pi t)} - \cos{(5.4 \pi t)} - \sin{(0.6 \pi t)} + 1.2 \sin{(3.1 \pi t)} \right)
\end{align*}
The input $u$ is a sum of sine- and cosine-functions. Every summand can be described by a dynamic system, e.g., the last summand $\tilde{u}(t)=1.2 \sin{(3.1 \pi t)}$ has the linear signal generator
\begin{align*}
	\tilde{u} = [1\,|\,0] \sigstate, \qquad \dot{\sigstate} = 
	3.1\pi \begin{mymatrix}
		 & 1\\
		-1 &
\end{mymatrix}		
	 \sigstate
	  \qquad \,\sigstate(0) = 1.2\begin{mymatrix}
	0 \\
	1
\end{mymatrix},
\end{align*}
analogously for the others. Superposing these single generators yield the linear signal generator for $u$.

\item[\emph{Case 2}] \emph{Nonlinear signal generator}.
\begin{align*} 
	u(t) =  \frac{1}{0.5- \exp{(2t)}} + 2 \exp{(-t)}
\end{align*}
The respective signal generator is nonlinear and reads
\begin{align*}
	u = [-0.5\,|\,2] \sigstate, \quad \dot{\sigstate} = 
	\begin{mymatrix}
		 -2 & \\
		 & -1
\end{mymatrix}		
	 \sigstate +
	 \begin{mymatrix}
		 -0.5 & 0& 0&0 \\
		 0 & 0& 0&0
	\end{mymatrix}		
	 \sigstate\supertens{2},
	  \quad \,\sigstate(0) =
	  \begin{mymatrix}
	4 \\
	1
\end{mymatrix}.
\end{align*}
\end{itemize}
Case 1 is particularly similar to a test case considered in \cite{morBre13}, \cite{inproc:morBenB12}, \cite{art:two-sided-hermm} for multi-moment matching. 

We use this benchmark example to demonstrate the feasibility of handling non-standard input dependencies. We employ two different discretization variants for the Burgers' equation to construct one FOM with linear input-dependency and one with nonlinear input-dependency. Note that the form of the input-dependency has nothing to do with the form of the signal generator. Both FOMs describe, up to a small discretization error,  the same dynamics and should thus serve as an equally valid basis for model reduction. The nonlinear input-dependency, however, cannot be directly treated by the moment matching methods from literature. We show that our input-tailored and input-weighted extensions from Section~\ref{sec:ext-var-appr-app} are applicable. They turn out to be (almost) independent of the underlying discretization.

\begin{table}[t] 
\renewcommand{\arraystretch}{1.25}
\centering\small{
\begin{tabular}{l|l|c}
{Expansion frequencies}
& \textit{AssM}{,} \textit{MultM} & $s_i \in \{0.03, \,  0.22 \}$\\
\hline
\multirow{2}{*}{Order moments}
& \textit{AssM}   & $\tilde{L}\,=3$, \, $L\,=2$ \\
& \textit{MultM}  & $q_1=3$, \, $q_2=2$ \\
\hline
\multirow{2}{*}{Tolerance} & \multirow{2}{*}{\textit{AssM}}  & $tol=10^{-3}$  \, (\textit{Case 1}) \\
		& & $tol=10^{-4}$ \, (\textit{Case 2}) \\ \hline
Resulting dimension & \textit{AssM}, \textit{MultM} & $n = 16$ 		
\medskip
\end{tabular}}
\caption{Reduction parameters for Burgers' equation (\textit{FOM} with $N =4000$).}
\label{tab:redpar-burger}
\end{table}

\begin{figure}[t]
\begin{tabular}{rll|l}
\begin{minipage}{0.022\textwidth}
{\vspace{0.8cm}
{\footnotesize \vspace{0.2cm} \rotatebox{90}{u(t)}\\ \vspace{0.4cm} \vspace{2.9cm}\\
\vspace{0.2cm} \rotatebox{90}{y(t)} \vspace{0.2cm} \\ \vspace{2.4cm} \\
\rotatebox{90}{Output error} 
}}
\end{minipage}
&
{\hspace{-0.4cm}
\begin{minipage}{0.43\textwidth}
\center
\hspace{0.5cm} {{\underline{Case 1}}} \\
\includegraphics[height = 0.85\textwidth, width = 1.0\textwidth]{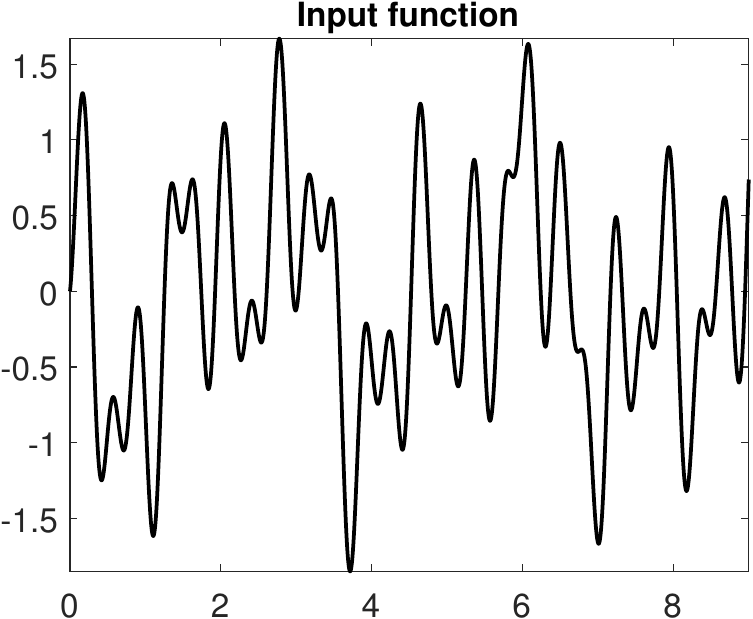} \\
\includegraphics[height = 0.85\textwidth, width = 1.0\textwidth]{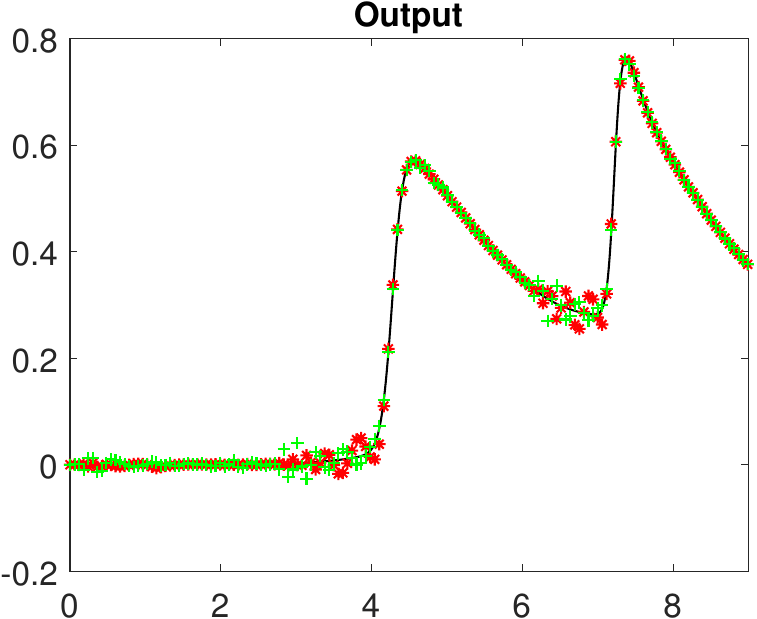} \\
\includegraphics[height = 0.85\textwidth, width = 1.0\textwidth]{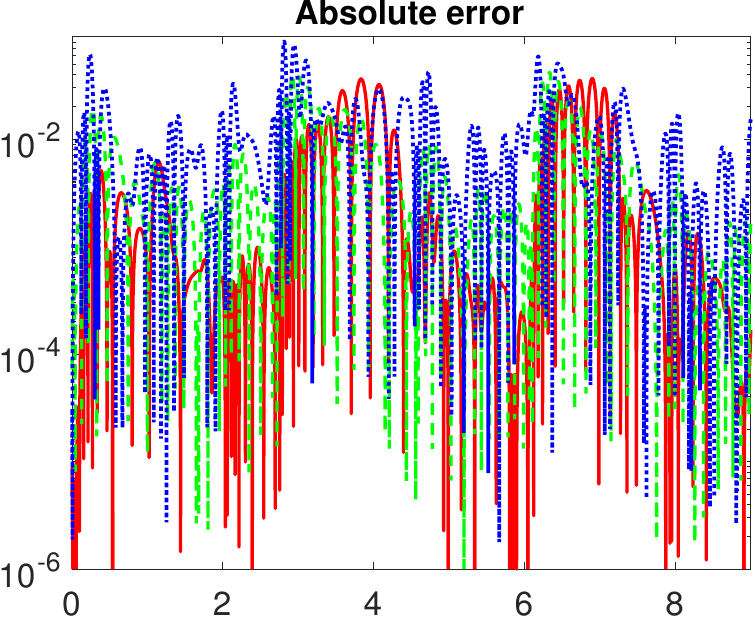} \\
 \hspace{0.5cm} {\footnotesize Time t}
\end{minipage}
}
&
\begin{minipage}{0.00\textwidth}
\end{minipage}
&
\begin{minipage}{0.43\textwidth}
\center
\quad {{\underline{Case 2}}} \\
\includegraphics[height = 0.85\textwidth, width = 1.0\textwidth]{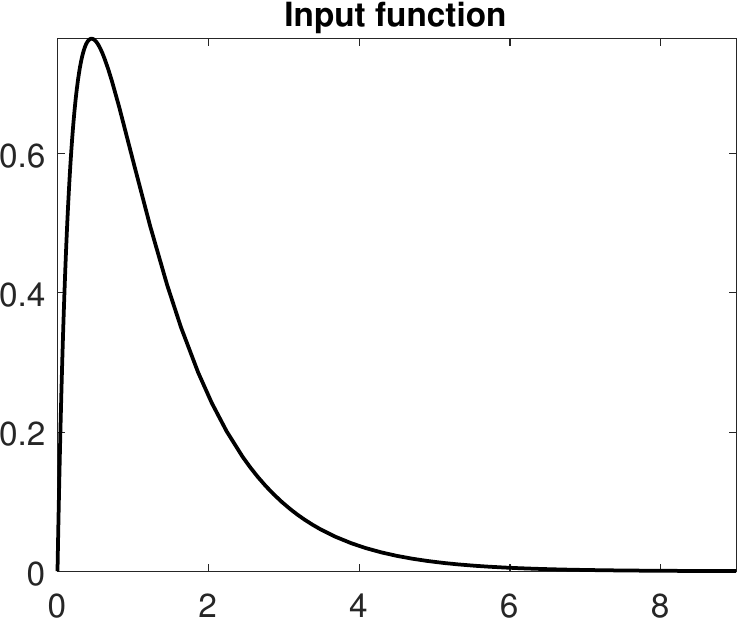} \\
\includegraphics[height = 0.85\textwidth, width = 1.0\textwidth]{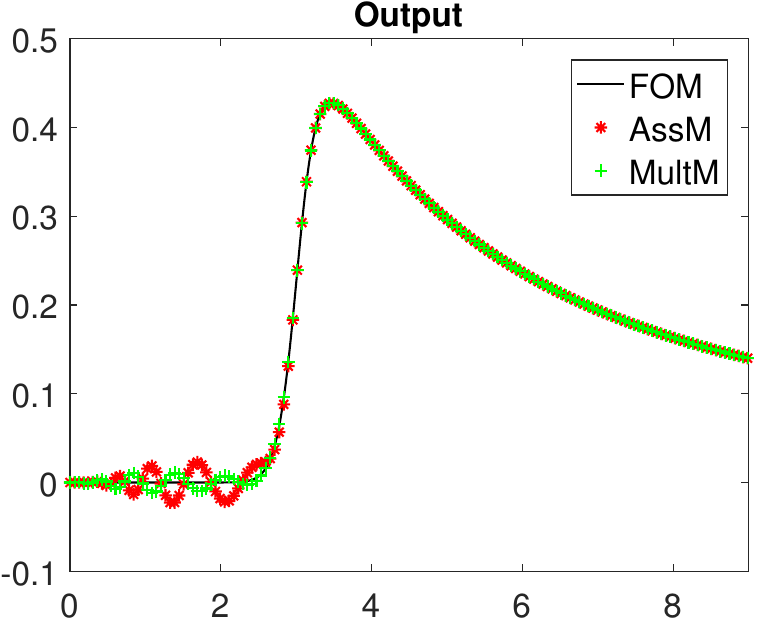} \\
\includegraphics[height = 0.85\textwidth, width = 1.0\textwidth]{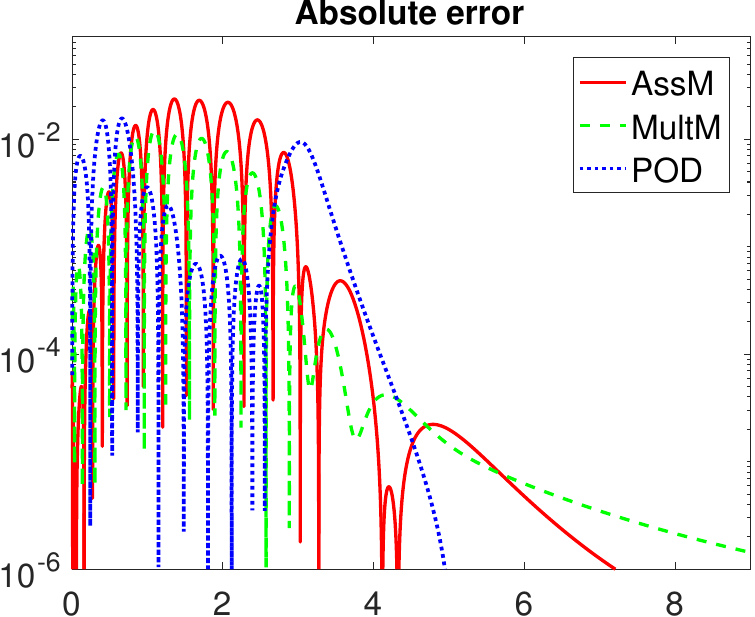}\\
 \hspace{0.5cm} {\footnotesize Time t}
\end{minipage}
\end{tabular}
\caption{Reduction results for Burgers' equation with linear input map. \textit{Top to bottom:} Input~$u$, output $y$, output errors. Dimensions: \textit{FOM}: $N = 4000$, \textit{ROM}: $n=16$ (cf.\ Table~\ref{tab:redpar-burger}).}
\label{fig:inputaware-burgdyn-oscil}
\end{figure}

The Burgers' equation is discretized in space with standard central finite differences and uniform mesh size $h$ implicitly defined by $h=1/(N+2)$ with $N$ inner grid points. Depending on the formulation of the nonlinearity -- advective form $v\,\partial_\spacevar  v$ or conservative form $0.5 \,\partial_\spacevar ( v^2)$ --, we get for the inner node values $v_i(t) \approx v(\spacevar_i,t)$ with $\xi_i = ih$, $1\leq i \leq N$ 
\begin{align*}
	\dot{v}_i &= -v_i  \frac{v_{i+1}-v_{i-1}}{2h} + \nu \frac{v_{i+1}-2 v_{i} + v_{i-1}}{h^2}, \,\, &&\rightarrow \textit{ FOM with linear input map}\\
	\dot{v}_i &= - \frac{v_{i+1}^2-v_{i-1}^2}{4h} \,\,\,\,\,+ \nu \frac{v_{i+1}-2 v_{i} + v_{i-1}}{h^2}, &&\rightarrow \textit{ FOM with nonlinear input map.}
\end{align*}
The discretized boundary conditions give $v_0 = u$ and $(v_{N+1}-v_{N})/h=0$, which we use to eliminate $v_0$ and $v_{N+1}$. This leaves us with the state $\vecbj{x}(t) = [v_1(t);v_2(t);\ldots;v_{N}(t)]$, the tridiagonal (viscosity-associated) system matrix $\mathbf{A}$ and $\tensorbj{E}= \tensorbj{I}_{N}$. The output matrix becomes $\tensorbj{C} = [0,\ldots0,1]$, as $y=v_{N+1}=v_N$ holds due to the boundary conditions. The remaining system matrices depend on the discretization variant. We obtain a quadratic-bilinear system (with linear input map as in \eqref{eq:qldae}) for the advective variant and a quadratic system with a quadratic (nonlinear) input map for the conservative variant.

Proceeding from the FOM with $N=4000$ and linear input-dependency, reduced models of dimension $n=16$ are investigated as example, see Table~\ref{tab:redpar-burger} for the parameters used in model reduction for \textit{AssM} (Algorithm~\ref{alg:mm-complete}) and \textit{MultM}. The respective results concerning output behavior and absolute reduction error over time are illustrated in Fig.~\ref{fig:inputaware-burgdyn-oscil}. As observed, all methods (\textit{AssM},  \textit{MultM} and \textit{POD}) perform comparably well, showing a similar error behavior with moderate numerical oscillations near steep gradients of the solution output in both input-scenario cases (for linear and nonlinear signal generator). This also applies to \textit{MpMo} which we omit here, as it does not provide any additional insight. Notably, \textit{POD} trained with the solution trajectory itself does not lead to significantly better results, which indicates that this benchmark example is rather hard to reduce for any kind of model reduction method.

\begin{figure}[t]
\begin{tabular}{rll|l}
\begin{minipage}{0.022\textwidth}
{\vspace{0.8cm}
{\footnotesize\rotatebox{90}{Output difference}\\ \vspace{1.6cm}\\
\rotatebox{90}{Output difference} 
}}
\end{minipage}
&
{\hspace{-0.4cm}
\begin{minipage}{0.43\textwidth}
\center
\hspace{0.5cm} {{\underline{Case 1}}} \\
\includegraphics[height = 0.85\textwidth, width = 1.0\textwidth]{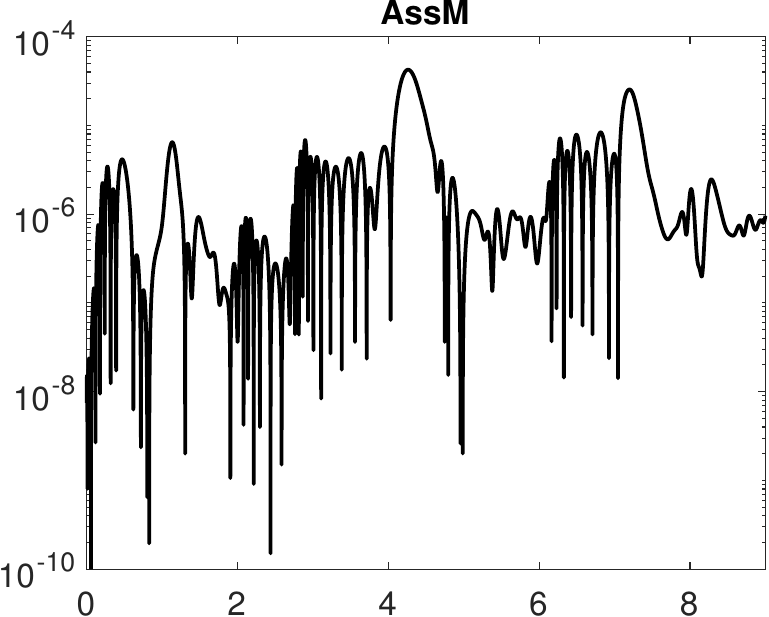} \\
\includegraphics[height = 0.85\textwidth, width = 1.0\textwidth]{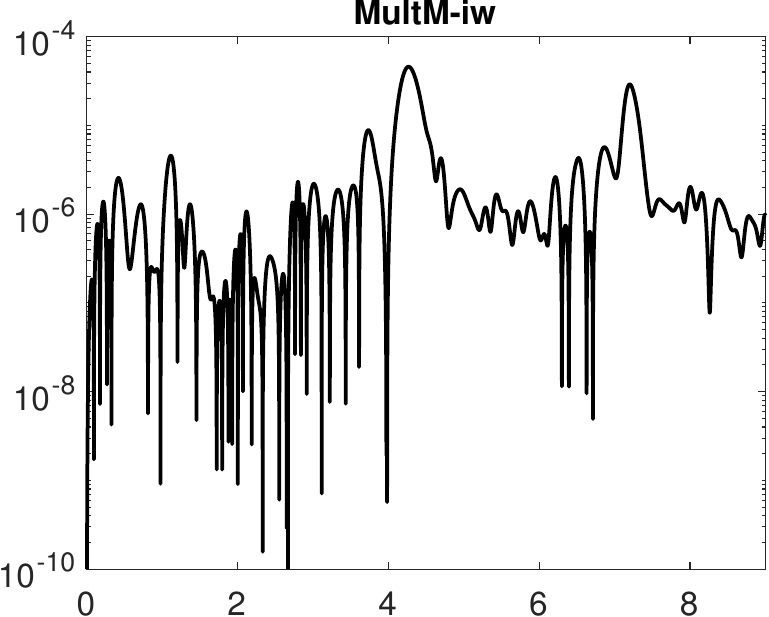}\\
 \hspace{0.5cm} {\footnotesize Time t}
\end{minipage}
}
&
\begin{minipage}{0.00\textwidth}
\end{minipage}
&
\begin{minipage}{0.43\textwidth}
\center
\quad {{\underline{Case 2}}}\\
\includegraphics[height = 0.85\textwidth, width = 1.0\textwidth]{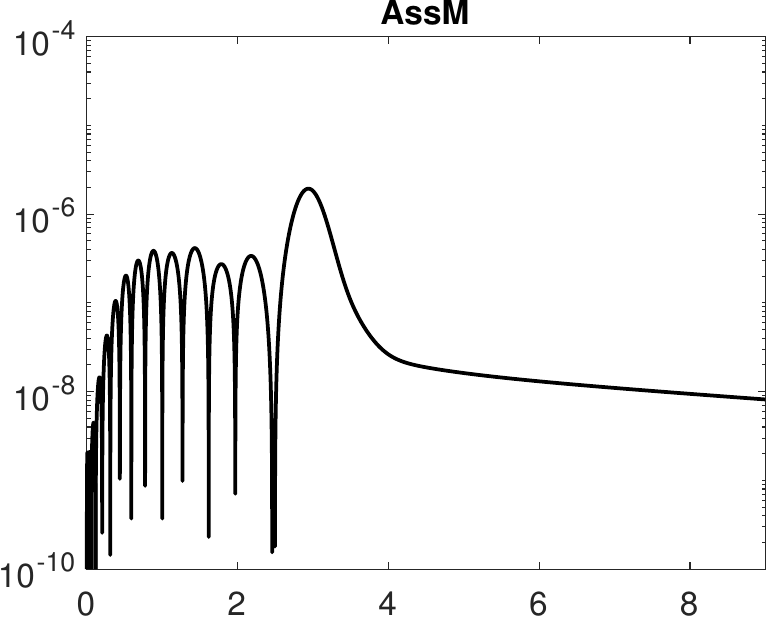} \\
\includegraphics[height = 0.85\textwidth, width = 1.0\textwidth]{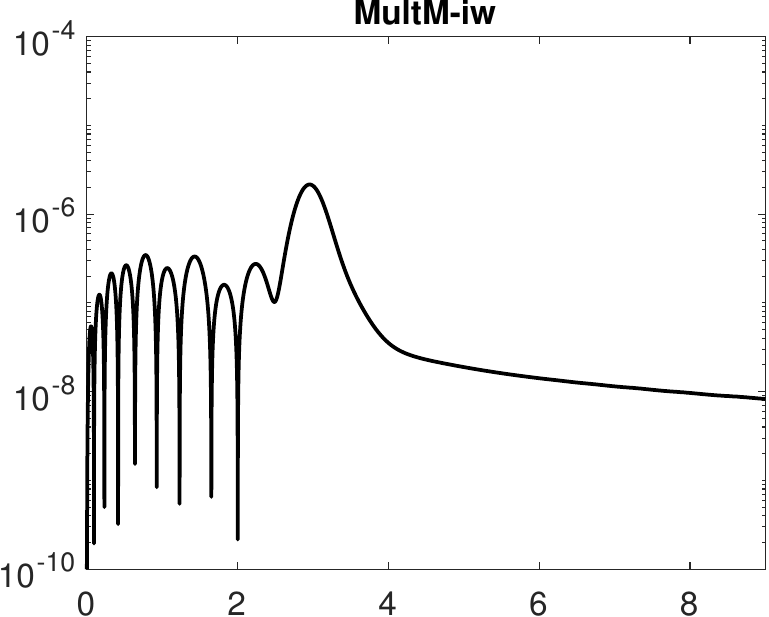}\\
 \hspace{0.5cm} {\footnotesize Time t}
\end{minipage}
\end{tabular}
\caption{Burgers' equation. Output difference of the ROMs associated to the two FOM variants (linear input map vs nonlinear input map)  for \textit{AssM} (top) and \textit{MultM-iw} (bottom). The difference is much smaller than the reduction error, cf., Fig.~\ref{fig:inputaware-burgdyn-oscil} and Fig.~\ref{fig:inputaware-burgdyn-MultM-err}.}
\label{fig:inputaware-burgdyn-AssM-comp}
\end{figure}

\begin{figure}[t]
\begin{tabular}{rll|l}
\begin{minipage}{0.022\textwidth}
{\vspace{0.5cm}
{\footnotesize \rotatebox{90}{Output error}
}}
\end{minipage}
&
{\hspace{-0.4cm}
\begin{minipage}{0.43\textwidth}
\center
\hspace{0.5cm} {{\underline{Case 1}}} \\
\includegraphics[height = 0.85\textwidth, width = 1.0\textwidth]{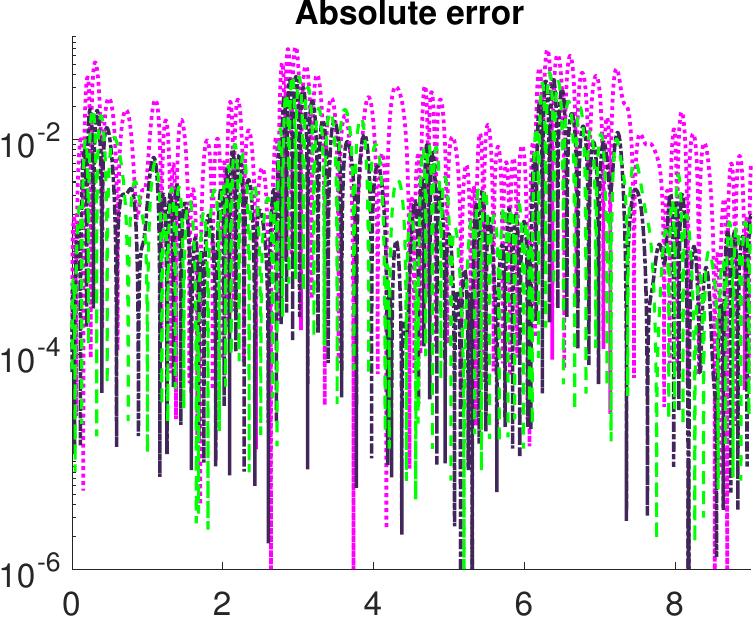}\\
 \hspace{0.5cm} {\footnotesize Time t}
\end{minipage}
}
&
\begin{minipage}{0.00\textwidth}
\end{minipage}
&
\begin{minipage}{0.43\textwidth}
\center
\quad {{\underline{Case 2}}}\\
\includegraphics[height = 0.85\textwidth, width = 1.0\textwidth]{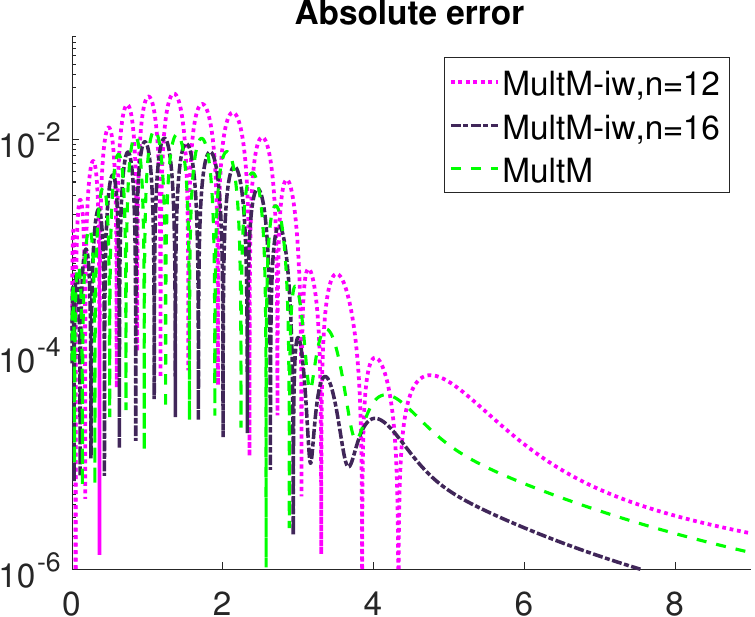}\\
 \hspace{0.5cm} {\footnotesize Time t}
\end{minipage}
\end{tabular}
\caption{Multi-moment matching results for Burgers' equation with linear input map. Parameters in Table~\ref{tab:redpar-burger} yield unweighted \textit{MultM} with $n=16$ and input-weighted {\textit{MultM-iw}} with $n=12$. {\textit{MultM-iw}} with $n=16$ is generated using instead $q_2=3$  for $s_1 = 0.03$.}
\label{fig:inputaware-burgdyn-MultM-err}
\end{figure}

As for the FOM with nonlinear input-dependency and the reduction parameters in Table~\ref{tab:redpar-burger}, \textit{AssM} yields the same results. The dimension of the resulting ROM is identical ($n=16$). The output responses of the ROMs associated to the two underlying discretization (FOM) variants only differ in the order of the discretization error which is much smaller -- by several orders of magnitude -- than the actual reduction error, as seen from Fig.~\ref{fig:inputaware-burgdyn-AssM-comp} in comparison to the output reduction error in Fig.~\ref{fig:inputaware-burgdyn-oscil}. In contrast to \textit{AssM}, \textit{MultM} cannot handle the  nonlinear input-dependency directly. Thus, we consider the input-weighted adaption of multi-moment matching (referred to as \textit{MultM-iw}) that we have proposed in Section~\ref{subsec:input-weighted}. For the input-weight, we exemplarily choose
\begin{align*}
u = \sigstatescal, \qquad \dot{\sigstatescal} = - \sigstatescal + {u}_F, && \, \sigstatescal(0) = 0,
\end{align*}
(i.e., $\tensorbj{C}_{\sigstatescal}=1$, $\tensorbj{G}_{\sigstatescal} = 0$, $\tensorbj{A}_{\sigstatescal} = -1$ and $\tensorbj{B}_{\sigstatescal} = 1$ in Definition~\ref{def:iw-sys}). Applying \textit{MultM-iw} to both FOM variants (with linear and nonlinear input map) we observe the same independence from the underlying discretization as for \textit{AssM}, cf.\ Fig.~\ref{fig:inputaware-burgdyn-AssM-comp}. 
In case of linear input maps, \textit{MultM-iw} typically leads to a ROM of smaller dimension than the standard unweighted method \textit{MultM}, because the input-weighted system, cf.\ Definition~\ref{def:iw-sys}, does not have any bilinear parts to be considered in the multi-moment matching. For the considered test scenarios with the reduction parameters in Table~\ref{tab:redpar-burger}, we get $n=12$ for \textit{MultM-iw} versus $n=16$ for \textit{MultM}. In this specific example of the Burgers' equation, the smaller dimension goes hand in hand with a slightly larger output error as seen in Fig.~\ref{fig:inputaware-burgdyn-MultM-err}. Also the choice of input weight might have an influences on the reduction results. However, we note that we have tested the method with altered reduction parameters (including different choices of input-weights) and observed comparable results when the reduced models are constructed to be of equal dimension, e.g., by incorporating an additional expansion frequency for \textit{MultM-iw}, see Fig.~\ref{fig:inputaware-burgdyn-MultM-err}.

\subsection{Chafee-Infante equation} \label{sec:chaf-eq}

The Chafee-Infante equation is a one-dim\-en\-sio\-nal convec\-tion-diffusion equation for $v = v(\spacevar,t)$ with a cubic nonlinearity in $v$. Following \cite{art:two-sided-hermm}, we introduce the augmented function $w$ by $w= v^2$ and consider an artificial differential equation describing $w$ by differentiating the algebraic relation to get $\partial_t w= 2 v\, \partial_t v$. By that a partial differential equation with only quadratic nonlinearities results. It reads
\begin{align*}
	\partial_t v(\spacevar,t) &= - v(\spacevar,t) \, w(\spacevar,t) + \partial_{\spacevar \spacevar} v(\spacevar,t) + v(\spacevar,t)  &&\text{in } (0,1) \times (0,T) \\
	\partial_t w(\spacevar,t) &= -2  w(\spacevar,t)^2 + 2 v(\spacevar,t) \, \partial_{\spacevar \spacevar} v(\spacevar,t) + 2 v(\spacevar,t)^2 &&\text{in } (0,1) \times (0,T) \\
	\gamma v(0,t) &+ (1-\gamma) \partial_\spacevar v(0,t) = u(t), \hspace{0.8cm}	\partial_\spacevar v(1,t) = 0 &&\text{in } (0,T)\\
	w(0,t)&= v(0,t)^2, \hspace{2.75cm}  \,\, w(1,t) = v(1,t)^2  &&\text{in } (0,T)\\
  v(\spacevar,0) &= v^0(\spacevar),  \hspace{3.1cm} w(\spacevar,0) = v^0(\spacevar)^2   &&\text{on } [0,1].
\end{align*}
The parameter $\gamma$ is varied between $\gamma=1$ and $\gamma=0$ in the test studies, which relates to Dirichlet- and Neumann-boundary conditions described by the input $u$ on the left boundary ($\spacevar=0$), respectively. The equations for $w(0,t)$ and $w(1,t)$ should be read as consistency conditions. Similarly as for the Burgers' equation, we discretize the system in space using central finite differences with a uniform mesh with $\tilde N$ inner grid points and eliminate the boundary node values by means of the boundary conditions. This leads to a quadratic-bilinear system of the form \eqref{eq:qldae} with state
$\vecbj{x}(t) =[v_1(t);\ldots;v_{\tilde N}(t); w_1(t);\ldots;w_{\tilde N}(t)]$ and $\tensorbj{E}= \tensorbj{I}_{N}$, $N=2\tilde N $ ($=1500$ here).

Two types of scenarios are considered for the Chafee-Infante equation, a standard input-output scenario (boundary-controlled scenario) in Section~\ref{subsec:chaf-bc} and an uncontrolled test scenario in Section~\ref{subsec:Chafunc}. The uncontrolled scenario, in which the dynamics is exclusively driven by the non-trivial initial conditions, comes from \cite{techrep:two-sided-mm}, cf., also \cite{art:two-sided-hermm}. 

\subsubsection{Boundary-controlled scenario} \label{subsec:chaf-bc}
Trivial initial conditions for $v^0$ are assumed. With $\gamma=1$ the input $u$ takes the role of a Dirichlet boundary condition. It is varied over the test cases by a linear scaling $\alpha \in \mathbb{R}$ according to
\begin{align*}
	u(t) &= \alpha \left[ \cos{(1.3 \pi t)} - \cos{(5.4 \pi t)} - \sin{(0.6 \pi t)} + 1.2 \sin{(3.1 \pi t)} \right] \hspace*{1.65cm}\\
 \emph{Case 1}: \qquad \alpha &= 1.\\
 \emph{Case 2}: \qquad \alpha&=0.125. 
\end{align*}
By construction, the corresponding signal generators of both cases coincide up to a scaling. (They are scalings of the signal generator in Case~1 of Section~\ref{sec:burgers-eq}.) As output we consider $y(t)= v(1,t)$, implying the output matrix $\tensorbj{C}= [\vecbj{0}_{1,\tilde N-1},1 , \vecbj{0}_{1,\tilde N}]$, since $v_{\tilde N}=v_{\tilde N +1}$ holds due to the boundary conditions.

\begin{table}[t]
\renewcommand{\arraystretch}{1.25}
\centering\small{
\begin{tabular}{l|l|c}
Expansion frequencies & \textit{AssM}{,} \textit{MultM} & $s_i \in \{1.5, \,   21.5, \,   48.3\}$\\ \hline
\multirow{3}{*}{Order moments}
& \textit{AssM} & $\tilde{L}=1$, \, $L=2$  \\
& \multirow{2}{*}{\textit{MultM}}  & \,$q_1=2$, \, $q_2= 2$ \quad ($s_i \in\{1.5, \, 21.5\}$)\\
&				  & \,$q_1=2$, \, $q_2= 1$ \quad ($s_3 = 48.3$) \, \phantom{1234}\\
\hline
\multirow{2}{*}{Tolerance} & \multirow{2}{*}{\textit{AssM} } & $tol=10^{-3}$ \, (\textit{Case 1}) \\ 
		  & 					 & $tol=10^{-4}$ \, (\textit{Case 2}) \\
\hline
Resulting dimension & \textit{AssM}{,} \,\textit{MultM} & $n=12$
\medskip
\end{tabular}}
\caption{Reduction parameters for controlled Chafee-Infante equation (\textit{FOM} with $N =1500$).}
 \label{tab:redpar-chafee}
\end{table}

\begin{figure}[t]
\begin{tabular}{rll|l}
\begin{minipage}{0.022\textwidth}
{\vspace{0.8cm}
{\footnotesize \vspace{0.4cm} \rotatebox{90}{y(t)}\\ \vspace{0.2cm} \vspace{2.9cm}\\
\rotatebox{90}{Output error} \\ \vspace{2.9cm} \\
\rotatebox{90}{Output error} 
}}
\end{minipage}
&
{\hspace{-0.4cm}
\begin{minipage}{0.43\textwidth}
\center
\hspace{0.5cm} {{\underline{Case 1}}} \\
\includegraphics[height = 0.85\textwidth, width = 1.0\textwidth]{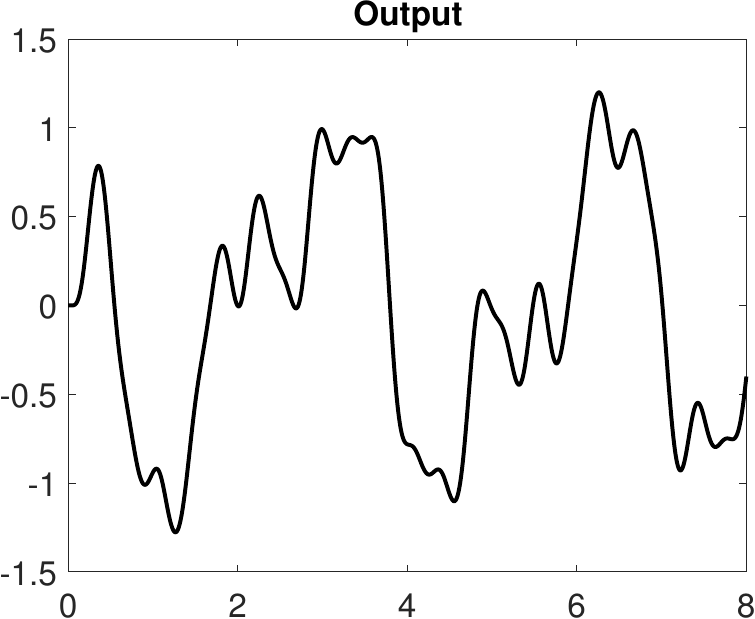} \\
\includegraphics[height = 0.85\textwidth, width = 1.0\textwidth]{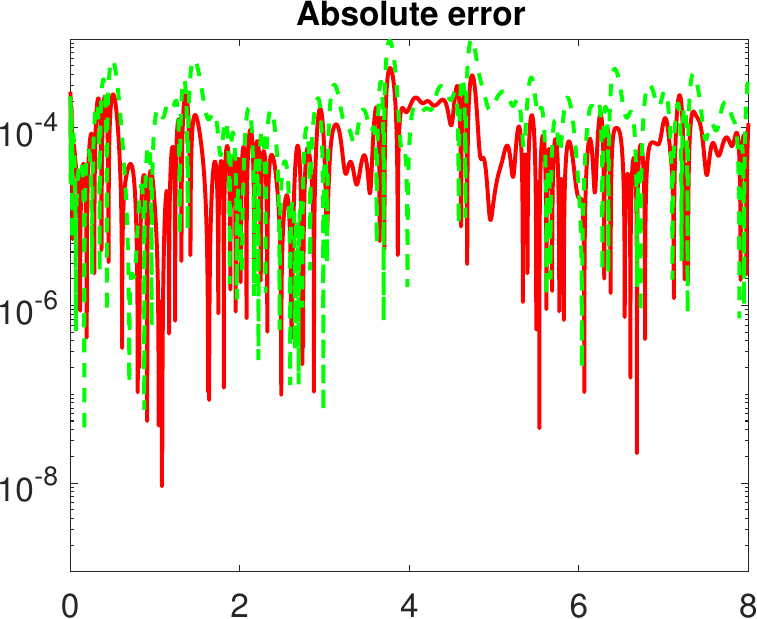} \\
\includegraphics[height = 0.85\textwidth, width = 1.0\textwidth]{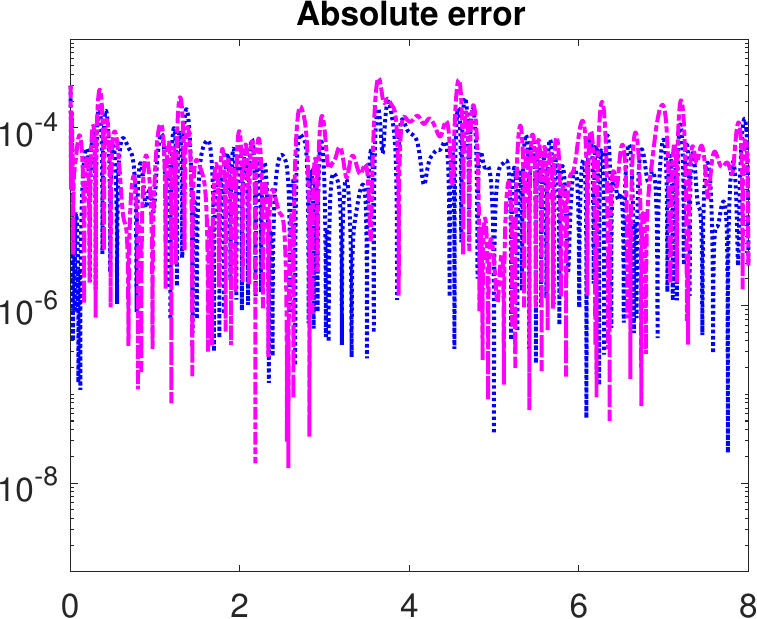} 
\end{minipage}
}
&
\begin{minipage}{0.00\textwidth}
\end{minipage}
&
\begin{minipage}{0.43\textwidth}
\center
\quad {{\underline{Case 2}}} \\
\includegraphics[height = 0.85\textwidth, width = 1.0\textwidth]{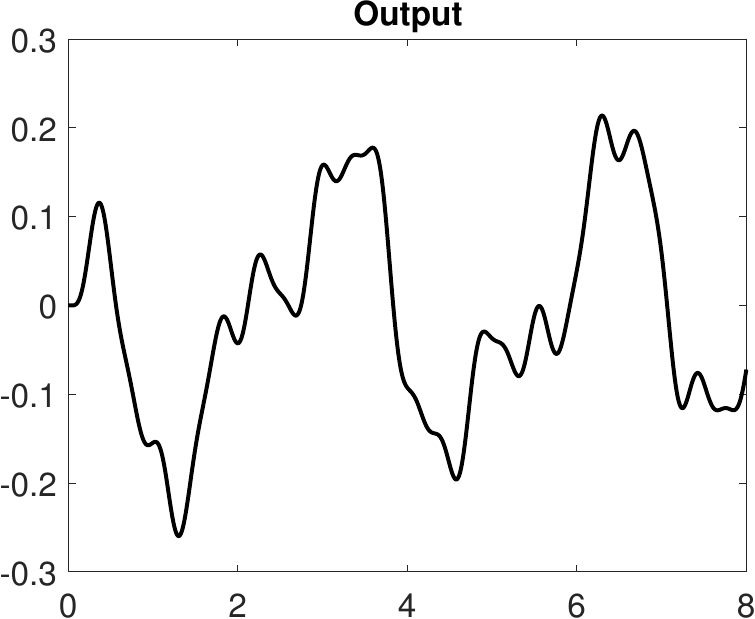} \\
\includegraphics[height = 0.85\textwidth, width = 1.0\textwidth]{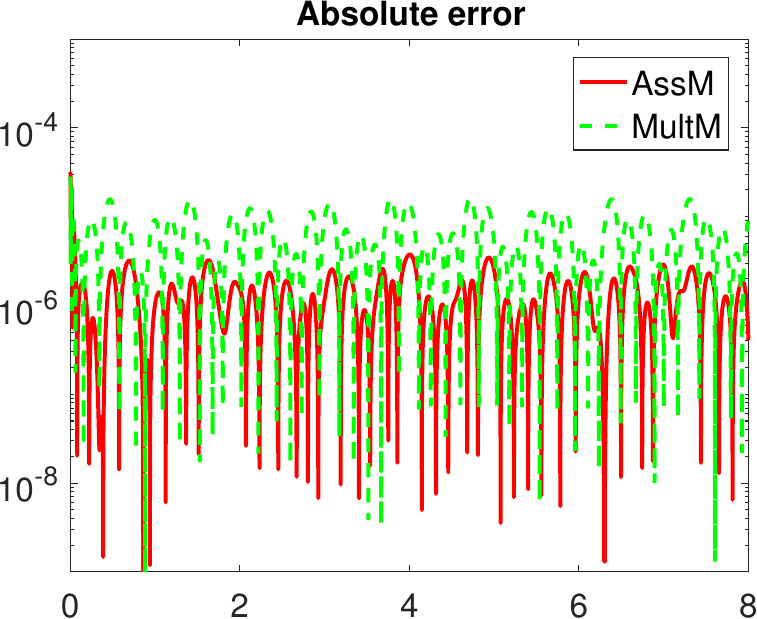} \\
\includegraphics[height = 0.85\textwidth, width = 1.0\textwidth]{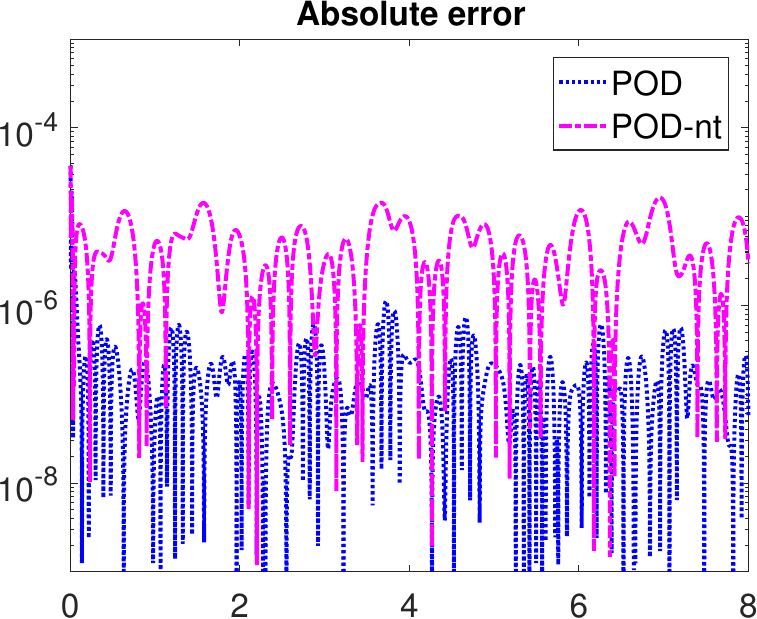} 
\end{minipage}
\end{tabular}
\caption{Reduction results for controlled Chafee-Infante equation. \textit{Top to bottom:} Output $y$ (with \textit{FOM}), output errors for case-independent \textit{AssM} and \textit{MultM}, output errors for case-dependent \textit{POD}. (\textit{POD-nt} is trained with the respective other case.) Dimensions: \textit{FOM}: $N=1500$, \textit{ROM}: $n=12$  (cf.\ Table~\ref{tab:redpar-chafee}).}
\label{fig:inputaware-chafdyn-oscil}
\end{figure}


Reduced models of dimension $n=12$ are constructed for  \textit{AssM} and \textit{MultM} with the parameters of Table~\ref{tab:redpar-chafee}. As in \cite{art:two-sided-hermm} for a similar case, \textit{MpMo} leads to very poor results, here even unstable simulations, and is thus not 
regarded. Also the direct application of proper orthogonal decomposition onto the full state results in significantly poor results, which is known to possibly happen when different physical variables are mixed
 \cite{art:amsallam-gal-wave}, \cite{art:morKunV01}. Hence, we use instead a block-structured version of proper orthogonal decomposition (also referred to as \textit{POD} here), in which two bases of dimension $6$ for the physical variables $[v_1(t);\ldots;v_{\tilde N}(t)]$  and $[w_1(t);\ldots;w_{\tilde N}(t)]$ are constructed separately and combined to a full block basis of dimension $12$.  As \textit{POD} depends on the chosen training trajectory, the two cases lead to two distinct reduced models. Results for which the training and solution (test) trajectory differ are indicated by the suffix '\textit{-nt}'. The training trajectory here is then from the other case. In contrast, \textit{AssM} leads to the same reduced model in both cases, as the input-tailored variational expansions coincide -- up to the scaling $\alpha$. Note that this scaling in the input is compensated by the scaling of \textit{tol} in the approximate moment matching, cf.\ Table~\ref{tab:redpar-chafee}.
 
As seen in Fig.~\ref{fig:inputaware-chafdyn-oscil}, the case-independent \textit{AssM} performs well for both cases, especially also better than \textit{MultM}. The results of proper orthogonal decomposition depend on the training trajectory. Worse results are observed for \textit{POD-nt} than for the perfectly trained model \textit{POD}. In particular, \textit{POD-nt} performs worse than our \textit{AssM}. 

\begin{remark}\label{rem-chaf-quadratize}
It should be mentioned that model order reduction can certainly also be applied directly to the original nonlinear problem formulation (without quadratric-bilinear reformulation). When this is done, the model reduction error is lower for any projec\-tion-based method. The drawback is that the resulting reduced models are not online-efficient and thus need to be augmented by complexity reduction in practice, cf., \cite{art:goyal-polynomial2019}, \cite{art:kramer-liftingPOD}. For completeness, corresponding results are presented in Appendix~\ref{app:nolin-mm-mor}.
\end{remark}

\subsubsection{Uncontrolled scenario} \label{subsec:Chafunc}

\begin{table}[t] \label{tab:chafunc}
\renewcommand{\arraystretch}{1.25}
\centering\small{
\begin{tabular}{l|l|c|c}
 &  & \textit{Smaller ROM} & \textit{Larger ROM} \\
\hline
\multirow{2}{*}{Expansion frequencies}
& \textit{AssM}, \textit{MultM} & {$s_1=4.77$} & {$s_1=4.77$} \\
& \textit{MpMo}  & $5$ IRKA-points & $9$ IRKA-points  \\
\hline
\multirow{2}{*}{Order moments}
 & \textit{AssM} & $\tilde{L}=2$, \, \,$L=2$ & $\tilde{L}=4$, \, \,$L=3$ \\
 & \textit{MultM} &  $q_1=3$, \, $q_2=3$  &  $q_1=6$, \, $q_2=4$ \\
\hline
Tolerance & \textit{AssM}  & $tol=5\cdot 10^{-5}$ & $tol =10^{-7}$\\
\hline
\multirow{2}{*}{Resulting dimension} & \textit{AssM}, \textit{MultM}& $n=10$ & $n=19$ \\ 
		& \textit{MpMo}  & $n=\phantom{1}9$ & $n=19$ \\

\end{tabular}}
\caption{Reduction parameters for uncontrolled Chafee-Infante equation to generate smaller and larger ROM. IRKA-points are obtained with the artificial single-output matrix $\vecbj C = [{\vecbj{1}}_{1,\tilde{N}},{\vecbj{0}}_{1,\tilde{N}}]/\tilde{N}$. (\textit{FOM} with $N =1500$).}
 \label{tab:redpar-Chafunc}
\end{table}


\begin{figure}[t]
\center
\begin{tabular}{rl}
\begin{minipage}{0.02\textwidth}
{\vspace{0.8cm}
{\footnotesize \vspace{0.4cm} \rotatebox{90}{\qquad State $v$}} 
}
\end{minipage}
&
\begin{minipage}{0.54\textwidth}
\hspace{-0.3cm}
\includegraphics[height = 0.75\textwidth, width = 1.0\textwidth]{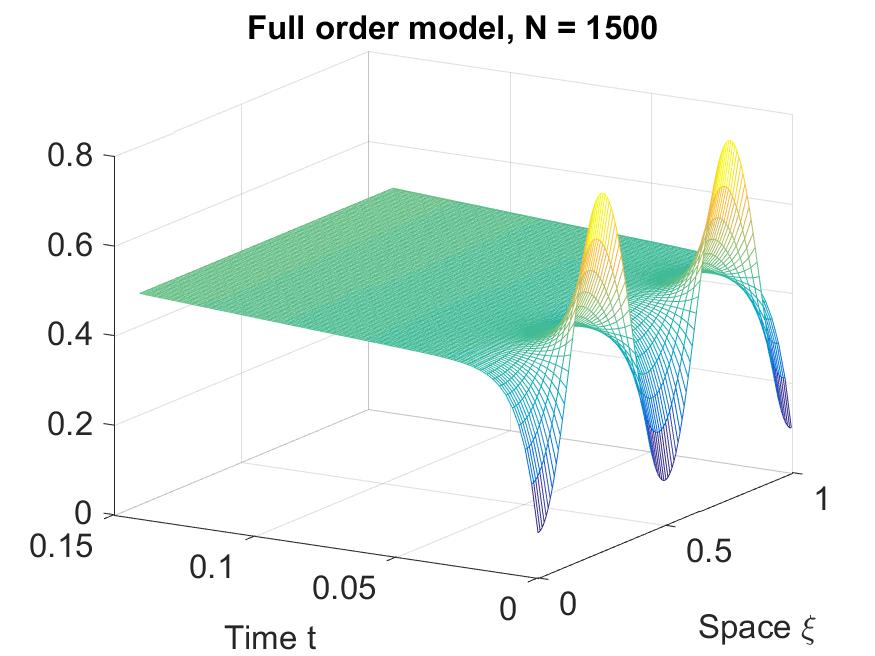}
\end{minipage}
\end{tabular}
\caption{Uncontrolled Chafee-Infante equation. \textit{FOM} solution over space and time, $N=1500$ and $t \in [0,0.15]$.}
\label{fig:chafunc_sol}
\end{figure}

\begin{figure}[t]
\begin{tabular}{rll|l}
\begin{minipage}{0.02\textwidth}
{\vspace{0.8cm}
{\footnotesize \vspace{0.6cm} \rotatebox{90}{Error in $v$}\\ \vspace{2.4cm}\\
\rotatebox{90}{Error in $v$} \\ \vspace{2.4cm} \\
\rotatebox{90}{\quad Error in $v$} \vspace{0.2cm}
}}
\end{minipage}
&
{\hspace{-0.3cm}
\begin{minipage}{0.425\textwidth}
\center
\hspace{0.5cm} {{\underline{Smaller reduced models}}} \\
\includegraphics[height = 0.75\textwidth, width = 1.0\textwidth]{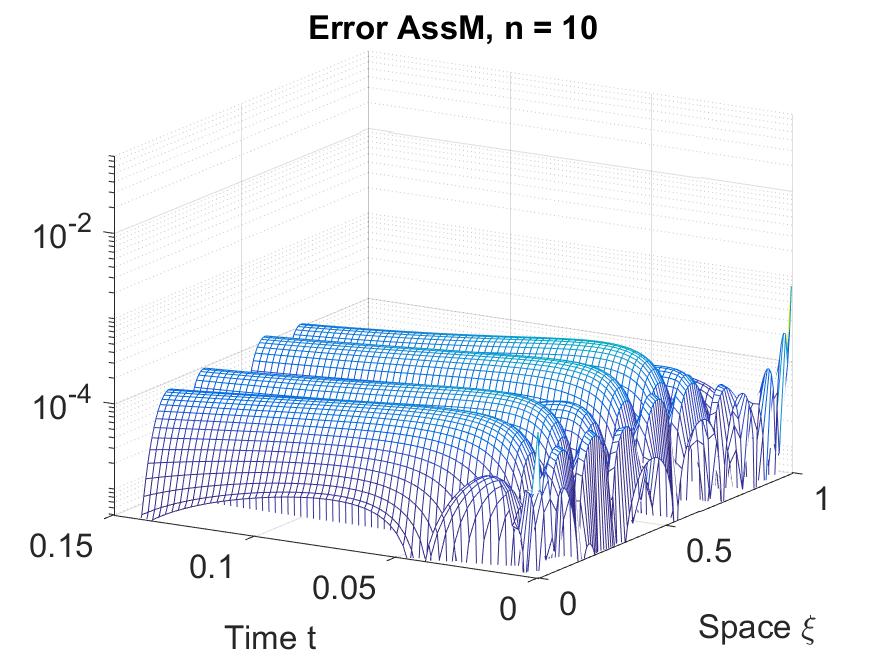} \\
\includegraphics[height = 0.75\textwidth, width = 1.0\textwidth]{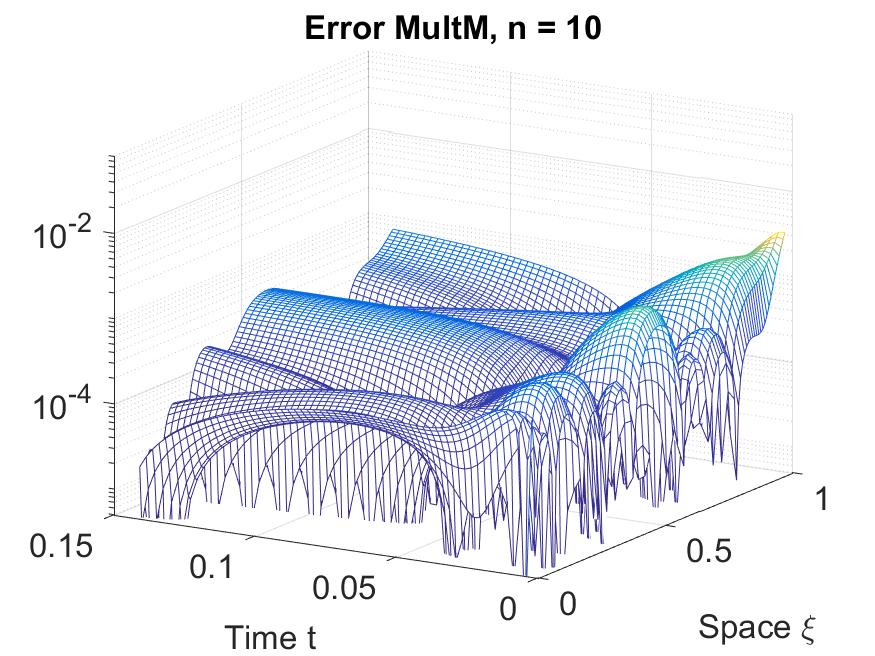} \\
\includegraphics[height = 0.75\textwidth, width = 1.0\textwidth]{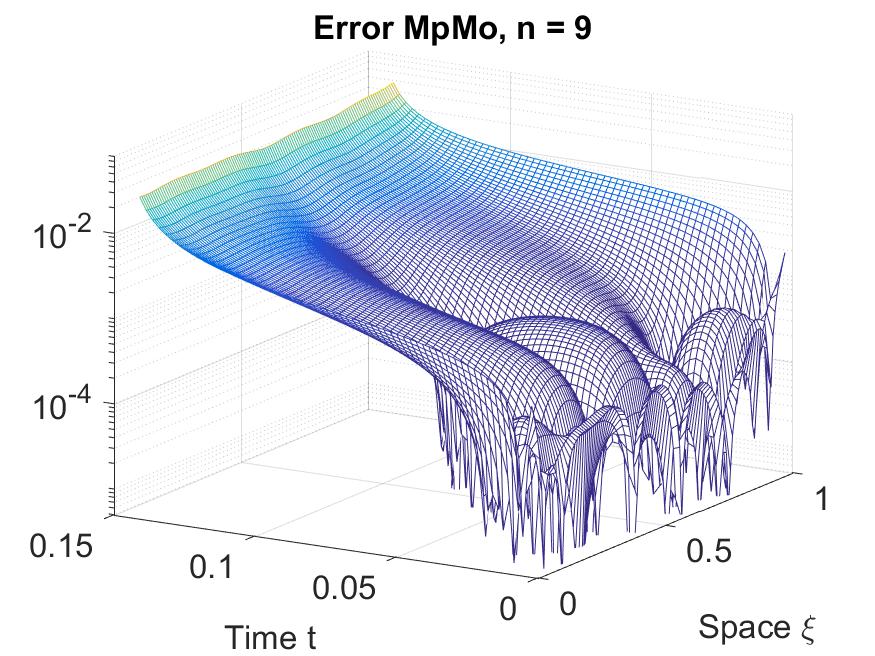} 
\end{minipage}
}
&
\begin{minipage}{0.00\textwidth}
\end{minipage}
&
\begin{minipage}{0.425\textwidth}
\center
\quad {{\underline{Larger reduced models}}} \\
\includegraphics[height = 0.75\textwidth, width = 1.0\textwidth]{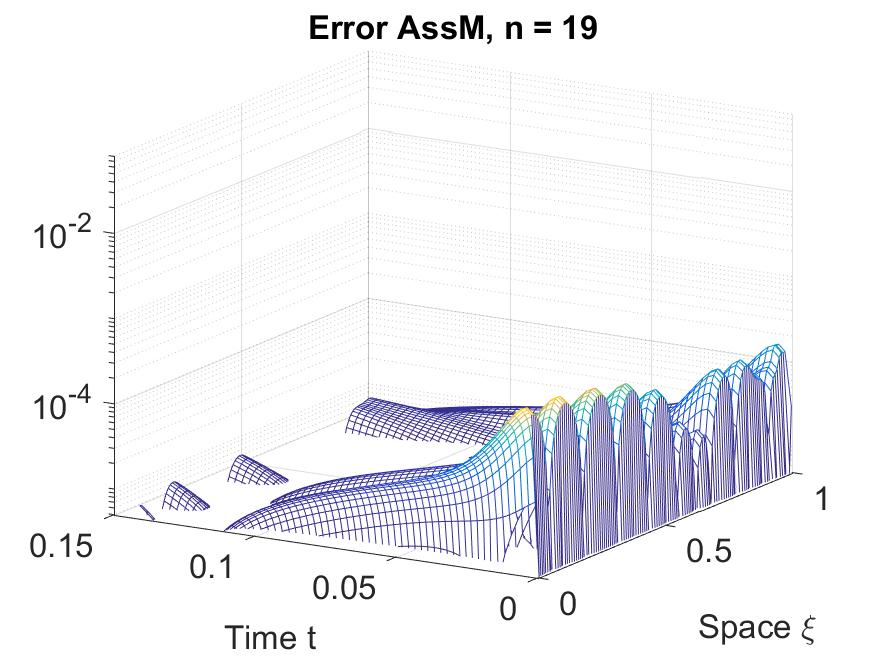} \\
\includegraphics[height = 0.75\textwidth, width = 1.0\textwidth]{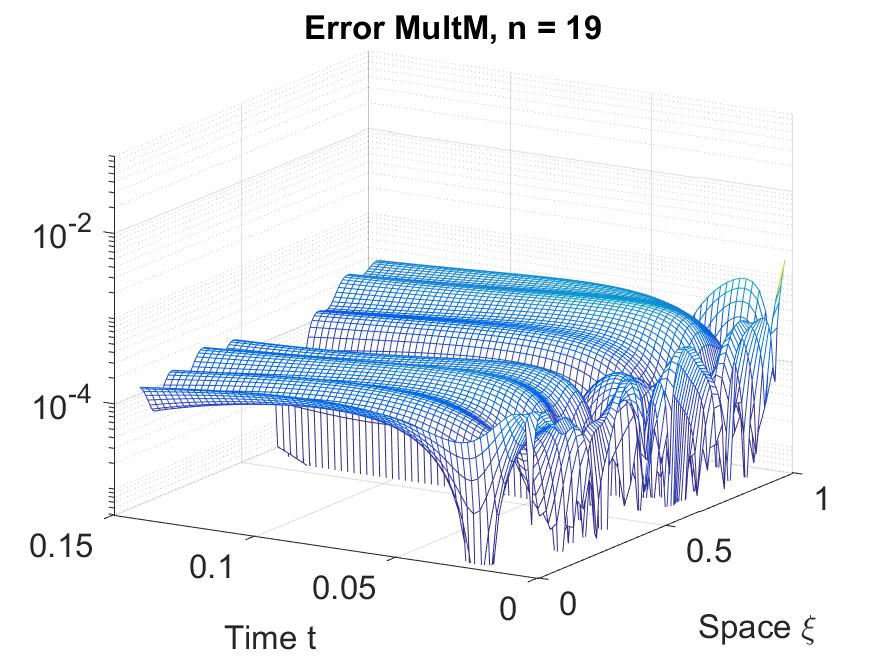} \\
\includegraphics[height = 0.75\textwidth, width = 1.0\textwidth]{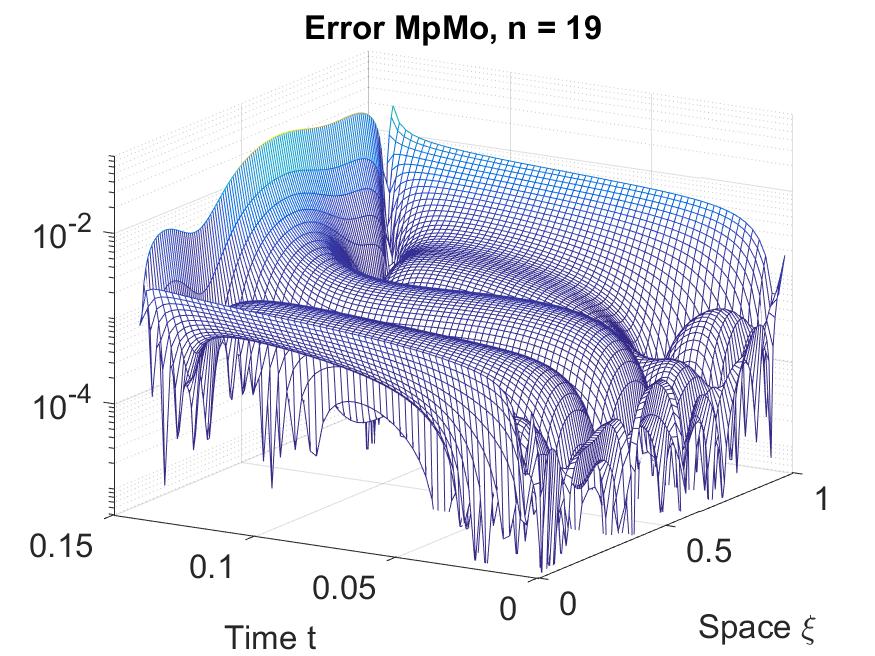} 
\end{minipage}
\end{tabular}
\caption{Reduction results for uncontrolled Chafee-Infante equation. Absolute errors for the full reconstructed state $v$ over space and time. \textit{Top to bottom:} \textit{AssM}, \textit{MultM} and \textit{MpMo}. Dimensions: \textit{FOM}: $N=1500$, \textit{ROM}: $n=10$ ($n=9$ for \textit{MpMo}) on the left, and $n=19$ on the right (cf.\ Table~\ref{tab:redpar-Chafunc}).}
\label{fig:chafunc_err}
\end{figure}

Trivial Neumann boundary conditions are assumed at $\spacevar = 0$, i.e., $\gamma = 0$ and $u(t) = 0$ for $t \geq 0$. The dynamics is exclusively driven by a non-zero initial condition given as
\begin{align*}
	v^0(x) = \frac{1}{10} + \frac{7}{10}  (\sin{((2 x + 1)\pi)})^2, \qquad x \in [0,1].
\end{align*}
The respective \textit{FOM} solution is illustrated in a space-time plot in Fig.~\ref{fig:chafunc_sol}. 

Whereas our \textit{AssM} can handle the set-up directly and without any problems, other system-theoretic methods cannot cope with the non-trivial initial conditions and require an auxiliary reformulation of the full order model. Following the ideas of \cite{art:two-sided-hermm}, we apply \textit{MultM} and \textit{MpMo} to the differential system for the shifted state $\tilde{\vecbj{x}}(t) = \vecbj{x}(t) - \vecbj{x}(0)$. By construction it has zero initial conditions and a constant input-term $\tensorbj B u$ given by $\tensorbj B = \tensorbj A \vecbj x(0)+ \tensorbj G  (\vecbj{x}(0))\supertens{2}$ and $u(t) = 1$ for $t \geq 0$.
The shifted system is also used for the heuristic determination of expansion frequencies by IRKA, where we take the single-output matrix $\vecbj C = [{\vecbj{1}}_{1,\tilde{N}},{\vecbj{0}}_{1,\tilde{N}}]/\tilde{N}$, which relates to the average over the original state, as in \cite{art:two-sided-hermm}, \cite{techrep:two-sided-mm}.

For the moment matching methods reduced models of two different sizes are generated with the parameters from Table~\ref{tab:chafunc}. The resulting absolute errors in the full state $v$ are illustrated in a space-time plot for $t \in [0,0.15]$, Fig.~\ref{fig:chafunc_err}. Regardless of the model size \textit{AssM} obviously performs best, followed by \textit{MultM} with errors that are up to an order of magnitude higher. Worst results are observed for \textit{MpMo}. Note that the \textit{MpMo} models even become unstable for $t \approx 0.2$, which is in accordance with the results reported in \cite{art:two-sided-hermm} for the one-sided method.

\newpage
\subsection{Nonlinear RC-ladder} \label{subsec:RC-eq}
This benchmark describes a nonlinear RC-ladder with $\tilde{N}$ capacitors and I-V diodes. The nonlinearity is due to the diode I-V characteristics, given by $g(v) = \exp{(40v)} - 1$ for voltages $v$. We use the same setup as in \cite{art:quadratic-bilinear-regular-krylov}, \cite{art:quadratic-bililinear-bt2017}, \cite{art:two-sided-hermm}, but also in \cite{inproc:fast-nonlin-mor-assoc-trafo}, \cite{art:mor-associated-transform-2016} a similar example has been studied. The node voltages $v_i$ ($2\leq i \leq \tilde{N}-1,$ and $\tilde{N}=500$) are described by 
\begin{align*}
	\dot{v}_1(t) &= -2 v_1(t) + v_2(t)- g(v_1(t)) - g(v_1(t)- v_2(t)) + u(t) \\
	\dot{v}_i(t) &= -2 v_i(t) + v_{i-1}(t) + v_{i+1} (t) + g(v_{i-1}(t) - v_i(t)) - g(v_i(t) - v_{i+1}(t))\\
	\dot{v}_{\tilde{N}}(t) &= -v_{\tilde{N}}(t) + v_{\tilde{N}-1}(t) + g(v_{\tilde{N}-1}(t) - v_{\tilde{N}}(t)).	
\end{align*}
The input $u$ corresponds to a current source. As detailed, e.g., in \cite{art:qlmor-gu-2011}, \cite{inproc:stahl}, the system can be recast as a quadratic-bilinear system of size $N = 2 \tilde{N}=1000$ in the new variables $x_1 = v_1$, and $x_i = v_{i-1}- v_i$ for $2\leq i \leq \tilde{N}$, and $x_{i} = \exp{(40 x_{i-\tilde{N}}})-1$ for $\tilde{N}+1 \leq i \leq 2\tilde{N}$. The output is chosen as $y = x_1$. The benchmark is treated with trivial initial conditions and two different cases of inputs:%
\begin{itemize}
\item[\emph{Case 1}] \emph{Exponential pulse}.
$u(t) = \exp{(-t)}$ with corresponding signal generator
\begin{align*}
	u = \sigstatescal, \qquad \dot{\sigstatescal} = -\sigstatescal, \qquad \sigstatescal(0) = 1.
\end{align*}
\item[\emph{Case 2}] \emph{Oscillation}.
$u(t) = 1+ \cos{(10 \pi t)}$ with corresponding signal generator
\begin{align*}
				u = [1\,|\,0\,|1] \sigstate, \qquad \dot{\sigstate} = 
	 \begin{mymatrix}
		0 & & \\
		& & 10\pi \\
		& -10\pi &
\end{mymatrix}		
	 \sigstate,
	  \qquad \,\sigstate(0) = \begin{mymatrix}
	1\\
	0 \\
	1
\end{mymatrix}.
\end{align*}
\end{itemize}
%
Case~1 is directly taken from \cite{morBre13}, \cite{techrep:two-sided-mm}, \cite{art:two-sided-hermm}, whereas Case~2 is modified from the reference to obtain a higher amplitude and frequency.

\begin{table}[t]
\renewcommand{\arraystretch}{1.25}
\centering\small{
\begin{tabular}{l|l|c}
\multirow{2}{*}{Expansion frequencies}
& \textit{AssM}{,} \textit{MultM} & $s_1=1.0$\\
& \textit{MpMo}  & $5$ IRKA-points (+1 only for $\Transfer{G}_1$) \\
\hline
\multirow{2}{*}{Order moments}
& \textit{AssM} & $\tilde{L}=3$, \, \,$L=2$  \\
& \textit{MultM} &  $q_1=5$, \, $q_2=2$ \\
\hline
Tolerance & \textit{AssM}  & $tol=6\cdot 10^{-4}$ \\
\hline
Resulting dimension & \textit{AssM}, \textit{MultM}, \textit{MpMo} & $n=11$ 
\medskip
\end{tabular}}
\caption{Reduction parameters for nonlinear RC-ladder (\textit{FOM} with $N =1000$).}
 \label{tab:redpar-RC}
\end{table}

\begin{figure}[t]
\begin{tabular}{rll|l}
\begin{minipage}{0.022\textwidth}
{\vspace{0.8cm}
{\footnotesize\rotatebox{90}{y(t)}\vspace{0.2cm}\\ \vspace{2.9cm}\\
\rotatebox{90}{Output error} \\ \vspace{2.9cm} \\
\rotatebox{90}{Output error} 
}}
\end{minipage}
&
{\hspace{-0.4cm}
\begin{minipage}{0.43\textwidth}
\center
\hspace{0.5cm} {{\underline{Case 1}}} \\
\includegraphics[height = 0.85\textwidth, width = 1.0\textwidth]{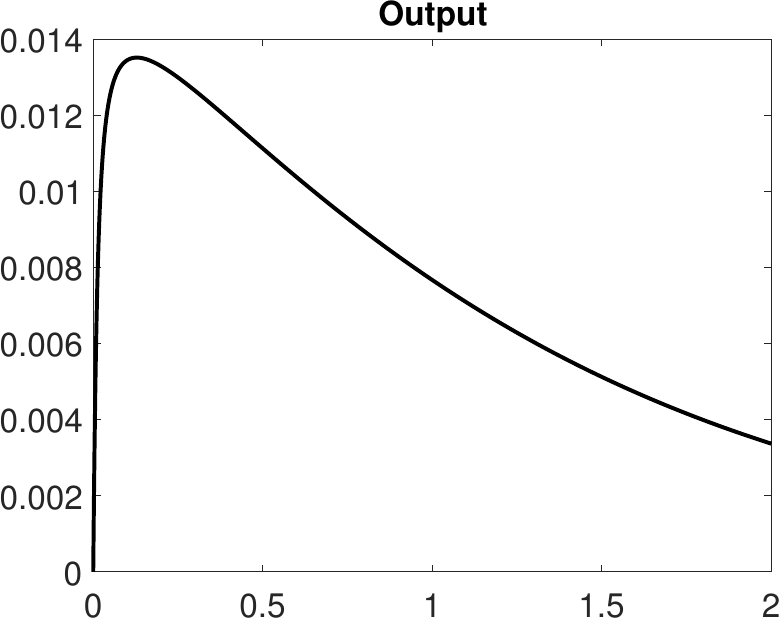} \\
\includegraphics[height = 0.85\textwidth, width = 1.0\textwidth]{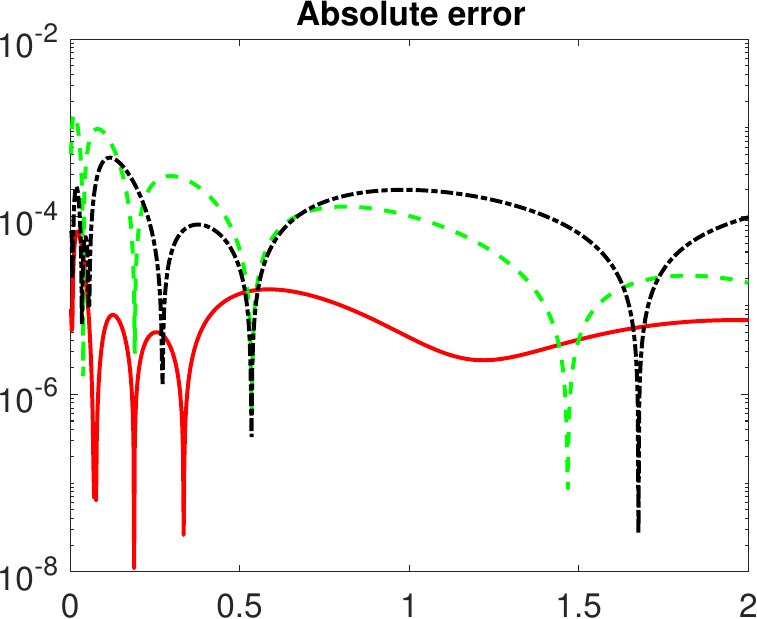} \\
\includegraphics[height = 0.85\textwidth, width = 1.0\textwidth]{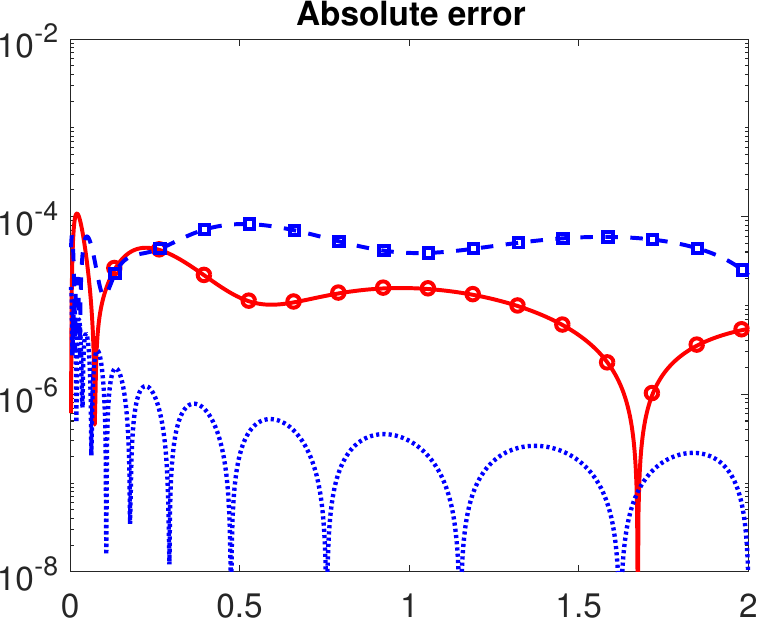} \\
 \hspace{0.5cm} {\footnotesize Time t}
\end{minipage}
}
&
\begin{minipage}{0.00\textwidth}
\end{minipage}
&
\begin{minipage}{0.43\textwidth}
\center
\hspace{0.5cm} {{\underline{Case 2}}} \\
\includegraphics[height = 0.85\textwidth, width = 1.0\textwidth]{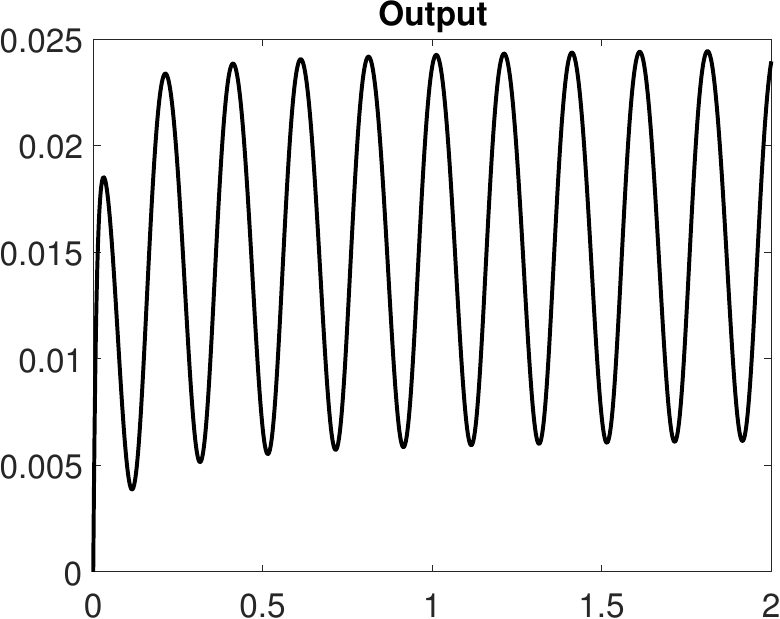} \\
\includegraphics[height = 0.85\textwidth, width = 1.0\textwidth]{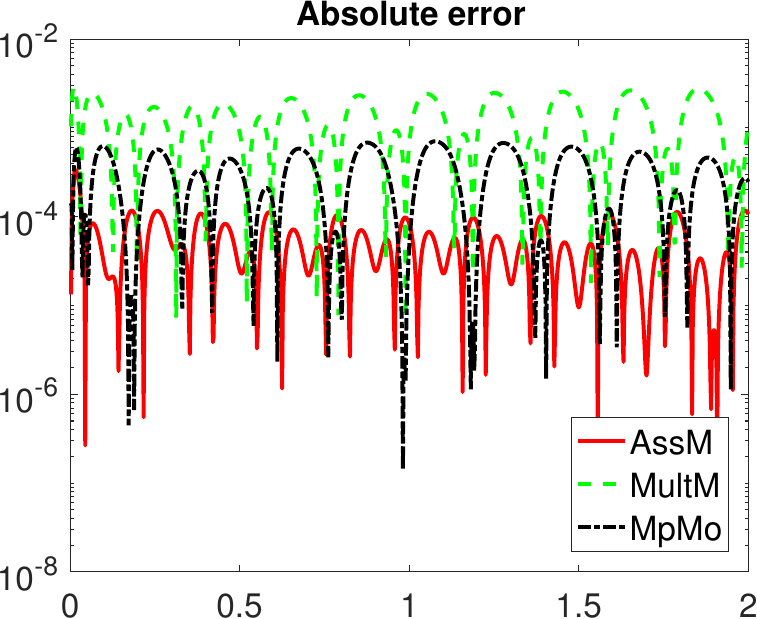} \\
\includegraphics[height = 0.85\textwidth, width = 1.0\textwidth]{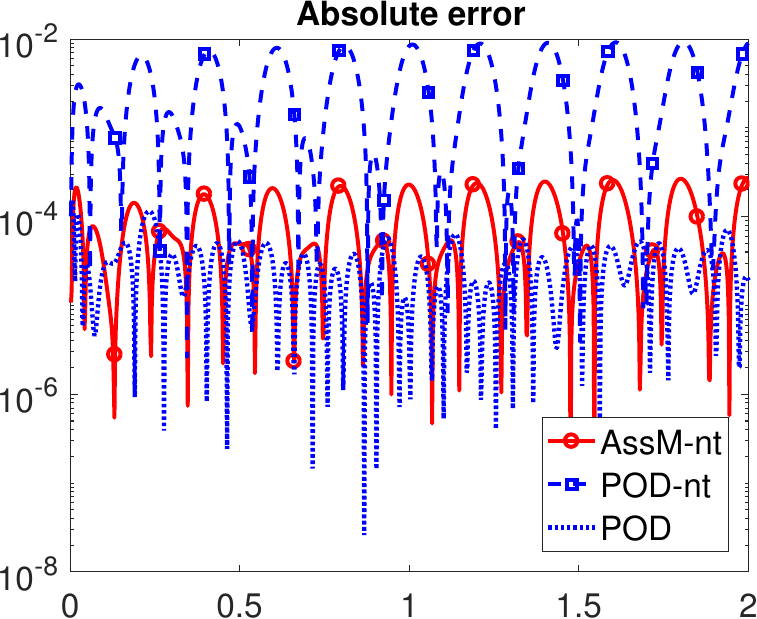} \\
 \hspace{0.5cm} {\footnotesize Time t}
\end{minipage}
\end{tabular}
\caption{Reduction results for nonlinear RC-ladder. \textit{Top to bottom:} Output $y$ (with \textit{FOM}), output errors of moment matching methods, output errors of \text{POD} and \textit{AssM} for varied training sets. Dimensions: \textit{FOM}: $N=1000$, \textit{ROM}: $n=11$ (cf.\ Table~\ref{tab:redpar-RC}).}
\label{fig:inputaware-RC-onefreq}
\end{figure}


The reduction parameters for \textit{AssM}, \textit{MultM} and \textit{MpMo} are summarized in Table~\ref{tab:redpar-RC}. Additionally, standard \textit{POD} is used to construct reduced models of same size, here $n=11$. As the inputs of the two cases differ nonlinearly from each other, both, \textit{POD} and \textit{AssM} lead to different reduced models depending on the case used in the reduction step. Results for which the training / input-tailoring scenario and the test case do not coincide are indicated by the suffix '\textit{-nt}'. We particularly take here the respective other case for training / input-tailoring.

As seen in  Fig.~\ref{fig:inputaware-RC-onefreq}, \textit{AssM} outperforms \textit{MultM} and \textit{MpMo} by about one order. \textit{POD} perfectly trained is yet superior to all three system-theoretic methods but falls off strongly when the training scenario differs from the test case. In contrast, our method shows to be much less sensitive to the training scenario. We observe reduction errors up to two orders smaller than for proper orthogonal decomposition, when training and test case do not match (compare \textit{AssM-nt} and \textit{POD-nt}).

\begin{remark}[Choice of signal generator]
There is no necessity to choose the signal generator in the reduction phase of \textit{AssM} such that it generates the signal of the test case, as we do it mostly throughout this paper. Robust (input-independent) other choices for signal generators are yet an open issue to us. Let us, however, note that in our experience the impact of the chosen signal generator in \textit{AssM} is not comparably strong as the impact of training trajectories in \textit{POD}, cf.\ Fig.~\ref{fig:inputaware-RC-onefreq}.
\end{remark}

\begin{table}[t]
\renewcommand{\arraystretch}{1.5}
\centering
\begin{tabular}{l|l|c}
Expansion frequencies & \textit{AssM}&  $s_i\in \{ 1.2 , \,  8.8, \,  37.7, \, 108.2 \}$ \qquad \qquad \quad \\
					  & \textit{AssM-inf}& $s_i \in \{ 0.2,  \,   1.3, \,    5.9, \,  20.0, \,   56.1, \,  121.3  \}$\\ \hline
{Order moments}& \textit{AssM}{,} \textit{AssM-inf}  &  $\tilde{L}=1$, \, $L=1$  \\
\hline
Tolerance & \textit{AssM} & $tol=5\cdot 10^{-4}$ \\
		  & \textit{AssM-inf} & $tol=\infty$ \qquad \qquad \\
\hline
Resulting dimension & \textit{AssM}{,} \textit{AssM-inf} & $n=12$
\medskip
\end{tabular}
\caption{Reduction parameters for nonlinear RC-ladder using multiple expansion frequencies (\textit{FOM} with $N =1000$).}
 \label{tab:redpar-RCherm}
\end{table}

\begin{figure}[t]
\begin{tabular}{rll|l}
\begin{minipage}{0.022\textwidth}
{\vspace{0.8cm}
{\footnotesize \vspace{0.4cm}\rotatebox{90}{Output error}
}}
\end{minipage}
&
{\hspace{-0.4cm}
\begin{minipage}{0.43\textwidth}
\center
\hspace{0.5cm} {{\underline{Case 1}}} \\
\includegraphics[height = 0.85\textwidth, width = 1.0\textwidth]{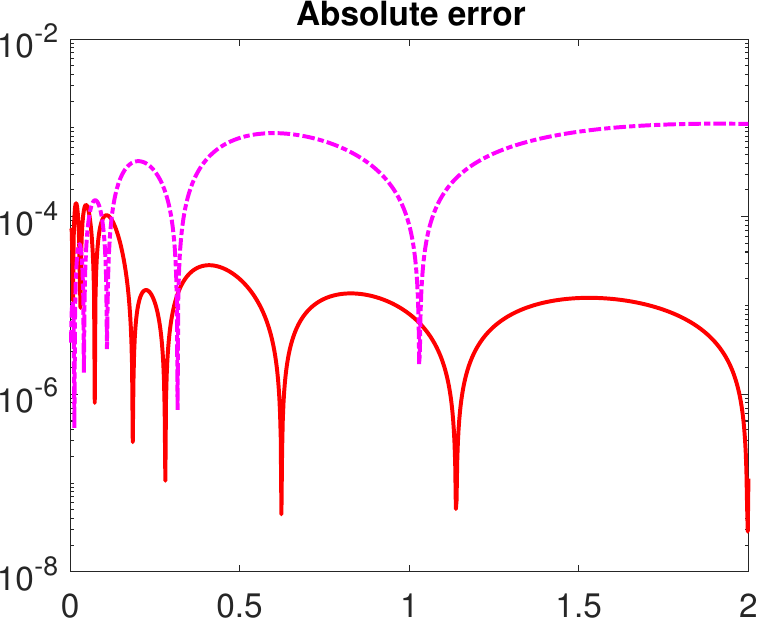} \\
 \hspace{0.5cm} {\footnotesize Time t}
\end{minipage}
}
&
\begin{minipage}{0.00\textwidth}
\end{minipage}
&
\begin{minipage}{0.43\textwidth}
\center
\quad {{\underline{Case 2}}} \\
\includegraphics[height = 0.85\textwidth, width = 1.0\textwidth]{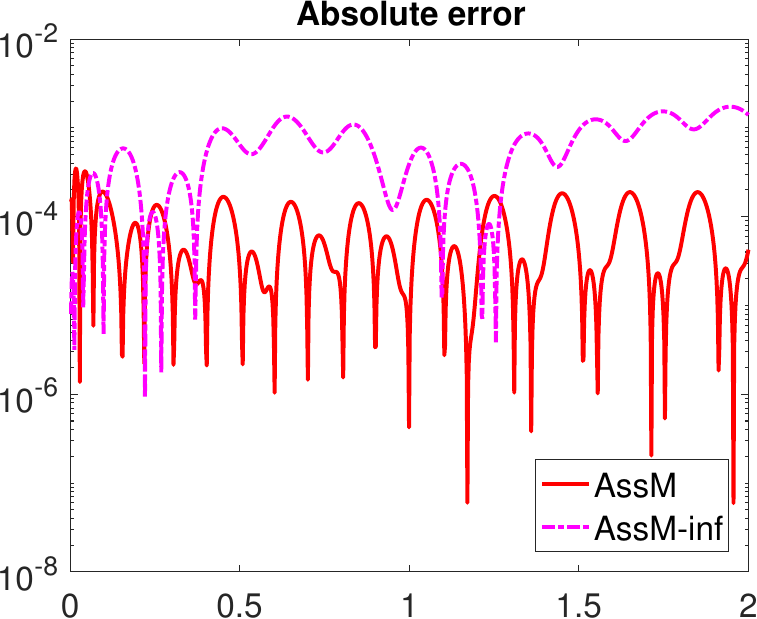} \\
 \hspace{0.5cm} {\footnotesize Time t}
 \end{minipage}
\end{tabular}
\caption{Reduction errors for nonlinear RC-ladder with multiple expansion frequencies. Dimensions: \textit{FOM}: $N=1000$, \textit{ROM}: $n=12$ (cf.\ Table~\ref{tab:redpar-RCherm}).}
\label{fig:inputaware-RC-multipfreq}
\end{figure}

In a last test, we showcase the importance of the newly proposed approximation condition \eqref{eq:tf2-approxspace-b}. This condition has no analogue in the former method \cite{inproc:fast-nonlin-mor-assoc-trafo}, \cite{art:mor-associated-transform-2016}. We repeat the test cases for the RC-ladder using our approach with altered parameters involving multiple expansion frequencies as described in Table~\ref{tab:redpar-RCherm}. The expansion frequencies are, as stated, found by applying IRKA onto the first transfer function of the Volterra series.  \textit{AssM} aims for the approximation condition \eqref{eq:tf2-approxspace-b} up to a small tolerance, whereas in \textit{AssM-inf} the approximation condition is ignored and instead more expansion frequencies are used. The latter therefore relates to \cite{inproc:fast-nonlin-mor-assoc-trafo}, \cite{art:mor-associated-transform-2016}. Although both models are of equal size, here $n=12$, \textit{AssM} leads to profoundly better results, as seen in Fig.~\ref{fig:inputaware-RC-multipfreq}.

\begin{remark}[Performance]
Our method \textit{AssM} yields low order high fidelity models that are competitive and overall similar to other system-theoretic model reduction methods as the multi-moment matching. It naturally extends to systems with non-standard input maps, which makes it in this respect similarly flexible as the trajectory-based methods like proper orthogonal decomposition. In contrast to trajectory-based methods, \textit{AssM} does not rely on pre-calculated full order model simulations. 

Concerning the offline phase the main part consists of solving Lyapunov-type equations, which makes \textit{AssM} more costly than simple multi-moment matching. Nonetheless,  \textit{AssM} is a profound enhancement over the approach in \cite{inproc:fast-nonlin-mor-assoc-trafo}, \cite{art:mor-associated-transform-2016} as it takes full advantage of the inherent tensor structure and the efficient Lyapunov solver from the \texttt{M.E.S.S.} Toolbox \cite{SaaKB16-mmess-1.0.1}. Its offline times are moderate. They are below six seconds for all tested scenarios, except for the Burgers' equation (Section~\ref{sec:burgers-eq}). There \textit{AssM} needs almost 40 seconds. However, note that the Burgers' equation is considered in a convection dominated regime which yields systems that are generally hard to reduce. For a more extensive study on computational times of \textit{AssM} in comparison to \cite{inproc:fast-nonlin-mor-assoc-trafo}, \cite{art:mor-associated-transform-2016} we refer to \cite{inproc:stahl}. 
\end{remark}

\section{Discussion and conclusion}

In this paper we suggested a new system-theoretic model reduction approach for quadratic-bilinear dynamical systems, which is based on a different perspective than the multivariate frequency-based ones. Instead of relying on input-output modeling, we used the notion of signal generator driven systems. By that input-tailored variational expansions were constructed for a large class of inputs. We compared our approach to the \textit{system-theoretic} multi-moment matching and the \textit{trajectory-based} proper orthogonal decomposition and observed competitive performance. Compared to the method in \cite{inproc:fast-nonlin-mor-assoc-trafo}, \cite{art:mor-associated-transform-2016}, which also utilizes univariate frequency representations, our method shows profound enhancements regarding analytical results and numerical performance. We stress that, in contrast to existing system-theoretic reduction methods, our method naturally extends to systems with non-trivial initial conditions and to systems with non-standard input dependencies such as, e.g., quadratic terms and time derivatives. As a byproduct of the latter, we proposed a modification for the input-output based system-theoretic methods that enables the handling of non-standard input dependencies.

We restricted the discussion in the main part to variational expansion terms up to order two. Nonetheless, the results are presented in a tensor notation allowing for convenient generalizations to higher order, as provided for the third order terms in Appendix~\ref{app:lem:mom-assoc-tf3}. Regarding higher order terms in the numerical implementation, the typical adjustments for dealing with tensors of order higher than two must of course be integrated, cf., \cite{art:kressner-kryl-tensor}, \cite{book:tensor-num-methods-quantum}. Other possible extensions of our approach could be more sophisticated automated choices of expansion frequencies including information about higher order frequencies representations by, e.g., greedy-type searches in the frequency domain, or generic (input-independent) signal generators as well as the handling of systems with more general nonlinearities. The development of a two-sided version of our input-tailored approach by the use of Petrov-Galerkin projection in the associated Lyapunov-problems \cite{book:hackbusch-tensor} could be another direction for future research.

\begin{appendix}
\renewcommand{\theequation}{\Alph{section}.\arabic{equation}}

\section{Univariate frequency representation of third order}\label{app:lem:mom-assoc-tf3}
Our approach uses univariate frequency representations tailored towards user-pre-defined families of inputs. This appendix provides the expressions and results associated to the third order term $\breve{\Transfer{W}}_3$. In particular, we state the respective extensions of Lemma~\ref{lem:assoc-trafo-repres} and Theorem~\ref{lem:mom-assoc-tf2}. The cascade- and tensor-structured pattern, which the second order terms and their moments evidently have, is preserved.

\begin{lemma}[Counterpart of Lemma~\ref{lem:assoc-trafo-repres}]\label{lem:assoc-trafo-repres_A}
Assume that the requirements of Theorem~\ref{theor:Pxu_quadratic} hold true. Then the associated frequency representation $\breve{\Transfer{W}}_3$ can also be formulated with the linear representation
	\begin{align*}
\breve{\Transfer{W}}_3(s) &= \breve{\tensorcj{C}}_3 \left(s \breve{\tensorcj{E}}_3 - \breve{\tensorcj{A}}_3 \right)^{-1} \breve{\veccj{b}}_3,\\
			 & \text{with} \quad
			\breve{\tensorcj{E}}_3 =
			\begin{mymatrix}
				\tensorcj{E} &  \\
				  		& \tensorcj{E}\supertens{2} & \\
				  	&		&				\tensorcj{E}{\supertens{3}} 
			\end{mymatrix}	
			,\quad	
				  		\breve{\tensorcj{A}}_3 =
			\begin{mymatrix}
				\tensorcj{A} & 2 \tensorcj{G}  &   \\
				  	&  	\circledmy{2}_{\tensorcj{E}} \tensorcj{A} &  \tensorcj{G} \otimes \tensorcj{E} \\
				  	&		&	\circledmy{3}_{\tensorcj{E}} \tensorcj{A} 
			\end{mymatrix} \\
			& \phantom{\text{with}} \quad	\breve{\veccj{b}}_3 = 
			\begin{mymatrix}
				\tensorbj{0} \\
				\tensorbj{0} \\
 	\veccj{b}{\supertens{3}}
			\end{mymatrix}, \quad 
						\breve{\tensorcj{C}}_3 = 
			\begin{mymatrix}
				\tensorbj{I}_M & \tensorbj{0} & \tensorbj{0}
			\end{mymatrix}.
	\end{align*}
\end{lemma} 
The linear state representation follows by straightforward calculus from Theorem~\ref{theor:Pxu_quadratic}. 

\begin{theorem}[Counterpart of Theorem~\ref{lem:mom-assoc-tf2}] \label{lem:mom-assoc-tf3}
	Assume the requirements of Theorem \ref{theor:Pxu_quadratic} and Lemma \ref{lem:assoc-trafo-repres_A} hold, and let for given $s_0 \in \mathbb{C}$ the matrix $\tensorcj{A}_{s_0} = - s_0 \tensorcj{E} + \tensorcj{A}$ be nonsingular. Then the moments $\moma^{(3)}_i$ of $\breve{\Transfer{W}}_3$ at $s_0$ are characterized by the recursion formula:
		\begin{eqnarray*}
		i = 0:\,\, & 		\circledmy{3}_{\tensorcj{E}} \tensorcj{A}_{s_0/3} \, \momc^{(3)}_0 &=\, -\veccj{b}{\supertens{3}} \\
		&	\circledmy{2}_{\tensorcj{E}} \tensorcj{A}_{s_0/2} \, \momb^{(3)}_0 &=\,  -(\tensorcj{G} \otimes \tensorcj{E})  \momc^{(3)}_{0} \\
		& \quad \tensorcj{A}_{s_0}  \moma^{(3)}_0  &=\,   -2 \tensorcj{G}  \momb^{(3)}_{0}\\
		i >0:\,\, & \circledmy{3}_{\tensorcj{E}} \tensorcj{A}_{s_0/3} \, \momc^{(3)}_i &=\,  \tensorcj{E}\supertens{3}\, \momc^{(3)}_{i-1}\\
		& \circledmy{2}_{\tensorcj{E}} \tensorcj{A}_{s_0/2} \, \momb^{(3)}_i &=\,  \tensorcj{E}\supertens{2}\, \momb^{(3)}_{i-1}  -(\tensorcj{G} \otimes \tensorcj{E})  \momc^{(3)}_{i}\\
		&	\quad \tensorcj{A}_{s_0} \, \moma^{(3)}_i &=\,  \tensorcj{E}\, \moma^{(3)}_{i-1}  -2 \tensorcj{G}  \momb^{(3)}_{i}.
	\end{eqnarray*}
Moreover, $\tensorbj{k}_i^{(3)}= [\moma^{(3)}_i;\momb^{(3)}_i; \momc^{(3)}_i]$ are the moments of $s \mapsto\left(s \breve{\tensorcj{E}}_3 - \breve{\tensorcj{A}}_3 \right)^{-1} \breve{\veccj{b}}_3$ at $s_0$.
\end{theorem}	
The proof follows similarly as the one of Theorem~\ref{lem:mom-assoc-tf2}.

\section{Proof of Theorem \ref{theor:Pxu_quadratic}} \label{app:proof-theorem-var-input}
This appendix provides the proof of Theorem~\ref{theor:Pxu_quadratic}.
The variational expansion w.r.t.\ the initial conditions in the theorem is particularly based on the following well-known result (Theorem~\ref{theor:Pxu_standard}) for that we state a proof for completeness.

\begin{theorem} \label{theor:Pxu_standard}
Consider the $\alpha$-dependent differential equation
\begin{align*}
	\dot{\veccj{w}}(t;\alpha) &= \veccj{f}(t,\veccj{w}(t;\alpha)) \qquad t \in (0,T)  \nonumber \\
		 \veccj{w}(0;\alpha) &= \hat{\veccj{b}} + \alpha  \veccj{b}  , \qquad \text{with } \hat{\veccj{b}} ,{\veccj{b}} \in \mathbb{R}^{M}
\end{align*}
for $T >0$ and a function $\veccj{f}$ being $N+1$-times continuously differentiable in $\veccj{w}$ and continuous in $t$. For $\alpha \in I$, where $I \subset \mathbb{R}$ is a bounded interval containing zero, the family of $\alpha$-dependent solutions $\veccj{w}(\cdot,\alpha)$ can be expanded as
 \begin{align*}
 	\veccj{w}(t;\alpha) = \veccj{w}_0(t) + \sum_{i=1}^N \alpha^i \veccj{w}_i(t) + \text{O}(\alpha^{N+1}), \qquad t \in [0,T).
 \end{align*}
\end{theorem} 
 \begin{proof}
With the regularity assumptions on the right hand side $\veccj{f}$, unique solutions are given by the Picard-Lindel\"of Theorem for all $\alpha \in I$. Moreover, the solution $\veccj{w}$ is $N+1$-times continuously differentiable in $\alpha$. For both statements we refer to, e.g., \cite[Sec.~5.4]{book:hartman-ode}, \cite[Sec.~1]{book:chicone2006ordinary}. Therefore, a Taylor series in $\alpha$ around $\alpha= 0$ gives
 \begin{align*}
 	\veccj{w}(t;\alpha) &= \veccj{w}_0(t) + \sum_{i=1}^N \alpha^i \veccj{w}_i(t) + \text{O}(\alpha^{N+1})\\
 	& \text{ with } \veccj{w}_i(t) := \frac{1}{i!} \frac{\partial^i}{\partial \alpha^i} \veccj{w}(t;\alpha)_{|\alpha=0}.
 \end{align*}
\end{proof}
Furthermore, we use the following technical result from \cite{art:H2mor-bilinear-breiten2012}, \cite{morBre13}.

\begin{lemma}\label{lem:kron-M-block}
	Let $\tensorbj{P}, \tensorbj{A} \in \mathbb{R}^{M,M}$, $\tensorbj{B} \in \mathbb{R}^{M,K}$, $\tensorbj{C} \in \mathbb{R}^{K,M}$, $\tensorbj{D} \in \mathbb{R}^{K,K}$, and let
	\begin{align*}
		\tensorbj{M} = 
		\begin{mymatrix}
			\tensorbj{I}_M \otimes 
			\begin{mymatrix}
				\tensorbj{I}_M \\
				\tensorbj{0}_{K,M}
			\end{mymatrix} &
			\tensorbj{I}_M \otimes
			\begin{mymatrix}
				\tensorbj{0}_{M,K} \\
				\tensorbj{I}_K				
			\end{mymatrix}
		\end{mymatrix}.
	\end{align*}
	Then it holds
	\begin{align*}
		\tensorbj{M}^T  \left( \tensorbj{P} \otimes 
		\begin{mymatrix}
			\tensorbj{A} & \tensorbj{B} \\
			\tensorbj{C} & \tensorbj{D}		
		\end{mymatrix} \right)
		 \tensorbj{M} =
		 	\begin{mymatrix}
		\tensorbj{P}\otimes \tensorbj{A} & \tensorbj{P}\otimes\tensorbj{B} \\
		\tensorbj{P}\otimes\tensorbj{C} & \tensorbj{P}\otimes\tensorbj{D}		
	\end{mymatrix}.
	\end{align*}
	Moreover, $\tensorbj{M}$ is a permutation matrix and therefore orthogonal, i.e., $\tensorbj{M}^{-1} = \tensorbj{M}^T$.
\end{lemma}
Let us now turn to the proof of Theorem~\ref{theor:Pxu_quadratic}.

\begin{proof}(Of Theorem \ref{theor:Pxu_quadratic})
Using Theorem \ref{theor:Pxu_standard} with $\veccj{f}(t,\veccj{w}) = \tensorcj{E}^{-1}(\tensorcj{A} \veccj{w} + \tensorcj{G} \veccj{w}\supertens{2})$ yields
 \begin{align*}
 	\veccj{w}(t;\alpha) = \sum_{i=1}^N \alpha^i \veccj{w}_i(t) + \text{O}(\alpha^{N+1}).
 \end{align*}
The term $\veccj{w}_0$ scaling with $\alpha^0$ drops out here as the solution for $\alpha = 0$ is $\veccj{w} \equiv \vecbj{0}$. Inserting this series representation into the differential equation and balancing equal powers in $\alpha$, we get
 \begin{align*}
 	\tensorcj{E} \dot{\veccj{w}}_1 &= \tensorcj{A} \veccj{w}_1, \qquad & \veccj{w}_1(0) = \veccj{b} \,\, \\
 	\tensorcj{E} \dot{\veccj{w}}_2 &= \tensorcj{A} \veccj{w}_2 +  \tensorcj{G} \veccj{w}_1\supertens{2}, \qquad &\veccj{w}_2(0) = \vecbj{0} \,\, \\
 	\tensorcj{E} \dot{\veccj{w}}_3 &= \tensorcj{A} \veccj{w}_3 +  \tensorcj{G} \left( \veccj{w}_1 \otimes \veccj{w}_2 +\veccj{w}_2 \otimes \veccj{w}_1 \right)  , \qquad & \veccj{w}_3(0) = \vecbj{0}.
 \end{align*}
To the equation for $\veccj{w}_1$ the standard Laplace-transform is applied, see, e.g., \cite{book:antoulas2005}, which gives the unique univariate frequency representation $\breve{\Transfer{W}}_1$ of $\veccj{w}_1$. Moreover, formally rewriting the equation for $\veccj{w}_1$ with the help of a Dirac impulse as
\begin{align*}
\tensorcj{E} \dot{\veccj{w}}_1 &= \tensorcj{A} \veccj{w}_1+ \veccj{b} \delta(t), \qquad \lim_{\bar{t}\uparrow 0} \veccj{w}_1(\bar{t}) = \vecbj{0}
\end{align*}
does not change its Laplace transform. Also multivariate frequency representations of $\veccj{w}_i$, $i=2,3$ can now be constructed following the standard procedure \cite{phd:Gu-mor-nonlinear}, \cite{book:nonlinear-system-theory-rugh}, \cite{art:mor-associated-transform-2016}. To construct the desired univariate associated frequency representations $\breve{\Transfer{W}}_i$, the Associated Transform  \cite{book:nonlinear-system-theory-rugh} is applied to the respective multivariate frequency representations of $\veccj{w}_i$. This step has already been performed for exactly our set of equations (using the Dirac impulse expression in the equation for $\veccj{w}_1$) in \cite{art:mor-associated-transform-2016}, \cite{inproc:fast-nonlin-mor-assoc-trafo}, see Remark \ref{rem:formal-volterra-dirac}. Therefore, our associated frequency representations coincide with their formally derived ones, and we can reuse their results. For $\breve{\Transfer{W}}_2$, the expression \eqref{eq:assoc-tf-2} equals \cite[eq.~(20)]{art:mor-associated-transform-2016}. To derive expression \eqref{eq:assoc-tf-3} for $\breve{\Transfer{W}}_3$, the following abbreviations are useful
\begin{align*}
\breve{\tensorcj{E}}_2 =
			\begin{mymatrix}
				\tensorcj{E} & \\
				  		& \tensorcj{E}\supertens{2}
			\end{mymatrix}	
	,\quad	
				  		\breve{\tensorcj{A}}_2 =
			\begin{mymatrix}
				\tensorcj{A} & \tensorcj{G} \\
				  	&  		\circledmy{2}_{\tensorcj{E}}\tensorcj{A}
			\end{mymatrix}, \quad 
			\breve{\veccj{b}}_2 = 
			\begin{mymatrix}
				\tensorbj{0} \\
				\veccj{b}\supertens{2}
			\end{mymatrix}, \quad 
						\breve{\tensorcj{C}}_2 = 
			\begin{mymatrix}
				\tensorbj{I}_M & \tensorbj{0}
			\end{mymatrix},
\end{align*}
cf., Lemma~\ref{lem:assoc-trafo-repres}. Then
expression \cite[eq.~(23)]{art:mor-associated-transform-2016} for $\breve{\Transfer{W}}_3$ reads in our notation
	\begin{align*} 
	\breve{\Transfer{W}}_3(s) &= (s \tensorcj{E} - \tensorcj{A} )^{-1}  \tensorcj{G} \nonumber \\
						     & \left[ (\breve{\tensorcj{C}}_2 \otimes \tensorbj{I}_M ) ( s \breve{\tensorcj{E}}_2 \otimes \tensorcj{E} - ( \breve{\tensorcj{A}}_2 \otimes \tensorcj{E} + \breve{\tensorcj{E}}_2 \otimes \tensorcj{A}))^{-1} (\breve{\veccj{b}}_2 \otimes \veccj{b}) \right. \nonumber \\
    & \, + \left. (\tensorbj{I}_M \otimes \breve{\tensorcj{C}}_2  ) ( s  \tensorcj{E} \otimes \breve{\tensorcj{E}}_2  - (\tensorcj{E} \otimes  \breve{\tensorcj{A}}_2 +  \tensorcj{A} \otimes \breve{\tensorcj{E}}_2 ))^{-1} (\veccj{b} \otimes  \breve{\veccj{b}}_2)  \right].
\end{align*}
It remains to prove that this is equivalent to \eqref{eq:assoc-tf-3}. First we show that
\begin{align}
	(\tensorbj{I}_M \otimes \breve{\tensorcj{C}}_2  ) ( s  \tensorcj{E} \otimes \breve{\tensorcj{E}}_2  - (\tensorcj{E} \otimes  \breve{\tensorcj{A}}_2 +  \tensorcj{A} \otimes \breve{\tensorcj{E}}_2 ))^{-1} (\veccj{b} \otimes  \breve{\veccj{b}}_2) \nonumber \\
	= (\breve{\tensorcj{C}}_2 \otimes \tensorbj{I}_M ) ( s \breve{\tensorcj{E}}_2 \otimes \tensorcj{E} - ( \breve{\tensorcj{A}}_2 \otimes \tensorcj{E} + \breve{\tensorcj{E}}_2 \otimes \tensorcj{A}))^{-1} (\breve{\veccj{b}}_2 \otimes \veccj{b}). \label{eq:assoc-tf-3-proofeq}
\end{align}
Using the respective orthogonal permutation matrix $\tensorbj{M}$ from Lemma \ref{lem:kron-M-block}, we get
\begin{align*}
(\tensorbj{I}_M \otimes \breve{\tensorcj{C}}_2  ) ( s  \tensorcj{E} \otimes \breve{\tensorcj{E}}_2  - (\tensorcj{E} \otimes  \breve{\tensorcj{A}}_2 +  \tensorcj{A} \otimes \breve{\tensorcj{E}}_2 ))^{-1} (\veccj{b} \otimes  \breve{\veccj{b}}_2) \\
	=
	(\tensorbj{I}_M \otimes \breve{\tensorcj{C}}_2) \tensorbj{M}  ( s  \tensorbj{M}^T (\tensorcj{E} \otimes \breve{\tensorcj{E}}_2 )\tensorbj{M} - \tensorbj{M}^ T(\tensorcj{E} \otimes  \breve{\tensorcj{A}}_2 +  \tensorcj{A} \otimes \breve{\tensorcj{E}}_2 )\tensorbj{M} )^{-1} \tensorbj{M}^T(\veccj{b} \otimes  \breve{\veccj{b}}_2)
\end{align*}
Then by Lemma \ref{lem:kron-M-block} we have 
\begin{align*}
	\tensorbj{M}^T (\tensorcj{E} \otimes \breve{\tensorcj{E}}_2)\tensorbj{M} &= 
	\begin{mymatrix}
				\tensorcj{E} \otimes \tensorcj{E} & \\
				  		& \tensorcj{E} \otimes\tensorcj{E}\supertens{2}
			\end{mymatrix}
			=
			\begin{mymatrix}
				\tensorcj{E}\supertens{2} & \\
				  		& \tensorcj{E}{\supertens{3}}
			\end{mymatrix}	= \breve{\tensorcj{E}}_2 \otimes \tensorcj{E} \\
			\tensorbj{M}^T(\tensorcj{A} \otimes \breve{\tensorcj{E}}_2  + \tensorcj{E} \otimes \breve{\tensorcj{A}}_2 )\tensorbj{M}
			&=
			 \begin{mymatrix}
		\circledmy{2}_{\tensorcj{E}} \tensorcj{A} &  \tensorcj{G} \otimes \tensorcj{E} \\
				  	&			\circledmy{3}_{\tensorcj{E}} \tensorcj{A} 
 \end{mymatrix}
 = ( \breve{\tensorcj{A}}_2 \otimes \tensorcj{E} +\breve{\tensorcj{E}}_2 \otimes \tensorcj{A}).
\end{align*}
A small calculation shows
\begin{align*}
	 \tensorbj{M}^T (\veccj{b} \otimes \breve{\veccj{b}}_2) &= 
	 \begin{mymatrix}
		\tensorbj{0} \\
 		\veccj{b}{\supertens{3}}
	 \end{mymatrix} =
			\breve{\veccj{b}}_2 \otimes  \veccj{b}\\
	(\tensorbj{I}_M \otimes \breve{\tensorcj{C}}_2) \tensorbj{M} &= \breve{\tensorcj{C}}_2 \otimes \tensorbj{I}_M	 = [\tensorbj{I}_{M^2} | \tensorbj{0} ],
\end{align*}
which together gives the equality \eqref{eq:assoc-tf-3-proofeq}. We therefore have
\begin{align*} 
	\breve{\Transfer{W}}_3(s) &= 2 (s \tensorcj{E} - \tensorcj{A} )^{-1}  \tensorcj{G} [\tensorbj{I}_{M^2} | \tensorbj{0} ]
 \left[ s \breve{\tensorcj{E}}_2 \otimes \tensorcj{E} - ( \breve{\tensorcj{A}}_2 \otimes \tensorcj{E} +\breve{\tensorcj{E}}_2 \otimes \tensorcj{A})\right] ^{-1} (\breve{\veccj{b}}_2 \otimes \veccj{b}) \nonumber \\
 	&= 2 (s \tensorcj{E} - \tensorcj{A} )^{-1} \tensorcj{G}  \left(s \tensorcj{E}\supertens{2} - \circledmy{2}_\tensorcj{E} \tensorcj{A} \right)^{-1}  (\tensorcj{G} \otimes \tensorcj{E})
 \left(s \tensorcj{E}{\supertens{3}} - \circledmy{3}_\tensorcj{E} \tensorcj{A} \right)^{-1} \veccj{b}{\supertens{3}},
\end{align*}
i.e., representation \eqref{eq:assoc-tf-3}. In the last step we just used the upper-triangular structure of the matrix in the squared brackets to be inverted to factorize the term.
\end{proof}
\section{Variational expansion w.r.t. multidimensionally para\-me\-trized initial conditions} 
\label{app:var-multdim-init}
In Theorem~\ref{theor:Pxu_quadratic} the variational expansions have been introduced for $\alpha$-dependent initial conditions where $\alpha\in \mathbb{R}$ is a scalar.
 This result is generalized to a parametrization of the initial conditions in a multidimensional linear space spanned by the column span of a matrix $\tensorcj{B}_0 \in \mathbb{R}^{M,K}$ in this appendix.
 
 We consider the $\vecbj{r}$-dependent dynamical system
\begin{align*}
	\tensorcj{E} \dot{\veccj{w}}(t;\vecbj{r}) &= \tensorcj{A} \veccj{w}(t;\vecbj{r}) + \tensorcj{G} \, (\veccj{w}(t;\vecbj{r}))\supertens{2}, \qquad t \in (0,T) \\
		 	\veccj{w}(0;\vecbj{r}) &= \tensorcj{B}_0 \vecbj{r}, \quad \text{ for } \vecbj{r} \in \mathbb{R}^K
\end{align*}
with $T >0$ and system matrices $\tensorcj{E},\tensorcj{G}$ as in Theorem~\ref{theor:Pxu_quadratic}. In this generalized setting, a variational expansion can be given as
 \begin{align*}
 	\veccj{w}(t;\vecbj{r}) =  \sum_{i=1}^{N} \veccj{w}_i(t) \vecbj{r}\supertens{i}  + \textit{higher order terms}.
 \end{align*}
The Laplace transforms $\breve{\Transfer{W}}_i$ of $\veccj{w}_i$ have an analogous form as in Theorem~ \ref{theor:Pxu_quadratic}. The only difference is that $\veccj{b}$ is replaced by $\tensorcj{B}_0$ at all instances.

To see that this holds true, we note that for each concrete choice of $\vecbj{r}$ one can define $\tilde{\veccj{b}}$ such that $\tilde{\veccj{b}}  = \tensorcj{B}_0 \vecbj{r}$. Then the state equation for $\veccj{w}$ can be written as
	\begin{align*} 
	\tensorcj{E} \dot{\veccj{w}}(t) &= \tensorcj{A} \veccj{w}(t) + \tensorcj{G} (\veccj{w}(t))\supertens{2} , \qquad
	 \veccj{w}(0) = \tilde{\veccj{b}}.
\end{align*}
Applying Theorem~\ref{theor:Pxu_quadratic} the solution $\veccj{w}$ is expanded in $\alpha$ for $\tilde{\veccj{b}} = \alpha \veccj{b}$. Afterwards a re-substitution of the '$\tilde{\veccj{b}}\supertens{i}$'-terms with the relation $\tilde{\veccj{b}}\supertens{i} = (\tensorcj{B}_0 \vecbj{r})\supertens{i}=  \tensorcj{B}_0 \supertens{i} \vecbj{r}\supertens{i}$ yields the claimed expressions.


\section{Model reduction on general nonlinear formulations}\label{app:nolin-mm-mor}
This work focuses on quadratic-bilinear systems. From a theoretical point of view, a large class of problems is covered, since many nonlinear systems can be recast in this form by introducing auxiliary variables. Whether these lifted (quadratic-bilinear) formulations should be used for model reduction is a topic in itself, cf., \cite{art:goyal-polynomial2019}, \cite{art:kramer-liftingPOD}. One advantage is certainly that no complexity reduction is needed to obtain online-efficient models. The main disadvantage, however, is that the lifted formulations are in general more difficult to approximate with projection-based model reduction, regardless of the method used. To illustrate this point (cf.\ Remark~\ref{rem-chaf-quadratize}), we reconsider the reduction of the benchmarks (Chafee-Infante equation and RC-ladder), this time with respect to their original nonlinear formulations. 

The application of \textit{POD} to nonlinear systems is straightforward. The generalization of moment matching methods to general nonlinear (not quadratic-bilinear) systems has been addressed to some extend in, e.g., \cite{art:goyal-polynomial2019}, \cite{inproc:feng2004}. We briefly explain the procedure to show its transferability to our \textit{AssM}.
The systems we consider can be written in the form
\begin{align} \label{eq:nlsys}
	\tensorbj{E} \dot{\vecbj{x}} = \vecbj{f}(\vecbj{x}) + \tensorbj{b} u, \qquad \vecbj{x}(0)= \vecbj{x}_0 = \vecbj{0}
\end{align}
with a smooth nonlinearity $\vecbj{f}$ fulfilling $\vecbj{f}(\vecbj{0}) = \vecbj{0}$ and a scalar input $u$. The reduced models are obtained by a Galerkin-projection of \eqref{eq:nlsys} with a reduction basis $\tensorbj{V}$. Let us stress that the resulting ROMs are of a general nonlinear form and not necessarily quadratic-bilinear. The basis $\tensorbj{V}$ is derived from moment matching conditions that are formulated in terms of variational expansions of the solution, in the same manner as for the quadratic-bilinear case.
As for the multi-moment matching methods (\textit{MultM}, \textit{MpMo}), they are based on a variational expansion of the solution in the input parameter $\alpha$, i.e., $\vecbj{x}= \vecbj{x}_0 + \alpha \vecbj{x}_1 + \alpha^2 \vecbj{x}_2  + \text{O}(\alpha^{3})$, for inputs $u(t) = \alpha v(t)$, cf., Section~\ref{subsec:mult-tf} and \cite{book:nonlinear-system-theory-rugh}. Using the variational expansion, the nonlinearity can be expanded in a Taylor series around $\vecbj{x}_0 = \vecbj{0}$ as 
$ \vecbj{f}(\vecbj{x}) =  \vecbj{f}'(\vecbj{0}) (\alpha \vecbj{x}_1+\alpha^2 \vecbj{x}_2)  + \alpha^2 \vecbj{f}''(\vecbj{0}) \vecbj{x}_1\supertens{2} + \text{O}(\alpha^{3})$,
where $\vecbj{f}'(\vecbj{0})$ and $\vecbj{f}''(\vecbj{0})$ denote tensor representations of the first and the second derivative of $\vecbj{f}$, respectively. Inserting the expansions in \eqref{eq:nlsys} and grouping terms of equal power in $\alpha$ gives
 \begin{align*}
 	\tensorbj{E} \dot{\vecbj{x}}_1 = \tensorbj{A} \vecbj{x}_1, + \vecbj{b}  v, \,  \quad \vecbj{x}_1(0)  = \vecbj{0}, \hspace{0.5cm} \text{and} \hspace{0.5cm} 
 	\tensorbj{E} \dot{\vecbj{x}}_2 = \tensorbj{A} \vecbj{x}_2 +  \tensorbj{G} \vecbj{x}_1\supertens{2}, \quad  \vecbj{x}_2(0) = \vecbj{0},
 \end{align*}
with $ \tensorbj{A}  =  \vecbj{f}'(\vecbj{0})$ and $\tensorbj{G} = \vecbj{f}''(\vecbj{0})$. 
From these equations for the variational terms $\vecbj{x}_1$ and $\vecbj{x}_2$, a reduction basis $\tensorbj{V}$ is constructed with the same methods as for the quadratic-bilinear case using the matrices $\tensorbj{E}$, $\tensorbj{A} $, $\tensorbj{G}$ and $\tensorbj{B}= \vecbj{b}$ and then applied in a Galerkin projection to \eqref{eq:nlsys}. 
Our moment matching approach (\textit{AssM}) is based on an input-tailored variational expansion. The derivation for the nonlinear system \eqref{eq:nlsys} is almost verbatim to the one for the quadratic-bilinear case given in Appendix~\ref{app:proof-theorem-var-input}, except for the fact  that a Taylor series expansion of the nonlinearity $\veccj{f}$ is needed at the beginning of the proof of Theorem~\ref{theor:Pxu_quadratic}. The resulting approach, which yields nonlinear ROMs satisfying an input-tailored moment matching condition, uses Algorithm~\ref{alg:mm-complete} with $\tensorbj{E}$, $ \tensorbj{A}  =  \vecbj{f}'(\vecbj{0})$, $\tensorbj{G} = \vecbj{f}''(\vecbj{0})$ and $\tensorbj{B}= \vecbj{b}$ (and otherwise the same parameters as in the quadratic-bilinear case) to determine a reduction basis $\tensorbj{V}$ and then performs a Galerkin projection of \eqref{eq:nlsys} with $\tensorbj{V}$.

\begin{figure}[tb]
\begin{tabular}{rll|l}
\begin{minipage}{0.022\textwidth}
{\vspace{0.8cm}
{\footnotesize\rotatebox{90}{}\vspace{0.2cm}\\ 
\rotatebox{90}{Output error} \\ 
}}
\end{minipage}
&
{\hspace{-0.4cm}
\begin{minipage}{0.43\textwidth}
\center
\hspace{0.5cm} {{\underline{Case 1}}} \\
\includegraphics[height = 0.85\textwidth, width = 1.0\textwidth]{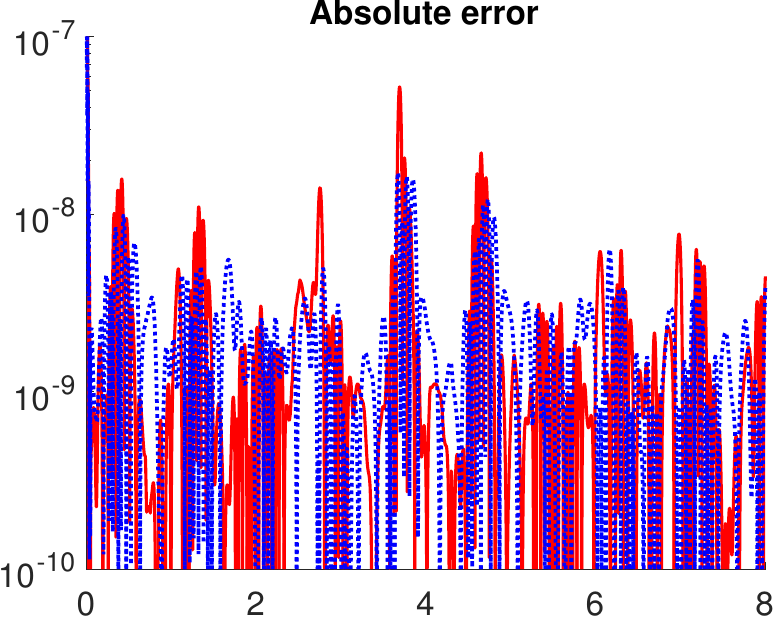} \\
 \hspace{0.5cm} {\footnotesize Time t}
\end{minipage}
}
&
\begin{minipage}{0.00\textwidth}
\end{minipage}
&
\begin{minipage}{0.43\textwidth}
\center
\hspace{0.5cm} {{\underline{Case 2}}} \\
\includegraphics[height = 0.85\textwidth, width = 1.0\textwidth]{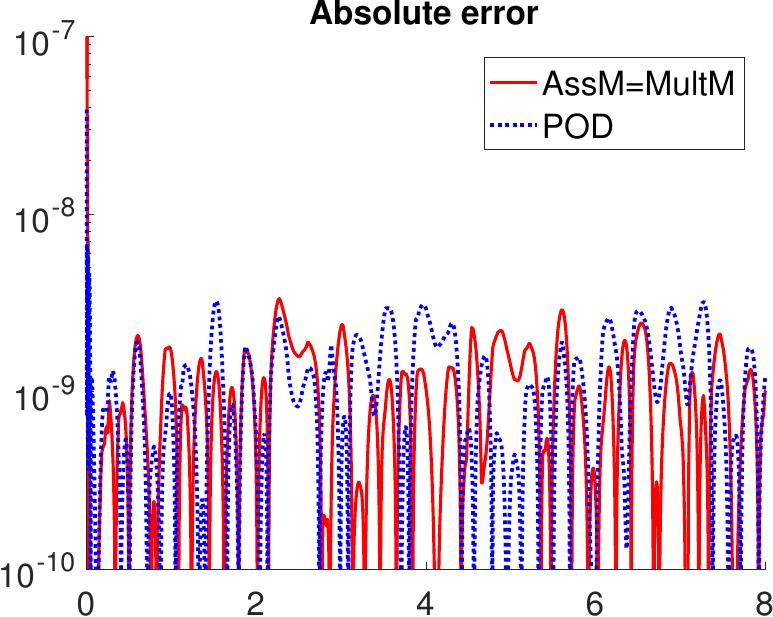} \\
 \hspace{0.5cm} {\footnotesize Time t}
\end{minipage}
\end{tabular}
\caption{Reduction errors for controlled Chafee-Infante equation in original nonlinear formulation. \textit{ROM}: $n=12$ (cf.\ results for quadratic-bilinear formulation in Fig.~\ref{fig:inputaware-chafdyn-oscil}).
\label{fig:inputaware-ChafNL}}
\vspace*{0.75cm}

\begin{tabular}{rll|l}
\begin{minipage}{0.022\textwidth}
{\vspace{0.8cm}
{\footnotesize\rotatebox{90}{}\vspace{0.2cm}\\ 
\rotatebox{90}{Output error} \\ 
}}
\end{minipage}
&
{\hspace{-0.4cm}
\begin{minipage}{0.43\textwidth}
\center
\hspace{0.5cm} {{\underline{Case 1}}} \\
\includegraphics[height = 0.85\textwidth, width = 1.0\textwidth]{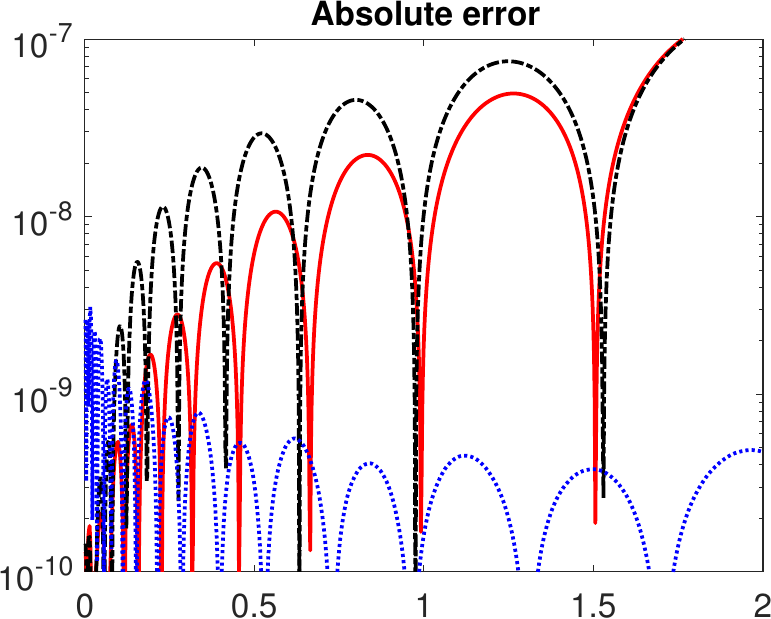} \\
 \hspace{0.5cm} {\footnotesize Time t}
\end{minipage}
}
&
\begin{minipage}{0.00\textwidth}
\end{minipage}
&
\begin{minipage}{0.43\textwidth}
\center
\hspace{0.5cm} {{\underline{Case 2}}} \\
\includegraphics[height = 0.85\textwidth, width = 1.0\textwidth]{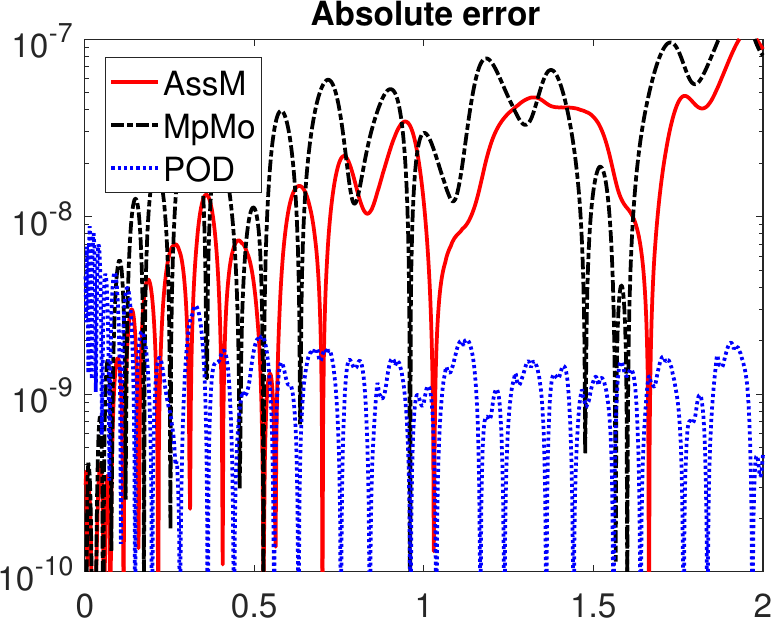} \\
 \hspace{0.5cm} {\footnotesize Time t}
\end{minipage}
\end{tabular}
\caption{Reduction errors for RC-ladder in original nonlinear formulation. \textit{ROM}: $n=10$. (cf.\ Table~\ref{tab:redpar-RCNL} and results for quadratic-bilinear formulation in Fig.~\ref{fig:inputaware-RC-onefreq}).
\label{fig:inputaware-RCNL}}
\end{figure}

In the benchmark example of the boundary-controlled Chafee-Infante equation the variational expansions have a trivial second term, as the equation has only a cubic nonlinearity and hence $\tensorbj{G}= \vecbj{f}''(\vecbj{0}) = \tensorbj{0}$ holds. Thus, the construction of the reduction basis $\tensorbj{V}$ degenerates for all moment matching methods to linear moment matching. In particular, \textit{AssM} is independent of the choice of signal generator and tolerance \textit{tol} and equals \textit{MultM} when $ \tilde{L}=q_1$ is chosen. To obtain a ROM of dimension $n=12$ as in Section~\ref{subsec:chaf-bc}, we use the same expansion frequencies as in Table~\ref{tab:redpar-chafee}, but with moment order $\tilde{L} = q_1=4$ (as $L=0$ by construction). The reduction errors of \textit{AssM} are compatible to the ones of the nonlinear \textit{POD} model of same dimension, as illustrated in Figure~\ref{fig:inputaware-ChafNL}. We stress that the reduction based on the original nonlinear formulation exceeds the fidelity of the reduction based on the quadratic-bilinear formulation by about two orders for \textit{AssM} and \textit{POD}  (cf.\ Fig.~\ref{fig:inputaware-chafdyn-oscil}).
A very similar observation can be made for the benchmark of the nonlinear RC-ladder, which allows for non-trivial variational expansions up to order two. The reduction results for \textit{AssM}, \textit{MpMo} and \textit{POD} on the nonlinear formulation show qualitatively the same trends as for the quadratic-bilinear formulation, but they are also quantitatively better for each method by about two orders of magnitude, with the nonlinear ROMs even slightly smaller ($n=10$, Fig.~\ref{fig:inputaware-RCNL}, Table~\ref{tab:redpar-RCNL}) than the respective quadratic-bilinear ones ($n=11$, Fig.~\ref{fig:inputaware-RC-onefreq}, Table~\ref{tab:redpar-RC}).

\begin{table}[bt]
\renewcommand{\arraystretch}{1.25}
\centering\small{
\begin{tabular}{l|l|c}
Expansion frequencies & \textit{AssM}, \textit{MpMo}& $5$ IRKA-points \\ \hline
{Order moments}& \textit{AssM} &  $\tilde{L}=1$, \, $L=1$  \\
\hline
Tolerance & \textit{AssM}  & $tol=5\cdot 10^{-4}$ \\
\hline
Resulting dimension & \textit{AssM}{,} \textit{MpMo} & $n=10$
\medskip
\end{tabular}}
\caption{Reduction parameters for RC-ladder in original nonlinear (not quadratic-bilinear) formulation.}
 \label{tab:redpar-RCNL}
\end{table}

\end{appendix}
\newpage
 
\bibliographystyle{alpha}
\bibliography{bib}

\end{document}